\documentclass[10pt]{amsart}
\usepackage{amsxtra, amsfonts, amsmath, amsthm, amstext, amssymb, amscd, mathrsfs, verbatim, color}
\usepackage{threeparttable}
\usepackage{mathtools}
\usepackage{tikz} 
\usepackage[ansinew]{inputenc}\usepackage[T1]{fontenc}
\usepackage[all,cmtip]{xy}
\usepackage[ngerman,french,english]{babel}

\theoremstyle{plain}

\newtheorem{theorem}{Theorem}[section]
\newtheorem{lemma}[theorem]{Lemma}
\newtheorem{definition-theorem}[theorem]{Definition-Theorem}
\newtheorem{proposition}[theorem]{Proposition}
\newtheorem{corollary}[theorem]{Corollary}

\newtheorem{conjecture}[theorem]{Conjecture}

\theoremstyle{definition}
\newtheorem{definition}[theorem]{Definition}
\newtheorem{example}[theorem]{Example}
\newtheorem{remark}[theorem]{Remark}
\newtheorem{notation}[theorem]{Notation}
\newcommand \bth[1] { \begin{theorem}\label{t#1} }
\newcommand \ble[1] { \begin{lemma}\label{l#1} }

\newcommand \bpr[1] { \begin{proposition}\label{p#1} }
\newcommand \bco[1] { \begin{corollary}\label{c#1} }
\newcommand \bde[1] { \begin{definition}\label{d#1}\rm }
\newcommand \bex[1] { \begin{example}\label{e#1}\rm }
\newcommand \bre[1] { \begin{remark}\label{r#1}\rm }

\newcommand \bnota[1] {\begin{notation}\label{n#1}\rm }
\newcommand {\ele} { \end{lemma} }

\newcommand {\epr} { \end{proposition} }
\newcommand {\eco} { \end{corollary} }
\newcommand {\ede} { \end{definition} }
\newcommand {\eex} { \end{example} }
\newcommand {\ere} { \end{remark} }
\newcommand {\enota} { \end{notation} }

\begin{document}
\setlength{\baselineskip}{1.2\baselineskip}
\title[Discrete series for type $H_{4}$]
{Discrete series for the graded Hecke algebra of type $H_{4}$}
\author[K.Y. Chan and S. Huang]{Kei Yuen Chan and Simeng Huang}

\address{The Department of Mathematics \\
The University of Hong Kong,
School of Mathematical Sciences \\
Fudan University}
\email{kychan1@hku.hk, 19110180025@fudan.edu.cn}

\selectlanguage{english} 
\begin{abstract}
This article confirms the prediction that the set of discrete series central characters for the graded (affine) Hecke algebra of type $H_4$ coincides with the set of the Heckman-Opdam central characters. Combining with previous cases of Kazhdan-Lusztig, Kriloff, Kriloff-Ram, Opdam-Solleveld, Ciubotaru-Opdam, this completes the classification of discrete series for all the graded Hecke algebras of positive parameters. Main tools include construction of calibrated modules and construction of certain minimally induced modules for discrete series. We also study the anti-sphericity, Ext-branching laws, weight structure and $W$-structure for some discrete series.

\end{abstract}
\maketitle

\setcounter{tocdepth}{1}
\tableofcontents
\part{Preliminaries}

\section{Introduction} 
\subsection{Background} 
The graded (affine) Hecke algebra is an important object in the study of representation theory of $p$-adic groups \cite{Lu89} as well as in the study of integrable systems \cite{OH97}. Discrete series is a basic building block in the representation theory of $p$-adic groups, affine Hecke algebras and their  graded version. All other irreducible representations can be built from discrete series via parabolic inductions. Thus, it is a fundamental problem to classify all discrete series. We first review some major progress of this realm in the literature:
\begin{enumerate}
    \item \cite{KL87} Kazhdan-Lusztig establishes a classification of discrete series (among other cases) for affine Hecke algebras in terms of data from nilpotent orbits. This essentially gives classification of discrete series in equal parameter cases for crystallographic cases.
    \item \cite{OH97} Heckman-Opdam determines a set of central characters expected to support discrete series.
    \item \cite{Kri99} Kriloff constructs all discrete series for type $H_3$.
    \item \cite{KR02} Kriloff-Ram constructs all discrete series for type $I_2(n)$.
    \item \cite{OS09} Opdam-Solleveld establishes an upper bound on the number of discrete series for affine Hecke algebras by analytic methods.
    \item \cite{OS10} Opdam-Solleveld classifies all discrete series central characters for graded  Hecke algebras of crystallographic types and affine Hecke algebras, including unequal positive parameter cases in non-simply laced types.
    \item \cite{Ch16} Chan establishes an upper bound on the number of discrete series to all graded Hecke algebras by algebraic methods.
    \item \cite{CO17} Ciubotaru-Opdam obtains another classification in the framework of elliptic representation theory.
\end{enumerate}

Despite several results above, the classification of discrete series for type $H_4$ was still incomplete, and this paper gives a resolution to this problem. We also compute rather concrete description of several discrete series.

In addition to obtain a parametrization set, one also wants to study some invariants of discrete series. In particular, Opdam-Solleveld \cite{OS09} and Ciubotaru-Kato-Kato \cite{CKK12} studied their formal degrees. However, the formal degree is not a well-defined invariant for non-crystallographic types, and so we shall not consider that. However, several work and classical examples e.g. \cite{OS09, Re00, Ch21, Ch24}, homological properties could be a replacement for harmonic analysis in some situations. Inspired by Prasad's proposal \cite{Pr18, Pr23}, we shall study some Ext-branching laws with the hope to use homological properties reflecting some hidden harmonic analysis.

Perhaps another interesting aspect in the representation theory for non-crystallographic types is to reveal some invisible geometry such as in the realm of Kazhdan-Lusztig theory, e.g. \cite{Al87, Cl06, EW14, Ma14}, and the Springer theory, e.g. \cite{LA82, GP00, AA08}. The $W$-module structure of discrete series is indeed expected to fit in the framework of Springer theory and to be computable from some version of Lusztig-Shoji algorithm for complex reflection groups, see \cite{KR02} and \cite{AA08}. For the viewpoint of Dirac cohomology theory \cite{BCT}, one sees \cite{Ci12, Ch13}.

After discrete series, one may ask the classification of tempered modules. In veiw of the elliptic representation theory \cite{Re01, Ci12} and computations for $H_4$ case \cite{Ch13}, we expect that other irreducible tempered representations are precisely parabolically induced from some discrete series.

\subsection{Statements of main results}

We refer the reader to Section \ref{ss gha} below for the notion of graded Hecke algebra. The Heckman-Opdam central character is defined in Definition \ref{def ho cc}, and is in terms of some combinatorial data of the root system. It is closely related to the notion of distinguished nilpotent orbits in the theory of Lusztig geometric graded Hecke algebra \cite{Lu88, Lu95}. Our first main result:

\begin{theorem} \label{thm ds equal to hocc} (=Theorem \ref{thm classification})
Let $\mathbb H$ be the graded Hecke algebra of type $H_4$ associated to a parameter function $c\neq 0$. Then, the set of discrete series central characters of $\mathbb H$ coincides with the set of Heckman-Opdam central characters of $\mathbb H$. Moreover, there are precisely 20 isomorphism classes of discrete series of $\mathbb H$.
\end{theorem}

One of our key theoretic input in proving Theorem \ref{thm ds equal to hocc} is the construction of the minimally induced module from discrete series in Theorem \ref{thm construct minimal ds}. 

The antisphericity of a module refers to the property that the module restricted to $W$ contains the sign representation of $W$. Its importance is related to the Whittaker model (see e.g. \cite{CS21} and references therein), and for graded Hecke algebras from $p$-adic groups, it is known that at each discrete series central character, there exists an anti-spherical discrete series. Since our study of discrete series has sustaintial analysis on weight spaces, it is desirable for us to reformulate the condition of antisphericity in terms of weights:

\begin{theorem} \label{thm minimal induce module} (=Theorem \ref{thm spherical anti-dominant})
Let $\mathbb H$ be an arbitrary graded Hecke algebra of a positive parameter function. Then an irreducible $\mathbb H$-module $X$ is anti-spherical if and only if $X$ has an anti-dominant weight.
\end{theorem}
Using Theorem \ref{thm minimal induce module}, we obtain new (non-regular) residual subspaces that are spherical tempered in the sense of \cite[Definition 1.6]{OH97}. 

Our proof of Theorem \ref{thm minimal induce module} uses the minimally induced modules constructed in Theorem \ref{thm construct minimal ds}. In type A case,  Theorem \ref{thm minimal induce module} is much known from the classification of irreducible representations of such algebra, and such formulation is already useful in some recent studies of Bernstein-Zelevinsky derivatives (see e.g. \cite[Proposition 3.5]{Ch25}). 

As a by-product of our study, we obtain the following interesting Ext-branching laws (see Section \ref{ss root system} for more notations):

\begin{theorem} (c.f. Theorem \ref{thm construct minimal ds}) \label{thm ext branching for discrete series}
Let $\mathbb H$ be an arbitrary graded Hecke algebra associated to a root system and a positive parameter function. Let $X$ be a discrete series of $\mathbb H$. Then there exists a parabolic subalgebra $\widetilde{\mathbb H}_I$ of $\mathbb H$ of corank $1$ and a discrete series $Y$ of $\widetilde{\mathbb H}_I$ such that 
\begin{enumerate}
\item $\mathrm{Hom}_{\widetilde{\mathbb H}_I}(X, Y)\neq 0$,
\item for $i \geq 1$, $\mathrm{Ext}^i_{\widetilde{\mathbb H}_I}(X,Y)=0$. 
\end{enumerate}
\end{theorem}

One may view Theorem \ref{thm ext branching for discrete series} as a module counterpart on the inductive definition of the Heckman-Opdam central characters (see \cite[Definition 1.4]{OH97}).




\subsection{Methods and strategy for Theorem \ref{thm ds equal to hocc}}

Note that the geometric methods in \cite{KL87} and \cite{Ka09} are not available, and the harmonic analysis method in \cite{OS09} is also only available for affine Hecke algebras. On the other hand, the algebraic methods in \cite{Kri99, KR02, Ch16} are applicable for our setting. 

We shall use $\mathbb H$ for the graded Hecke algebra of type $H_4$ (see Section \ref{ss gha}). It follows from \cite{Ch16} that the upper bound for the number of discrete series for $\mathbb H$ is 20. The remaining task is to show there exist 20 isomorphism classes of discrete series for $\mathbb H$. 

The first input is the Heckman-Opdam central characters, predicting the first places to look for discrete series. The second input is the calibrated module, developed in \cite{KR02} by Kriloff-Ram. Perhaps it is quite fortunate that after some extensive search, we can find one calibrated discrete series at each Heckman-Opdam central character (of type $H_4$). 

We now briefly explain how to search for calibrated discrete series.  A calibrated module is built from a combinatorial data, called a skew local region (see Definition \ref{def skew}), and has very nice homological properties (see Lemma \ref{lem restrict calibrated}). We use lemmas in Section \ref{ss lem checking skew} to obtain skew local regions. This is done through sections in Part \ref{part construct cal ds}. Although the proof seems to be simple, this is only the case after several trials and computations to locate those skew local regions.

Those calibrated modules account for 17 discrete series. The remaining three will be non-calibrated. Some explicit construction using local regions (and techniques of intertwining operators) on those three modules could still be possible, but rather complicated -- not even mentioning how to give a transparent proof. Nevertheless, more structure information will appear in the PhD thesis of the second-named author. 

We shall take another approach to show there are three more discrete series. It is done by an inductive construction of discrete series. Such construction is originally done in \cite{Ch14} to give an inductive proof of Ext-vanishing of discrete series. More precisely, it is shown in Theorem \ref{thm construct minimal ds}, each discrete series is a simple quotient of certain parabolically induced module from a corank $1$ graded Hecke algebra. Such parabolically induced module should usually (if not always) have simple structure in the sense that in only has one non-tempered composition factor, and its dual is a Langlands standard module. With such simple structure, we shall call those induced modules to be a minimally induced module for a discrete series.

Thus the next task is to look for those minimally induced modules. Again, this is first done by several computations before spotting the two correct candidates, which are given in Sections \ref{ss minimal induce chi 17} and \ref{ss skew 16 ds}. 

The final task is to show that the two minimally induced modules contain 3 new discrete series. To do so, we consider $X$ to be a minimally induced module.  It is expected that all discrete series appear in the cosocle of $X$. Thus, the space
\begin{align} \label{eqn space for ds}
\mathrm{Hom}_{\mathbb H}(X, X^{\star})
\end{align}
should reflect the number of discrete series, where $X^{\star}$ is the dual of $X$. One then applies Frobenius reciprocity and the geometric lemma to analyse the space (\ref{eqn space for ds}). With several tricks such as splitting of certain map for induced modules, the computation becomes quite pleasant to allow us to obtain a lower bound for the space (\ref{eqn space for ds}) in the cases that we are interested in, see Propositions \ref{prop existence ds 17} and \ref{prop existence ds 16}. Now one compares the lower bound with the expected upper bound on the number of discrete series, and obtains equalities.

\subsection{Structure of this paper}

Section \ref{sec notation prelim} reviews some basic notations, including two important sets $P(\gamma)$ and $Z(\gamma)$ in defining Heckman-Opdam central characters. Section \ref{s skew local regions} reviews local regions and construction of calibrated modules, while Section \ref{s root system}  reviews the root system of type $H_4$. Sections \ref{ss ds regular} to \ref{ss ds z>2} show the existence of calibrated discrete series at each Heckman-Opdam central character of type $H_4$, and we also compute several structures. Section \ref{s inductive construct ds} introduces minimally induced modules for discrete series, Section \ref{s splitting i} studies splitting an induced map, and Section \ref{ss geo lem} reviews the geometric lemma and the second adjointness theorem. We then use those tools to establish the existence of three more discrete series in Sections \ref{s exist ds 17} and \ref{s exist ds 16}. The classification result is given in Section \ref{s classify ds}, and Section \ref{s sph ds} proves Theorems \ref{thm minimal induce module} while Section \ref{s ext branching ds} discusses some branching laws of discrete series. In Section \ref{s springer correspond}, we shall describe a conjectural $W$-structure on discrete series and the relation to the Springer correspondence.


\subsection{Acknowledgements} The authors would like to thank Dan Ciubotaru, Eric Opdam, Gordan Savin and Peter Trapa for valuable discussions on this project over many years. The authors would also like to thank Ruben La for some communications on Iwahori-Matsumoto involutions. Several computations in this article are done in SageMath. This project is supported in part by the Research Grants Council of the Hong Kong Special Administrative Region, China (Project No: 17305223, 17308324) and the National Natural Science Foundation of China (Project No. 12322120). 

\section{Notations and Preliminaries} \label{sec notation prelim}  

\subsection{Root system} \label{ss root system}
 Let $V_0$ be a real vector space equipped with an inner product. Let $R \subset V_0$ be a \textit{(reduced) root system}, which we mean:
 \begin{enumerate}
    \item for each $\alpha \in R$, $R$ is closed under the reflection corresponding to $\alpha$;
    \item if $\alpha \in R$, then $2\alpha \notin R$.
 \end{enumerate}
 We remark that we do not assume $R$ necessarily spans $V_0$. Fix a set $\Delta$ of simple roots in $R$. Let $R^+$ be the set of positive roots determined by $\Delta$. For each $\alpha \in R$, let $s_{\alpha}$ be the corresponding reflection. Let $W$ be the reflection group generated by $s_{\alpha}$, for all $\alpha \in R$. Let $l: W \rightarrow \mathbb Z_{\geq 0}$ be the length function and let $w_o$ be the longest element in $W$. Let $V_0^{\vee}=\mathrm{Hom}_{\mathbb R}(V_0, \mathbb R)$. For each $\alpha \in R$, let $\alpha^{\vee} \in V_0^{\vee}$ such that 
\[  s_{\alpha}(v)=v-\alpha^{\vee}(v) \alpha .
\]
For $\gamma_1, \gamma_2 \in V_0$, we write $\gamma_1<\gamma_2$ if $\gamma_1\neq \gamma_2$ and $\gamma_2-\gamma_1$ can be written as a linear combination of simple roots in $\Delta$ with all non-negative coefficeints.

Let $R^{\vee}$ be the collection of all $\alpha^{\vee}$ for $\alpha \in R$. Let $V=\mathbb C\otimes_{\mathbb R}V_0$ and let $V^{\vee}=\mathbb C\otimes_{\mathbb R}V_0^{\vee}$. For $v \in V$ and $\chi \in V^{\vee}$, we also sometimes write $\langle v, \chi \rangle$ for $\chi(v)$. We shall also call elements in $V^{\vee}$ to be a {\it character} of $V$. A character $\chi \in V^{\vee}$ is said to be \textit{regular} if the stabilizer group $\mathrm{Stab}_{W}(\chi)=\{w\in W|w\chi=\chi\}$ only contains the trivial element.



The \textit{fundamental weight}  $\varpi_{\alpha}\in V_0$ associated to $\alpha \in \Delta$ is the element in $V_0$ such that $\langle \varpi_{\alpha}, \beta^{\vee} \rangle =\delta_{\alpha, \beta}$ for all $\beta \in \Delta$.

The \textit{fundamental co-weight} $\varpi_{\alpha}^{\vee} \in V_0^{\vee}$ associated to $\alpha \in \Delta$ is the element in $V_0^{\vee}$ such that $\langle \beta, \varpi_{\alpha}^{\vee}\rangle =\delta_{\alpha, \beta}$ for all $\beta \in \Delta$.


A root system $R$ is said to be \textit{crystallographic} if $\langle \alpha, \beta^{\vee} \rangle \in \mathbb{Z}$ for all $\alpha,\beta\in R$. Otherwise, we say that $R$ is of \textit{non-crystallographic} type. The irreducible noncrystallographic root systems include $I_{2}(n)$ for $n=5$ and $n\geq7$, $H_{3}$ and $H_{4}$. 



\subsection{Graded Hecke algebra} \label{ss gha}

Let $c:R\rightarrow \mathbb{C}$ be a parameter function such that $c(\alpha)=c(\alpha')$ if $\alpha$ and $\alpha'$ are in the same $W$-orbit. We write $c_{\alpha}$ for $c(\alpha)$. The parameter function $c$ is said to be {\it positive} if $c_{\alpha}>0$ for all $\alpha \in R$.

\begin{definition}\cite[Section 4]{Lu89} (also see \cite{Kri99}) \label{def graded hecke alg} The \textit{graded Hecke algebra} $\mathbb{H}=\mathbb{H}(V, R, \Delta, c)$  \footnote{The grading of the algebra originally comes from the geometry, but we shall not consider 'graded' modules and we shall not discuss the grading on $\mathbb H$.} is the associative complex algebra with an unit generated by the symbols $\{t_{w}:w\in W\}$ and $\{f_{v}:v\in V\}$ satisfying the following relations:\\
\begin{itemize}
\item[(i)] The map $w\mapsto t_{w}$ from $\mathbb{C}[W]=\bigoplus_{w\in W}\mathbb{C}w\rightarrow\mathbb{H}$ is an algebra injection.
\item[(ii)] The map $v\mapsto f_{v}$ from $S(V)\rightarrow\mathbb{H}$ is an algebra injection. For simplicity, we shall simply write $v$ for $f_{v}$ from now on.
\item[(iii)] The generators satisfy the following relation: for all $\alpha \in \Delta$ and for all $v \in V$,
$$t_{s_{\alpha}}v-s_{\alpha}(v)t_{s_{\alpha}}=c_{\alpha}\langle v, \alpha^{\vee} \rangle .$$
\end{itemize}
We shall later simply write $t_{\alpha}$ for $t_{s_{\alpha}}$.
\end{definition}

For a root system $R$, set $V_R$ to be the complex vector subspace of $V$ spanned by $R$. Then we say that a graded Hecke algebra $\mathbb H$ is \textit{associated to $R$ and of a parameter function $c$} if $\mathbb H=\mathbb H(V_R, R, \Delta, c)$.

For the remaining of this section, we fix a set of data $(V, R, \Delta, c)$ and a graded Hecke algebra $\mathbb{H}$.



\subsection{Tempered $\mathbb{H}$-modules and discrete series}

All modules we consider in this article will be over $\mathbb C$ and we shall not explicitly mention this later.  For $\chi \in V^{\vee}$, let $\mathbb{C}_{\chi}$ be the one-dimensional $S(V)$-module given by $v.x=\chi(v)x$ for $v\in V$ and $x\in \mathbb C_{\chi}$.

 Let $M$ be an $\mathbb{H}$-module and let $\chi\in V^{\vee}$. The $\chi$-\textit{weight space} and the \textit{generalized $\chi$-weight space} of $M$ are respectively 
$$M_{\chi}=\{m\in M: v.m=\chi(v)m,\forall v\in V\},$$ $$M_{\chi}^{gen}=\{m\in M: (v-\chi(v))^{k}.m=0,\forall v\in V,\exists k\in \mathbb{Z}_{\geq 0}\}.$$
Then we have $$M=\bigoplus_{\chi} M_{\chi}^{gen}$$ and we say $\chi$ is a \textit{weight} of $M$ if $M_{\chi}^{gen}\neq0$. 

\begin{definition} \label{def temp}
An irreducible $\mathbb{H}$-module $M$ is said to be \textit{tempered} if $M$ is irreducible and $$\mathrm{Re}~\chi(\varpi_{\alpha})\leq 0$$ for all weights $\chi$ of $M$ and all $\alpha \in \Delta$.
\end{definition}

\begin{definition} \label{def discrete series}
 An irreducible $\mathbb{H}$-module $M$ is said to be an \textit{essentially discrete series} if $$\mathrm{Re}~\chi(\varpi_{\alpha})<0$$ for all weights $\chi$ of $M$ and all 
 $\alpha \in \Delta$. Under the condition that $\mathrm{dim}~V=|\Delta|$, we shall simply call an essentially discrete series to be a {\it discrete series}.
\end{definition}

The terminologies in Definitions \ref{def temp} and \ref{def discrete series} are based on their correspondence to the harmonic analysis definition for $p$-adic groups (via the so-called Casselman criteria).

\subsection{Parabolic subalgebras} \label{ss parabolic subalg}

For any subset $I\subseteq \Delta$, let $W_{I}$ be the parabolic subgroup of $W$ generated by $\{s_{\alpha}: \alpha \in I\}$. Let $V_I$ be the complex vector space spanned by $I$ and let 
\[ V_I^{\bot} = \left\{ \gamma\in V: \alpha^{\vee}(\gamma)=0, \mbox{for all $\alpha \in I$} \right\} .\] 
Let $R_{I}$ be the subroot system of $R$ spanned by $I\subseteq \Delta$ i.e.
\[   R_I=R\cap V_I .
\]
Let $\mathbb{H}_{I}$ be the subalgebra of $\mathbb{H}$ generated by $t_{\alpha}$ $(\alpha \in I)$ and all $v\in V\subset S(V)$. Let $V_I^{\bot, \vee}=\mathrm{Hom}_{\mathbb C}(V_I^{\bot}, \mathbb C)$. Let $\widetilde{\mathbb{H}}_{I}$ be the subalgebra of $\mathbb H$ generated by $t_{\alpha}$ for $\alpha \in I$ and all $\alpha \in I \subset R\subset V$. Let $S(V_I^{\bot})$ be the subalgebra of $\mathbb H$ generated by $v \in V_I^{\bot}$. Since $V=V_I\oplus V_I^{\bot}$, there is a natural isomorphism $\mathbb H_I$ and $\widetilde{\mathbb H}_I\otimes  S(V_I^{\bot})$. For a $\widetilde{\mathbb H}_I$-module $\widetilde{U}$ and a character $\chi \in V_I^{\bot, \vee}$, we shall view $\widetilde{U}\otimes \mathbb C_{\chi}$ as an $\mathbb H_I$-module via the above isomorphism.

Tempered modules can be constructed from parabolic induction from discrete series in general:

\begin{proposition} \label{prop temp induced from ds}
Let $X$ be a tempered $\mathbb{H}$-module. Then there is a subset $I\subseteq \Delta$ and an $\widetilde{\mathbb H}_I$ discrete series $\widetilde{U}$ and $\chi \in V_I^{\bot, \vee}$ satisfying $\mathrm{Re}(\chi(\alpha))=0$ for all $\alpha \in \Delta-I$ such that $X$ is an irreducible subquotient module of $\mathrm{Ind}_{\mathbb{H}_I}^{\mathbb H}(\widetilde{U}\otimes \mathbb C_{\chi})$. 
\end{proposition}

The above result is well-known, see e.g. \cite[Lemma 7.4.92]{Ch14} and \cite[Lemma 5.1.1]{BC13}.

\subsection{Langlands classification}

\begin{theorem} \cite{Ev96} (also see \cite[Theorem 2.4] {KR02}) \label{thm langlands class gha}
Let $X$ be an irreducible $\mathbb{H}$-module.
\begin{enumerate}
\item[(1)]  There is a subset $I\subseteq \Delta$, a tempered  $\widetilde{\mathbb{H}}_I$-module $\widetilde{U}$ and a character $\chi \in V_I^{\bot, \vee}$ satisfying $\mathrm{Re}(\chi(\alpha))>0$ for all  $\alpha \in \Delta - I$ such that $X$ is the unique irreducible quotient module of $\mathrm{Ind}_{\mathbb{H}_I}^{\mathbb H}(\widetilde{U}\otimes \mathbb C_{\chi})$.
\item[(2)] Moreover, suppose $I'\subset \Delta$,  a tempered $\widetilde{\mathbb H}_{I'}$-module $\widetilde{U}'$ and $\chi' \in V^{\bot, \vee}$ are another data such that $\mathrm{Re}(\chi'(\alpha))>0$ for all $\alpha \in \Delta-I'$ and $X$ is unique irreducible quotient of $\mathrm{Ind}_{\mathbb H_{I'}}^{\mathbb H}(\widetilde{U}'\otimes \mathbb C_{\chi'})$. Then $I=I'$ and $\widetilde{U}\otimes \mathbb C_{\chi}\cong \widetilde{U}'\otimes \mathbb C_{\chi'}$.
\end{enumerate}
\end{theorem}

Let $U=\widetilde{U}\otimes \mathbb C_{\chi}$ above. We shall call $(U,I)$ to be the {\it Langlands parameter} of $X$, and $\mathrm{Ind}_{\mathbb H_I}^{\mathbb H}U$ is a {\it standard module}.



\subsection{Criteria for discrete series}

\begin{proposition} \label{prop for ds criteria}
Let $n=|\Delta|$. Let $\gamma_1, \ldots, \gamma_n \in R^{\vee}$ and let $c_1, \ldots, c_n \in \mathbb R$. Let $\chi=c_1\gamma_1+\ldots +c_n\gamma_n$. Suppose the following two conditions hold:
\begin{enumerate}
 \item $\gamma_1,\ldots, \gamma_n$ are linearly independent; 
 \item $c_k\gamma_k(\varpi_{\alpha})\leq0$ for all $k \in\{1,\ldots, n\}$ and all $\alpha \in \Delta$. 
 \end{enumerate}
 Then $\chi(\varpi_{\alpha})< 0$ for all $\alpha \in \Delta$. 
\end{proposition}
\begin{proof}
Since all $c_k\gamma_k(\varpi_{\alpha})\leq0$, it follows that, for all $\alpha \in \Delta$,
\[   \chi(\varpi_{\alpha})\leq 0
\]
It remains to show that the above inequalities are strict. Suppose not. Then, we have $\chi(\varpi_{\beta})=0$ for some $\beta \in \Delta$. We must then also have $c_k\gamma_k(\varpi_{\beta})=0$ for all $k$. This then implies that $\mathrm{span}\left\{ \gamma_1, \gamma_2, \ldots, \gamma_n \right\}$ is spanned by those $\alpha^{\vee}$ for $\alpha \in \Delta-\left\{ \beta \right\}$. This contradicts (1).

\end{proof}

\subsection{Central characters}

It is known (see e.g. \cite[Proposition 4.5]{Lu89}) that the center of $\mathbb H$ is isomorphic to $S(V)^W$, the set of $W$-invariant elements in $S(V)$.

For each $\chi \in V^{\vee}$ and $w \in W$, we write $(w.\chi)(v)=\chi(w^{-1}.v)$, and $W\chi$ to be the $W$-orbit of $\chi$. For a finite-dimensional $\mathbb H$-module $M$, we say that $M$ has \textit{the central character} $W\chi$ if for any $z \in S(V)^W$ and any $x\in M$, $z.x=\chi(z)x$. Note that, by Schur's lemma, an irreducible $\mathbb H$-module always has a central character.

Recall that $V_R=\mathbb C\otimes_{\mathbb R}R$.  A central character $W\chi$ is said to be a {\it discrete series central character} for $\mathbb H(V_R, R, \Delta, c)$ if there exists a discrete series of $\mathbb H(V_R, R, \Delta, c)$ with the central character $W\chi$.

\subsection{Heckman-Opdam central character}

 For $\chi \in V^{\vee}$, define 
\[ Z(\chi)=\{\alpha \in R^+:\chi(\alpha)=0\} \] 
and
\[   P(\chi)=\left\{ \alpha \in R :  \chi(\alpha)=c_{\alpha} \right\} .   \]

\begin{definition} \cite{OH97} \label{def ho cc}
A central character $W\chi \in V_R^{\vee}/W$ is said to be a {\it Heckman-Opdam central character} of $\mathbb H(V_R, R, \Delta, c)$  if $\chi$ satisfies the following condition:
\begin{align} \label{eqn HO cc}
|P(\chi)|  =2|Z(\chi)|+|\Delta| .
\end{align}
Note that (\ref{eqn HO cc}) is independent of a choice of a representative in $W\chi$.
\end{definition}

The complete list of the Heckman-Opdam central character is obtained in \cite{OH97}, in which those characters are called  residual or distinguished points. For type $H_4$, it is given in \cite[Table 4.14]{OH97}, also see \cite[Table 9]{Ch13}. We also point out a generalization in \cite{Ch13}, which fits well in the elliptic representation theory.

\section{Skew local regions} \label{s skew local regions}

We shall continue to work on general graded Hecke algebra $\mathbb H$ in this section, and fix notations $R, V, W$ in Section \ref{sec notation prelim} .

\subsection{Local region}

The inversion set of $w \in W$ is defined as
\[  R(w) = \left\{ \alpha \in R^+: w(\alpha) \in R^- \right\}.
\]

\begin{definition} \label{def local region}
\begin{enumerate}
\item Let $\chi\in V^{\vee}$. For $J\subseteq P(\chi)$, define $$\mathcal{F}^{(\chi,J)}=\{w\in W: \quad R(w)\cap Z(\chi)=\emptyset,\   w(\beta) \in R^-, \mbox{ for } \beta \in J,\  w(\beta) \in R^+  \mbox{ for } \beta \in P(\chi)\setminus J \} .$$
A \textit{local region} is a pair $(\chi,J)$ with $\chi\in V^{\vee}$ and $J\subseteq P(\chi)$ such that $\mathcal{F}^{(\chi,J)}\neq\emptyset$.
\item For $\gamma \in R$, write $\mathrm{sign}(\gamma)=+$ if $\gamma \in R^+$ and write $\mathrm{sign}(\gamma)=-$ if $\gamma \in R^-$.
\end{enumerate}
\end{definition}

The following lemma follows from definitions:

\begin{lemma} \label{lem sign characterize}
Let $(\chi, J)$ and $(\chi, J')$ be local regions with $J\neq J'$. Then,
\begin{enumerate}
    \item for any $w_1, w_2 \in \mathcal F^{(\chi, J)}$, $\mathrm{sign}(w_1(\beta))=\mathrm{sign}(w_2(\beta))$ for any $\beta\in P(\chi)$. 
    \item for $w_1 \in \mathcal F^{(\chi, J)}$ and $w_2 \in \mathcal F^{(\chi, J')}$, $\mathrm{sign}(w_1(\beta))\neq \mathrm{sign}(w_2(\beta))$ for some $\beta \in P(\chi)$.
\end{enumerate}
\end{lemma}


Lemma \ref{lem sign characterize} motivates the following sign characterization of a local region:

\begin{definition}
Let $P(\chi)=\left\{ \beta_1, \ldots, \beta_r\right\}$ with the fixed ordering. Let $(\chi, J)$ be a local region. We shall say that $\mathcal F^{(\chi, J)}$ has the {\it sign}:
\[  (\mathrm{sign}(w\beta_1), \ldots, \mathrm{sign}(w\beta_r)) ,
\]
where $w$ is any element in $\mathcal F^{(\chi, J)}$. 
\end{definition}

We record the following lemma. While we do not need this for the following proofs, it provided guidance on finding a local region:

\begin{lemma} \label{lem local region guide}
Let $\alpha \in Z(\chi)$ and let $\beta \in P(\chi)$. Then the followings hold:
\begin{enumerate}
\item $s_{\alpha}(\beta) \in P(\chi)$. 
\item Suppose $\langle \beta, \alpha^{\vee} \rangle <0$. Then, for $w \in W$ satisfying $R(w)\cap Z(\chi)=\emptyset$, if $w(s_{\alpha}(\beta))< 0$, then $w(\beta)<0$.
\end{enumerate}
\end{lemma}

\begin{proof}
The first assertion follows from $\langle s_{\alpha}(\beta),\chi\rangle =\langle \beta, s_{\alpha}(\chi) \rangle=\langle \beta,\chi\rangle$. 

For the second assertion, the condition in the lemma forces that $s_{\alpha}(\beta)=\beta+k\alpha$ for some $k>0$. We also have that $w\alpha>0$. Hence, if $w(\beta)>0$, $w(\beta)+kw(\alpha)>0$ and so $w(s_{\alpha}(\beta))>0$.
\end{proof}


\subsection{Module structure in a local region}

 We now introduce a useful tool in the study of the weight spaces of a finite-dimensional $\mathbb H$-module $M$. This is an important ingredient in constructing calibrated modules in Theorem \ref{thm construct calibrated module} below.

\begin{definition} Let $M$ be a finite dimensional $\mathbb{H}$-module. Let $\alpha \in \Delta$. Let $\chi \in V^{\vee}$ be a weight of $M$ such that $\chi(\alpha)\neq 0$. Define $\tau_{\alpha}: M_{\chi}^{gen}\rightarrow M_{s_{\alpha}\chi}^{gen}$ given by:
$$ m\mapsto(t_{\alpha}-\frac{c_{\alpha}}{\alpha})m.$$
Here $\alpha$ acts on $M_{\chi}^{gen}$ by $\chi(\alpha)$ times a unipotent transformation and so the operator $\alpha$ on $M_{\chi}^{gen}$ is invertible.
\end{definition}

A significance of the operator $\tau_{\alpha}$ comes from the following two properties:

\begin{lemma} \label{rmk local operators} \cite[Proposition 2.5]{KR02} (also see \cite{Sol12, Ch17})
Let $M$ be a finite-dimensional $\mathbb H$-module. Let $\chi \in V^{\vee}$ and $\alpha \in \Delta$ such that $\chi(\alpha)\neq 0$.
\begin{enumerate}
\item[(1)] As operators on $M^{gen}_{\chi}$, $p\tau_{\alpha}=\tau_{\alpha}s_{\alpha}(p)$ for all $p\in S(V)$;
\item[(2)] As operators on $M^{gen}_{\chi}$, $\tau_{\alpha}\tau_{\alpha}=\frac{(c_{\alpha}+\alpha)(c_{\alpha}-\alpha)}{\alpha (-\alpha)}$. Thus both maps $\tau_{\alpha}: M_{\chi}^{gen}\rightarrow M_{s_{\alpha}\chi}^{gen}$ and $\tau_{\alpha}: M_{s_{\alpha}\chi}^{gen}\rightarrow M_{\chi}^{gen}$ are invertible if and only if $\chi(\alpha)\neq\pm c_{\alpha}$.   
\end{enumerate}
\end{lemma}

\begin{lemma} \cite[Corollary 2.6]{KR02} Let $M$ be a finite dimensional $\mathbb{H}$-module. Let $\chi\in V^{\vee}$ and let $J\subseteq P(\chi)$. Then, for $w,w'\in\mathcal{F}^{(\chi,J)}$,
$$\mathrm{dim}(M_{w\chi}^{gen})=\mathrm{dim}(M_{w'\chi}^{gen}) .$$
\end{lemma}

\subsection{Skew region and $\mathbb{H}^{(\chi,J)}$}

\begin{definition} \label{def skew} For $\alpha, \beta \in \Delta$, let $R_{\alpha \beta}$ be the rank two subroot system of $R$ generated by simple roots $\alpha$ and $\beta$. A local region $(\chi,J)$ is \textit{skew} if\\
(1) for all $\gamma \in\Delta$ and $w\in\mathcal{F}^{(\chi,J)}$, $w\chi(\gamma)\neq0$;\\
(2) for all pairs of simple roots $\alpha,\beta$ and $w\in\mathcal{F}^{(\chi,J)}$ such that $\{\gamma\in R_{\alpha \beta}: w\chi(\gamma)=0\}\neq\emptyset$, the set $\{\gamma\in R_{\alpha\beta}: w\chi(\gamma)=\pm c_{\gamma}\}$ contains more than two elements.
\end{definition}

\begin{theorem} \cite[Theorem 4.5]{KR02} \label{thm construct calibrated module} Let $(\chi,J)$ be a skew local region. Define $$\mathbb{H}^{(\chi,J)}=\mathrm{span}_{\mathbb C}\{v_{w}: w\in\mathcal{F}^{(\chi,J)}\}$$ where the symbols $v_{w}$ are a labeled basis of the vector space $\mathbb{H}^{(\chi,J)}$. Then the following actions turn $\mathbb{H}^{(\chi,J)}$ into an irreducible $\mathbb{H}$-module:\\
\begin{enumerate}
\item[(1)] for each $w\in\mathcal{F}^{(\chi,J)}$ and $p \in V$, $pv_{w}=w\chi(p)v_{w}$\\
\item[(2)] for $w \in \mathcal F^{(\chi, J)}$ and $\alpha \in \Delta$, $t_{\alpha}v_{w}=\frac{c_{\alpha}}{w\chi(\alpha)}v_{w}+(1+\frac{c_{\alpha}}{w\chi(\alpha)})v_{s_{\alpha}w},\forall \alpha \in\Delta$ where we set $v_{s_{\alpha}w}=0$ if $s_{\alpha}w\notin\mathcal{F}^{(\chi,J)}$.
\end{enumerate}
\end{theorem}
\begin{proof}
Since the construction is a crucial part in this article, we sketch a proof. The proof goes by checking the module relations of Definition \ref{def graded hecke alg} are compatible with the defining actions in (1) and (2) above. 

For Definition \ref{def graded hecke alg}(2), it is straightforward from the defining action in (1). For Definition \ref{def graded hecke alg}(3), it follows from a straightforward computation using Lemma \ref{rmk local operators}(1) and the relation
\[  \tau_{\alpha}.v_w =t_{\alpha}.v_w-\frac{c_{\alpha}}{w\chi(\alpha)}v_w
\]
for $\alpha \in \Delta$ and $w \in \mathcal F^{(\chi, J)}$.

For Definition \ref{def graded hecke alg}(1), the quadratic relation appeals to the classification of modules of $\mathbb H$ attached to $A_1$, or more directly by a computation using Lemma \ref{rmk local operators}(2); and the braid relation follows from the existence of corresponding modules from the classification of modules of $\mathbb H$ attached to $I_2(n)$ in \cite{KR02}. 
\end{proof}

It follows from the definition of a skew local region and the construction in Theorem \ref{thm construct calibrated module}:
\begin{corollary} \cite{KR02} \label{cor property of cal rep}
Let $M=\mathbb H^{(\chi, J)}$ for some skew local region $(\chi, J)$. Then, for any weight $\chi'$ of $M$, 
\begin{enumerate}
\item $M_{\chi'}^{gen}=M_{\chi'}$;
\item $\mathrm{dim}~M_{\chi'}=1$.
\end{enumerate}
\end{corollary}

Indeed, the condition (1) of Corollary \ref{cor property of cal rep} characterizes those modules $\mathbb H^{(\chi, J)}$ for a skew local region $(\chi, J)$ (see \cite[Theorems 4.3 and 4.7]{KR02}). We shall call those modules $\mathbb H^{(\chi, J)}$ to be \textit{calibrated} $\mathbb H$-modules.

\begin{corollary}
Let $(\chi_1, J_1)$ and let $(\chi_2, J_2)$ be two skew local regions. Then $\mathbb H^{(\chi_1, J_1)} \cong \mathbb H^{(\chi_2, J_2)}$ if and only if there exists $w \in W$ such that $w\chi_1=\chi_2$ and $w(J_1)=J_2$.
\end{corollary}

\begin{proof}
For the if direction, one can identify the weight spaces of $\mathbb H^{(\chi_1, J_1)}$ and $\mathbb H^{(\chi_2, J_2)}$, and then it follows from the construction in Theorem \ref{thm construct calibrated module} to obtain a natural isomorphism. For the only if direction, if such $w \in W$ does not exists, then the weight spaces of $\mathbb H^{(\chi_1, J_1)}$ and $\mathbb H^{(\chi_2, J_2)}$ are not the same.
\end{proof}




\subsection{Restriction of a calibrated module} 

The following is one important homological property for a calibrated module, which will be used few times in later sections.

\begin{lemma} \label{lem restrict calibrated}
Let $M$ be a calibrated $\mathbb H$-module. Let $I \subset \Delta$. Then $M|_{\mathbb H_I}$ is semisimple and is a direct sum of irreducible calibrated modules.
\end{lemma}

\begin{proof}
Let $\chi$ be a weight of $M$. We pick a non-zero weight vector $v$ in $M_{\chi}$. Then it suffices to show that the space generated by $v$ is irreducible. Now, one applies the operators $\tau_{\alpha}$ with $\alpha \in I$ (which is well-defined by Definition \ref{def skew} (1)) to obtain a module $M'$ i.e. $M'$ is spanned by $\tau_{\alpha_r}\ldots \tau_{\alpha_1}.v$ with $\alpha_1, \ldots, \alpha_r\in I$, whenever the expression $\tau_{\alpha_r}\ldots \tau_{\alpha_1}.v$ is well-defined.

It is clear from construction that the module is invariant under $\mathbb H_I$-action. It remains to show that $M'$ is irreducible as an $\mathbb H_I$-module. To this end, we consider any non-zero weight vector $v'$ in $M'$, and then $v'$ takes the form $\tau_{\alpha_r}\ldots \tau_{\alpha_1}.v$ such that for each $i=1,\ldots, r$,
\[   \langle   \alpha_i, s_{\alpha_{i-1}}\ldots s_{\alpha_1}(\chi) \rangle\neq 0 
\]
and
\begin{align} \label{eqn condition local region}
\langle  \alpha_i, s_{\alpha_{i-1}}\ldots s_{\alpha_1}(\chi)\rangle \neq \pm c_{\alpha_i}
\end{align}
Here the first condition is guaranteed by Definition \ref{def skew}(1) and the second condition is guaranteed by the defining action (2) in Theorem \ref{thm construct calibrated module}.

By using (\ref{eqn condition local region}), $\tau_{\alpha_1}\ldots \tau_{\alpha_r}.v'$ is a non-zero scalar multiple of $v$.  This implies that any non-zero irreducible submodule of $M'$ contains $v$, and so is equal to $M'$. In other words, $M'$ is irreducible.
\end{proof}

\subsection{Lemmas for checking a local region to be skew} \label{ss lem checking skew}

\begin{lemma} \label{lem rule out orthogonal condition}
    Let $\beta,\beta',\alpha$ be three different roots in $R$. If $\beta'=s_{\alpha}(\beta)=\beta-\langle \beta,\alpha^{\vee} \rangle\alpha$ and $\mathrm{sign}(\beta')=-\mathrm{sign}(\beta)$, then $\alpha$ is not a simple root.
\end{lemma}
\begin{proof}
If $\alpha$ is a simple root, then $s_{\alpha}$ only changes the signs of $\pm\alpha$. But $\beta\neq\pm\alpha$, thus $\mathrm{sign}(\beta')=\mathrm{sign}(\beta)$, which is a contradiction. 
\end{proof}

\begin{lemma} \label{lem roots with diff sign}
   Let $R'$ be a rank $2$ subroot system of $R$. Let $\beta,\beta',\alpha$ are three different roots in $R$. If $\beta'=s_{\alpha}(\beta)$, $\mathrm{sign}(\beta')=-\mathrm{sign}(\beta)$ and $\alpha\in R'$, then $\beta,\beta'\in R'$. 
\end{lemma}
\begin{proof}
Without loss of generality, we may assume $\beta$ is in $R^+$. Suppose $\beta \notin R'$. When we write $\beta=\sum_{\gamma \in \Delta}c_{\gamma}\alpha$ into a linear combination of simple roots, $c_{\gamma'}>0$ for some $\gamma' \notin R'$. But then, since $\alpha \in R'$, writing $s_{\alpha}(\beta)$ into the linear combination of  simple roots also has a positive coefficient for $\gamma'$. This contradicts that $\beta' \in R^-$.
\end{proof}

\begin{lemma} \label{lem no rank 2 subroot system}
Let $\beta, \beta', \alpha \in R$. Suppose the following two conditions holds:
\begin{enumerate}
\item $\beta, \beta', \alpha$ are linearly independent; and
\item $\mathrm{sign}(\beta)=-\mathrm{sign}(s_{\alpha}(\beta))$ and $\mathrm{sign}(\beta')=-\mathrm{sign}(s_{\alpha}(\beta'))$.
\end{enumerate}
Then $\alpha$ is not in any rank $2$ subroot system of $R$.
\end{lemma}

\begin{proof}
Suppose $\alpha$ is in some rank 2 subroot system $R'$ of $R$. Then Lemma \ref{lem roots with diff sign} implies that $\beta, \beta'$ are also in $R'$. This then contradicts that $\beta, \beta', \alpha$ are linearly independent.
\end{proof}

\section{Root system of type $H_4$} \label{s root system}


\subsection{Root system of type $H_{4}$} We remark that the convention for the root system of type $H_4$ follows \cite[Section 2.13]{Hu90} and also agrees with the convention in SageMath for the command {\small \textup{CoxeterGroup(['H',4])}}, which differs from some other sources such as \cite{Gr74} and \cite[Section 11.2]{GP00}. Let $V_0=V_0^{\vee}=\mathbb{R}^{4}$ and the non-degenerate pairing $\langle , \rangle$ be the standard inner product on $\mathbb{R}^4$. Let $R=R^{\vee}$ consist of vectors obtained from $\sqrt{2}(1,0,0,0)$, $\sqrt{2}(1/2,1/2,1/2,1/2)$ and $\sqrt{2}(a,1/2,b,0)$ by even permutations of coordinates and arbitrary sign changes with $a=(\sqrt{5}+1)/4$ and $b=(\sqrt{5}-1)/4$. One has $R\subset \mathbb R^4$ to form a root system of type $H_{4}$, and we fix a set of simple root $\Delta=$
$$\{\alpha_{1}=\sqrt{2}(-1/2,-a,0,b), \alpha_{2}=\sqrt{2}(1/2,b,-a,0), \alpha_{3}=\sqrt{2}(-a,1/2,b,0), \alpha_{4}=\sqrt{2}(a,-1/2,b,0)\}.$$
As usual, $||\alpha_1||^2=||\alpha_2||^2=||\alpha_3||^2=||\alpha_4||^2=2$. The corresponding Dynkin diagram is:

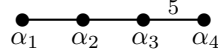
\setlength{\unitlength}{0.8cm}
\begin{picture}(12,1.3)
\thicklines
\put(9.4,0.9){\footnotesize $5$}
\put(7,0.8){\circle*{0.2}}
\put(8,0.8){\circle*{0.2}}
\put(9,0.8){\circle*{0.2}}
\put(10,0.8){\circle*{0.2}}
\put(7,0.8){\line(1,0){1}}
\put(8,0.8){\line(1,0){1}}
\put(9,0.8){\line(1,0){1}}
\put(6.8,0.4){$\alpha_1$}
\put(7.8,0.4){$\alpha_2$}
\put(8.8,0.4){$\alpha_3$}
\put(9.8,0.4){$\alpha_4$}
\end{picture}

The set of (a total of 60) positive roots in $R$ can be explicitly described as: 
\[R^{+}=\sqrt{2}\{
(a,1/2,b,0),(a,-1/2,b,0),(-a,1/2,b,0),(a,1/2,-b,0),\]\[(1/2,b,a,0),(1/2,-b,a,0),(-1/2,b,a,0),(1/2,b,-a,0),\]
\[(b,a,1/2,0),(-b,a,-1/2,0),(-b,a,1/2,0),(b,a,-1/2,0),\]\[(1,0,0,0),(0,1,0,0),(0,0,1,0),(0,0,0,1),
(\pm1/2,\pm1/2,\pm1/2,1/2),\]\[(\pm1/2,\pm a,0,b),(\pm a,\pm b,0,1/2),(\pm b,\pm1/2,0,a),\]\[
(0,\pm1/2,\pm a,b),(0,\pm a,\pm b,1/2),(0,\pm b,\pm1/2,a),\]\[(\pm a,0,\pm1/2,b),(\pm b, 0,\pm a,1/2),(\pm1/2,0,\pm b,a)\}; \]
For $\alpha \in R$ and $\beta^{\vee} \in R^{\vee}$, we define $\beta^{\vee}(\alpha)=\langle \alpha, \beta^{\vee}\rangle$. The fundamental co-weights and the fundamental weights are\[\varpi_{\alpha_1}=\sqrt{2}(0,0,0,2a)=(4a+2)\alpha_1+(8a+3)\alpha_2+(12a+4)\alpha_3+(10a+3)\alpha_4,\]
\[\varpi_{\alpha_2}=\sqrt{2}(1/2,a,0,3a+1/2)=(8a+3)\alpha_1+(16a+6)\alpha_2+(24a+8)\alpha_3+(20a+6)\alpha_4,\]
\[ \varpi_{\alpha_3}=\sqrt{2}(1/2,a+1/2,a,4a+1)=(12a+4)\alpha_1+(24a+8)\alpha_2+(36a+12)\alpha_3+(30a+9)\alpha_4,\]
\[\varpi_{\alpha_4}=\sqrt{2}(a,a,a,3a+1)=(10a+3)\alpha_1+(20a+6)\alpha_2+(30a+9)\alpha_3+(24a+8)\alpha_4. \]
We shall also write $\varpi_i$ for $\varpi_{\alpha_i}$ later. Let $s_i=s_{\alpha_i}$ be the reflection corresponding to $\alpha_i$. The following braid relations are satisfied: $$(s_1s_3)^2=(s_1s_4)^2=(s_2s_4)^2=(s_1s_2)^3=(s_2s_3)^3=(s_3s_4)^5=1 $$
and so the group generated by reflections $s_i$ forms the reflection group $W$ of type $H_4$.

\part{Construction of calibrated discrete series} \label{part construct cal ds}

In this part, we shall consider $R$ to be the root system associated to $H_4$ and $\mathbb H$ to be the graded Hecke algebra of type $H_4$. Since all the roots in $R$ are in the same $W$-orbit, the parameter function has to be a constant function. Note that two non-zero parameter functions give isomorphic graded Hecke algebras and so we may specify $c$ to be certain constant function. From now on, we shall take $c\equiv 1/2$ for the graded Hecke algebra of type $H_4$. In particular, we have, for all $\alpha \in \Delta$, $t_{s_{\alpha}}\alpha+\alpha t_{s_{\alpha}} =1$.

\section{Discrete series at regular central characters} \label{ss ds regular}


\subsection{Regular central character} 

We first recall the following construction theorem for irreducible $\mathbb H$-modules of regular central characters.

\begin{proposition} \label{prop regular character skew}
    Let $\chi$ be a regular character. Then the followings hold:
    \begin{enumerate}
    \item  each local region $(\chi, J)$ (for some $J \subset P(\chi)$) is skew. 
    \item for any irreducible $\mathbb H$-module $M$ with the central character $W\chi$,
    \[   M \cong \mathbb H^{(\chi, J)}
    \]
    for some $J \subset P(\chi)$.
    \end{enumerate}
\end{proposition}
\begin{proof}
    The first assertion follows directly from $Z(\chi)=\emptyset$, and the second assertion follows from \cite[Theorems 4.3 and 4.7]{KR02}.
\end{proof}

\begin{theorem} \label{thm regular ds}
For $j=1, \ldots, 12$, let $\chi_j \in V^{\vee}$ be the element uniquely determined by the condition $P(\chi_j)$ equal to the set given in the second column of Table \ref{table regular discrete series}.

Then there exists a skew local region $(\chi_j, J_j)$ such that $\mathcal F^{(\chi_j, J_j)}$ has the sign given by the third column of Table \ref{table regular discrete series}. Moreover, $\mathbb H^{(\chi_j, J_j)}$ is a discrete series and its dimension is equal to the corresponding value in the forth column of Table \ref{table regular discrete series}.
\\

\begin{table}
    \centering
\begin{tabular}{|c|c|c|c|c|}
\hline
$\chi_j$ & $P(\chi_j)=\{\beta_j^1,\beta_j^2,\beta_j^3,\beta_j^4\}$ & $\mathcal F^{(\chi_j, J_j)}$ & Dimension \\
\hline
$\chi_{1}$ & $\alpha_{1},\alpha_{2},\alpha_{3},\alpha_{4}$ & $(-,-,-,-)$ & $1$\\
\hline
$\chi_{2}$ & $\alpha_{1}+\alpha_{2},\alpha_{2}+\alpha_3,\alpha_{3}+2a\alpha_4,2a\alpha_3+\alpha_4$ & $(-,-,-,-)$ & $14$\\
\hline
$\chi_{3}$ & $\alpha_{1},\alpha_{2}+\alpha_3,\alpha_{3}+2a\alpha_4,2a\alpha_3+\alpha_4$ & $(-,-,-,+)$ & $4$\\
\hline
$\chi_{4}$ & $\alpha_{1},\alpha_{2},\alpha_{3}+2a\alpha_4,2a\alpha_3+\alpha_4$ & $(-,-,-,-)$ & $5$\\
\hline
$\chi_{5}$ & $\alpha_{1}+\alpha_{2}+\alpha_3,\varepsilon_1,2a\alpha_2+2a\alpha_3+\alpha_4,\varepsilon_3$ & $(-,-,-,-)$ & $55$\\
\hline
$\chi_{6}$ & $\alpha_{1}+\alpha_{2},\varepsilon_1,2a\alpha_2+2a\alpha_3+\alpha_4,\varepsilon_3$ & $(-,+,-,-)$ & $20$\\
\hline
$\chi_{7}$ & $\alpha_{1},\varepsilon_1,2a\alpha_2+2a\alpha_3+\alpha_4,\varepsilon_3$ & $(-,-,-,-)$ & $30$\\
\hline
 $\chi_{8}$ & \tiny{$\alpha_{1}+\alpha_{2}+\alpha_3+2a\alpha_4,2a\alpha_{2}+2a\alpha_3+2a\alpha_4,2a\alpha_1+2a\alpha_2+2a\alpha_3+\alpha_4,\alpha_2+(2a+1)\alpha_3+2a\alpha_4$} & $(-,-,-,-)$ & $115$\\
\hline
$\chi_{9}$ & \tiny{$2a\alpha_{1}+2a\alpha_{2}+2a\alpha_3+2a\alpha_4,\alpha_1+\alpha_{2}+(2a+1)\alpha_3+2a\alpha_4,\alpha_2+(2a+1)\alpha_3+(2a+1)\alpha_4,\varepsilon_2$} & $(-,-,-,-)$ & $240$\\
\hline
$\chi_{10}$ & \tiny{$2a\alpha_{1}+2a\alpha_{2}+2a\alpha_3+\alpha_4,\alpha_{2}+(2a+1)\alpha_3+(2a+1)\alpha_4,\varepsilon_2,\alpha_1+\alpha_2+(2a+1)\alpha_3+2a\alpha_4$} & $(-,-,-,+)$ & $86$\\
\hline
$\chi_{11}$ & $\beta_{11}^1,\beta_{11}^2,\beta_{11}^3,\beta_{11}^4$ & $(+,-,-,-)$ & $284$\\
\hline
$\chi_{12}$ & $\beta_{12}^1,\beta_{12}^2,\beta_{12}^3,\beta_{12}^4$ & $(-,-,-,-)$ & $409$\\
\hline
\end{tabular}
    \caption{Regular discrete series central characters \\
$\varepsilon_1=\sqrt{2}(1,0,0,0)=\alpha_2+\alpha_3+2a\alpha_4,$\\
    $\varepsilon_2=\sqrt{2}(0,1,0,0)=2a\alpha_2+(2a+1)\alpha_3+2a\alpha_4,$\\
    $\varepsilon_3=\sqrt{2}(0,0,1,0)=2a\alpha_3+2a\alpha_4,$\\
    $\varepsilon_4=\sqrt{2}(0,0,0,1)=4a\alpha_1+(6a+1)\alpha_2+(8a+2)\alpha_3+(6a+2)\alpha_4,$\\
    $\beta_{11}^1=2a\alpha_{1}+2a\alpha_{2}+(2a+1)\alpha_3+2a\alpha_4,$\\
    $\beta_{11}^2=\alpha_1+(2a+1)\alpha_{2}+(2a+1)\alpha_3+2a\alpha_4,$\\
    $\beta_{11}^3=\alpha_1+\alpha_2+(2a+1)\alpha_3+(2a+1)\alpha_4,$\\
    $\beta_{11}^4=2a\alpha_2+4a\alpha_3+(2a+1)\alpha_4,$\\
    $\beta_{12}^1=2a\alpha_{1}+(2a+1)\alpha_{2}+(2a+1)\alpha_3+2a\alpha_4,$\\
    $\beta_{12}^2=2a\alpha_1+2a\alpha_{2}+(2a+1)\alpha_3+(2a+1)\alpha_4,$\\
    $\beta_{12}^3=\alpha_1+(2a+1)\alpha_2+(2a+1)\alpha_3+(2a+1)\alpha_4,$\\
    $\beta_{12}^4=2a\alpha_2+4a\alpha_3+(2a+1)\alpha_4.$\\}
    \label{table regular discrete series}
\end{table}

\end{theorem}

\begin{proof}
Recall that we have identified $R$ and $R^{\vee}$. For each $\chi_j$, write $\chi_j =c_1\beta_j^1+c_2\beta_j^2+c_3\beta_j^3+c_4\beta_j^4$. It is a direct computation to check that there exists a pair $(\chi_j, J_j)$ such that $\mathcal F^{(\chi_j, J_j)}\neq \emptyset$ and $\mathcal F^{(\chi_j, J_j)}$ has the sign $(-\mathrm{sign}(c_1),-\mathrm{sign}(c_2),-\mathrm{sign}(c_3),-\mathrm{sign}(c_4))$. The linear combinations of each $\chi_j$ in terms of $\beta_j^1, \beta_j^2, \beta_j^3, \beta_j^4$ can be found in an appendix (Section \ref{s rdsc}).

Now it follows from Proposition \ref{prop for ds criteria} that $\mathbb H^{(\chi_j, J_j)}$ is a discrete series. By Theorem \ref{thm construct calibrated module}, it follows that the dimension of $\mathbb H^{(\chi_j, J_j)}$ is equal to the number of elements in $\mathcal F^{(\chi_j, J_j)}$. The latter set is computed in SageMath by using the definition of $\mathcal F^{(\chi_j, J_j)}$. 
\end{proof}





\section{Discrete series at central characters $\chi$ with $|Z(\chi)|=1$}

\subsection{Classification of discrete series with $|Z(\chi)|=1$}

\begin{theorem} \label{thm ds z=1}
For $j=13,14,15$, define $\chi_j \in V^{\vee}$ determined uniquely by the sets $Z(\chi_j)$ and $P(\chi_j)$ given respectively in the second and third columns of Table \ref{table discrete series z=1}. Then there exists $J_j\subset P(\chi_j)$ such that $(\chi_j, J_j)$ is a skew local region and the sign of $\mathcal{F}^{(\chi_j,J_j)}$ is given by the forth column of Table \ref{table discrete series z=1}. Moreover, $\mathbb H^{(\chi_j, J_j)}$ is a discrete series.

\begin{table}
    \centering
\begin{tabular}{|c|c|c|c|c|c|}
\hline
$\chi_j$ & $Z(\chi_j)$ & $P(\chi_j)$ & $\mathcal F^{(\chi_j, J_j)}$ & $Dimension$ \\
\hline
$\chi_{13}$ & $\alpha_4$ & $\beta_{13}^1,\beta_{13}^2,\beta_{13}^3,\beta_{13}^4,\beta_{13}^5,\beta_{13}^6$ & $(-,-,-,+,-,+)$ & $81$\\
\hline
$\chi_{14}$ & $\alpha_4$ & $\beta_{14}^1,\beta_{14}^2,\beta_{14}^3,\beta_{14}^4,\beta_{14}^5,\beta_{14}^6$ & $(-,+,-,+,-,-)$ & $10$\\
\hline
$\chi_{15}$ & $\alpha_2$ & $\beta_{15}^1,\beta_{15}^2,\beta_{15}^3,\beta_{15}^4,\beta_{15}^5,\beta_{15}^6$ & $(-,+,-,+,-,-)$ & $9$\\
\hline
\end{tabular}\\
 \caption{Discrete series central characters  \\
 $\beta_{13}^1=\sqrt{2}(0,1,0,0)$, $\beta_{13}^2=\sqrt{2}(a,1/2,b,0)$, \\ $\beta_{13}^3=\sqrt{2}(-a,0,1/2,b)$, 
$\beta_{13}^4=\sqrt{2}(0,-1/2,a,b)$, \\
$\beta_{13}^5=\sqrt{2}(-1/2,-1/2,-1/2,1/2)$, $\beta_{13}^6=\sqrt{2}(0,-a,-b,1/2)$; \\
$\beta_{14}^1=\sqrt{2}(-b,a,1/2,0)$, 
$\beta_{14}^2=\sqrt{2}(1/2,b,a,0)$, \\
$\beta_{14}^3=\sqrt{2}(b,a,-1/2,0)$, $ \beta_{14}^4=\sqrt{2}(a,1/2,-b,0)$, \\
$\beta_{14}^5=\sqrt{2}(-a,0,-1/2,b)$, $\beta_{14}^6=\sqrt{2}(1/2,-a,0,b)$  \\
$ \beta_{15}^1=\sqrt{2}(-1/2,-a,0,b)$, $\beta_{15}^2=\sqrt{2}(0,-1/2,-a,b)$, \\
$\beta_{15}^3=\sqrt{2}(1/2,-b,a,0)$, $\beta_{15}^4=\sqrt{2}(1,0,0,0)$,\\ $\beta_{15}^5=\sqrt{2}(-1/2,b,a,0)$, $\beta_{15}^6=\sqrt{2}(b,a,-1/2,0)$.
 }
    \label{table discrete series z=1}
\end{table}

\end{theorem}

We only show the existence of calibrated discrete series at $W\chi_{13}$ in Section \ref{ss proof thm z=1}. The argument for $W\chi_{14}$ and $W\chi_{15}$ is similar by using computations in Section \ref{ss useful formulas} and more explicit structure is given in Section \ref{ss explicit module str}.

\subsection{Some useful formulas} \label{ss useful formulas}
We shall record some useful formulas. In below, we shall give some good choices of linear combinations of those $\beta_j^k$ so that one can apply the simple trick of Lemma \ref{lem no rank 2 subroot system}. 

\begin{align} \label{eqn chi13 beta relation}
\beta_{13}^2=s_4(\beta_{13}^1)=\beta_{13}^1+\alpha_4, \quad \beta_{13}^4=s_4(\beta_{13}^3)=\beta_{13}^3+\alpha_4, \quad 
\beta_{13}^6=s_4(\beta_{13}^5)=\beta_{13}^5+2b\alpha_4,
\end{align}

\begin{align*}
\chi_{13}&=c_1\beta_{13}^1+c_2\beta_{13}^2+c_3\beta_{13}^3+c_5\beta_{13}^5 =\frac{\sqrt{2}}{6a+2}(a, \frac{3}{2}a+\frac{1}{2}, \frac{1}{2}a+\frac{1}{4}, 6a+\frac{7}{4}),
\end{align*}
where $c_1=\frac{1}{6a+2}>0, c_2=\frac{11a+3}{6a+2}>0, c_3=\frac{3a+2}{6a+2}>0, c_5=\frac{8a+4}{6a+2} >0$.


\begin{align} \label{eqn relation betas}
\beta_{14}^2=s_4(\beta_{14}^1)=\beta_{14}^1+\alpha_4, \quad \beta_{14}^4=s_4(\beta_{14}^3)=\beta_{14}^3+2b\alpha_4, \quad 
\beta_{14}^6=s_4(\beta_{14}^5)=\beta_{14}^5+2a\alpha_4,
\end{align}

\begin{align*}
\chi_{14}&=k_1\beta_{14}^1+k_3\beta_{14}^3+k_5\beta_{14}^5+k_6\beta_{14}^6=\frac{\sqrt{2}}{8a+2}(a+\frac{1}{2}, 2a+1, a, 10a+\frac{7}{2}), 
\end{align*}
where  $k_1=\frac{20a+5}{8a+2}>0, k_3=\frac{10a+1}{8a+2}>0, k_5=\frac{4a+2}{4a+1}>0, k_6=\frac{13a+3}{4a+1}>0$.

\[  \beta_{15}^2=s_2(\beta_{15}^1)=\beta_{15}^1+\alpha_2, \quad \beta_{15}^4=s_2(\beta_{15}^3)=\beta_{15}^3+\alpha_2, \quad 
\beta_{15}^6=s_2(\beta_{15}^5)=\beta_{15}^5+2a\alpha_2,
\]

\begin{align*}
\chi_{15}&=l_1\beta_{15}^1+l_3\beta_{15}^3+l_5\beta_{15}^5+l_6\beta_{15}^6 =\frac{\sqrt{2}}{4a+2}(a+\frac{1}{2}, 2a+\frac{1}{2}, 2a, 11a+3), \\
\end{align*}
where  $l_1=\frac{34a+11}{4a+2}>0, l_3=\frac{18a+7}{4a+2}>0, l_5=\frac{3a+2}{2a+1}>0, l_6=\frac{21a+6}{2a+1}>0$.

\subsection{Proof of Theorem \ref{thm ds z=1}} \label{ss proof thm z=1}

The basic strategy is similar to Theorem \ref{thm regular ds}, which we write $\chi_{13}$ into a linear combination of some roots in $P(\chi_{13})$. There are two subtleties, compared with the regular cases:
\begin{enumerate}
\item $P(\chi_{13})$ contains more than 4 roots, and so is not linearly independent. 
\item A local region for a non-regular character is not necessarily skew.
\end{enumerate}
For (1), we select some elements of $P(\chi_{13})$ in  Section \ref{ss useful formulas} to write $\beta_{13}$ into their linear combinations. For (2), we shall need to use Lemmas \ref{lem rule out orthogonal condition} and \ref{lem no rank 2 subroot system} to locate suitable regions to be skew.

\begin{proposition} $\mathbb{H}^{(\chi_{13},J_{13})}$ is a discrete series.
\end{proposition} 
\begin{proof} 

\noindent
\underline{Step 1: Find the skew local region}\\
We pick $J_{13}=\left\{ \beta_{13}^1, \beta_{13}^2, \beta_{13}^3, \beta_{13}^5\right\}$. Then one checks that $s_3s_4s_3s_4s_2s_3s_1w_o\in \mathcal F^{(\chi_{13}, J_{13})}$. Thus $\mathcal F^{(\chi_{13}, J_{13})}\neq\emptyset$, and so $(\chi_{13}, J_{13})$ is a local region. Moreover, $\mathcal{F}^{(\chi_{13}, J_{13})}$ has the sign $(-,-,-,+,-,+)$. 

We next show that $\mathcal F^{(\chi_{13}, J_{13})}$ is skew. Let $w \in \mathcal F^{(\chi_{13}, J_{13})}$. Since $Z(\chi_{13})=\left\{ \alpha_4\right\}$, $Z(w\chi_{13})=\left\{ w\alpha_4\right\}$. 

We first consider the first condition of Definition \ref{def skew}. To this end, $w(\beta_{13}^4)=s_{w\alpha_4}(w\beta_{13}^3)$, but $\mathrm{sign}(w\beta_{13}^4)=-\mathrm{sign}(w\beta_{13}^3)$. Hence, by Lemma \ref{lem rule out orthogonal condition},  $w\alpha_4$ is not a simple root and so it checks the first condition of Definition \ref{def skew}.

We now consider the second condition of Definition \ref{def skew}. By the formulas in Section \ref{ss useful formulas}, $w\beta_{13}^3, w\beta_{13}^5, w\alpha_4$ are linearly independent. Moreover, since $\mathcal F^{(\chi_{13}, J_{13})}$ has the sign $(-,-,-,+,-,+)$, formulas in Section \ref{ss useful formulas} also give
\[ \mathrm{sign}(w\beta_{13}^3)=-\mathrm{sign}(w\beta_{13}^4)= -\mathrm{sign}(s_{w\alpha_4}(w\beta_{13}^3))
\]
and
\[ \mathrm{sign}(w\beta_{13}^5)=-\mathrm{sign}(w\beta_{13}^6)= -\mathrm{sign}(s_{w\alpha_4}(w\beta_{13}^5)) .
\]
Then, by Lemma \ref{lem no rank 2 subroot system}, $w\alpha_4$ cannot be in any rank $2$ subroot system. Since $Z(w\chi_{13})=\left\{ w\alpha_4\right\}$, we then have
\[  \left\{ \alpha \in R': \langle \alpha, w\chi_{13}\rangle =0 \right\} =\emptyset
\]
    for any rank $2$ subroot system $R'$ of $R$. This verifies the second condition of Definition \ref{def skew}. \\

\noindent
\underline{Step 2: Show $\mathbb H^{(\chi_{13}, J_{13})}$ is a discrete series.}

By Theorem \ref{thm construct calibrated module} and Step 1, $\mathbb H^{(\chi_{13}, J_{13})}$ exists. It follows from Section \ref{ss useful formulas} that 
\[  \chi_{13}=c_1\beta^1_{13}+c_2\beta^2_{13}+c_3\beta^3_{13}+c_5\beta^5_{13} 
\]
and so for $w \in \mathcal F^{(\chi_{13}, J_{13})}$, 
 $c_1w\beta^1_{13}<0,\ c_2w\beta^2_{13}<0,\ c_3w\beta^3_{13}<0,\ c_5w\beta^5_{13}<0$. Thus,  $(w\chi_{13})(\varpi_i)<0$ for $i=1,2,3,4$ by Proposition \ref{prop for ds criteria}. This implies that $\mathbb H^{(\chi_{13}, J_{13})}$ is a discrete series.
\end{proof}

\subsection{Explicit structure of $\mathbb H^{(\chi_{13}, J_{13})},\mathbb H^{(\chi_{14}, J_{14})}, \mathbb H^{(\chi_{15}, J_{15})}$} \label{ss explicit module str} We use three graphs to show all elements in $\mathcal F^{(\chi_{13}, J_{13})}, \mathcal F^{(\chi_{14}, J_{14})}, \mathcal F^{(\chi_{15}, J_{15})}$. With the descriptions in Theorem \ref{thm construct calibrated module}, we have quite explicit structures of $\mathbb{H}$-actions on $\mathbb H^{(\chi_{j}, J_{j})}$ ($j=13,14,15$). In the following graphs, the arrow means length increasing and the simple reflection means left multiplication by the simple reflection (for example, in the following diagram: $w^3_{13}=s_2w^1_{13}$ and $l(w^3_{13})=l(w^1_{13})+1$).

\begin{enumerate}
\item The 81 elements in $\mathcal F^{(\chi_{13}, J_{13})}$ are $\{w_{13}^1,\cdots,w_{13}^{81}\}$. Here $w_{13}^{1}=s_1s_2s_3s_4s_3s_4s_1s_2s_3s_4$\\$s_1s_2s_3s_4s_2s_3s_4s_1s_2s_3s_2s_1w_o$. We also record that the corresponding element in SageMath of $w_{13}^1$ is $W[10980]$. (Here $W$=CoxeterGroup(['H'],4) in SageMath code.)

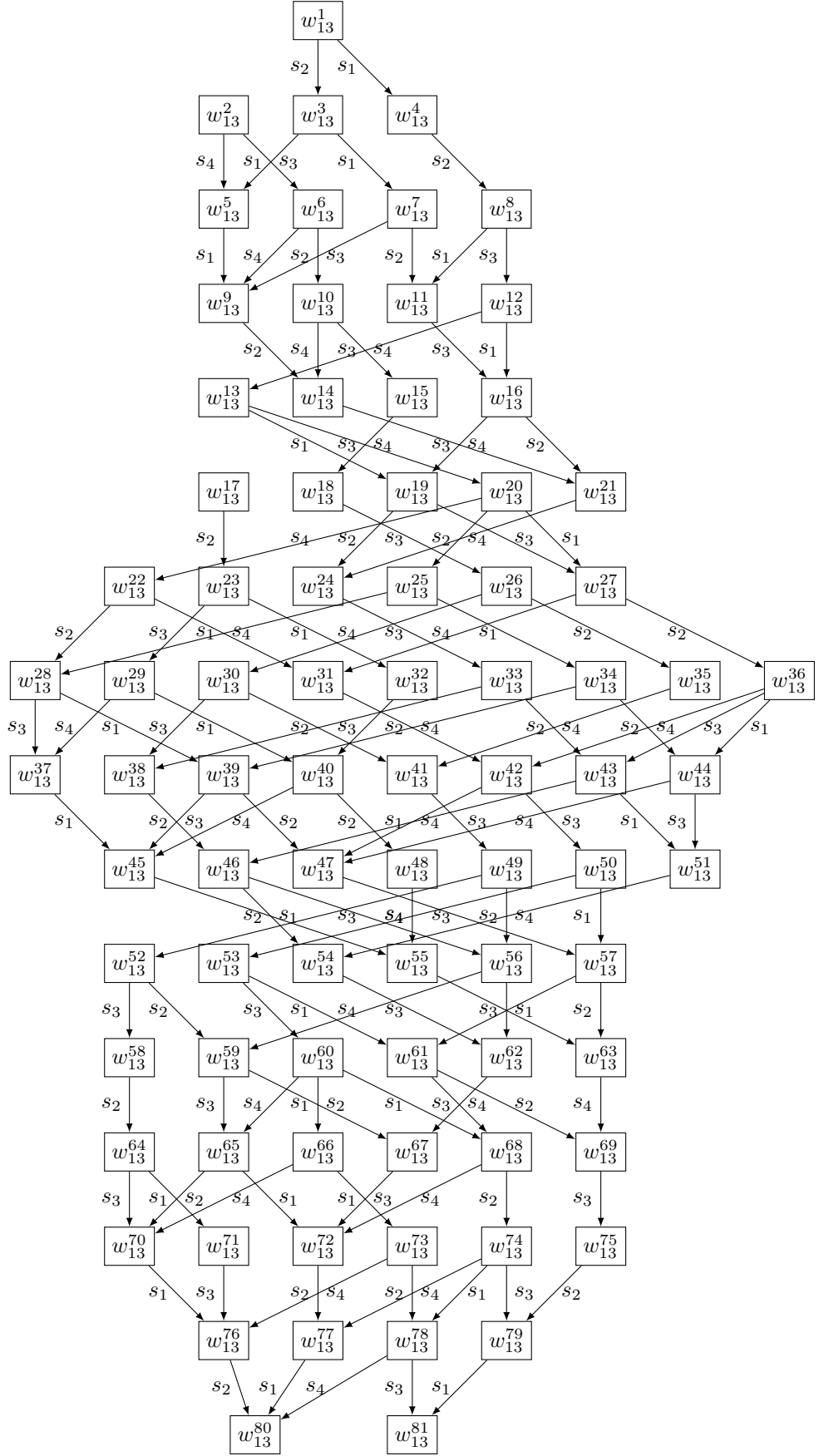
\tikzset{int/.style={draw,minimum size=1em}}
\begin{tikzpicture}[auto,>=latex]
 \node [int] (1node) {$w_{13}^{80}$};
 \node (2) [right of=1node, node distance=1.5cm, coordinate] {1};
 \node (3) [above of=1node, left of=1node, node distance=1.5cm, coordinate] {1};
 \node (4) [above of=1node, node distance=1.5cm, coordinate] {1};
 \node (5) [above of=2, node distance=1.5cm, coordinate] {2};
 \node (6) [right of=5, node distance=1.5cm, coordinate] {5};
 \node [int] (2node) [right of=2, distance=0cm] {$w_{13}^{81}$}; 
 \node [int] (3node) [right of=3, distance=0cm] {$w_{13}^{76}$};
 \node [int] (4node) [right of=4, distance=0cm] {$w_{13}^{77}$};
 \node [int] (5node) [right of=5, distance=0cm] {$w_{13}^{78}$};
 \node [int] (6node) [right of=6, distance=0cm] {$w_{13}^{79}$};
 \path[<-] (1node) edge node[left] {$s_{2}$} (3node);
 \path[<-] (1node) edge node[left] {$s_{1}$} (4node);
 \path[<-] (1node) edge node[left] {$s_{4}$} (5node);
 \path[<-] (2node) edge node[left] {$s_{3}$} (5node);
 \path[<-] (2node) edge node[left] {$s_{1}$} (6node);
 \node (7) [above of=3, left of=3, node distance=1.5cm, coordinate] {3};
 \node (8) [above of=3, node distance=1.5cm, coordinate] {3};
 \node (9) [above of=4, node distance=1.5cm, coordinate] {4};
 \node (10) [above of=5, node distance=1.5cm, coordinate] {5};
 \node (11) [above of=6, node distance=1.5cm, coordinate] {6};
 \node (12) [above of=6, right of=6, node distance=1.5cm, coordinate] {6};
 \node [int] (7node) [right of=7, distance=0cm] {$w_{13}^{70}$}; 
 \node [int] (8node) [right of=8, distance=0cm] {$w_{13}^{71}$};
 \node [int] (9node) [right of=9, distance=0cm] {$w_{13}^{72}$};
 \node [int] (10node) [right of=10, distance=0cm] {$w_{13}^{73}$};
 \node [int] (11node) [right of=11, distance=0cm] 
 {$w_{13}^{74}$};
 \node [int] (12node) [right of=12, distance=0cm] {$w_{13}^{75}$};
 \path[<-] (3node) edge node[left] {$s_{1}$} (7node);
 \path[<-] (3node) edge node[left] {$s_{3}$} (8node);
 \path[<-] (3node) edge node[right] {$s_{4}$} (10node);
 \path[<-] (4node) edge node[left] {$s_{2}$} (9node);
 \path[<-] (4node) edge node[right] {$s_{4}$} (11node);
 \path[<-] (5node) edge node[left] {$s_{2}$} (10node);
 \path[<-] (5node) edge node[right] {$s_{1}$} (11node);
 \path[<-] (6node) edge node[right] {$s_{3}$} (11node);
 \path[<-] (6node) edge node[right] {$s_{2}$} (12node);
 \node (13) [above of=7, node distance=1.5cm, coordinate] {7};
 \node (14) [above of=8, node distance=1.5cm, coordinate] {8};
 \node (15) [above of=9, node distance=1.5cm, coordinate] {9};
 \node (16) [above of=10, node distance=1.5cm, coordinate] {10};
 \node (17) [above of=11, node distance=1.5cm, coordinate] {11};
 \node (18) [above of=12, node distance=1.5cm, coordinate] {12};
 \node [int] (13node) [right of=13, distance=0cm] {$w_{13}^{64}$}; 
 \node [int] (14node) [right of=14, distance=0cm] {$w_{13}^{65}$};
 \node [int] (15node) [right of=15, distance=0cm] {$w_{13}^{66}$};
 \node [int] (16node) [right of=16, distance=0cm] {$w_{13}^{67}$};
 \node [int] (17node) [right of=17, distance=0cm] 
 {$w_{13}^{68}$};
 \node [int] (18node) [right of=18, distance=0cm] {$w_{13}^{69}$};
 \path[<-] (7node) edge node[left] {$s_{3}$} (13node);
 \path[<-] (7node) edge node[right] {$s_{2}$} (14node);
 \path[<-] (7node) edge node[right] {$s_4$} (15node);
 \path[<-] (8node) edge node[left] {$s_{1}$} (13node);
 \path[<-] (9node) edge node[right] {$s_{1}$} (14node);
 \path[<-] (9node) edge node[right] {$s_{3}$} (16node);
 \path[<-] (9node) edge node[right] {$s_{4}$} (17node);
 \path[<-] (10node) edge node[left] {$s_{1}$} (15node);
 \path[<-] (11node) edge node[left] {$s_{2}$} (17node);
 \path[<-] (12node) edge node[left] {$s_{3}$} (18node);
 \node (19) [above of=13, node distance=1.5cm, coordinate] {13};
 \node (20) [above of=14, node distance=1.5cm, coordinate] {14};
 \node (21) [above of=15, node distance=1.5cm, coordinate] {15};
 \node (22) [above of=16, node distance=1.5cm, coordinate] {16};
 \node (23) [above of=17, node distance=1.5cm, coordinate] {17};
 \node (24) [above of=18, node distance=1.5cm, coordinate] {18};
 \node [int] (19node) [right of=19, distance=0cm] {$w_{13}^{58}$}; 
 \node [int] (20node) [right of=20, distance=0cm] {$w_{13}^{59}$};
 \node [int] (21node) [right of=21, distance=0cm] {$w_{13}^{60}$};
 \node [int] (22node) [right of=22, distance=0cm] {$w_{13}^{61}$};
 \node [int] (23node) [right of=23, distance=0cm] 
 {$w_{13}^{62}$};
 \node [int] (24node) [right of=24, distance=0cm] {$w_{13}^{63}$};
 \path[<-] (13node) edge node[left] {$s_{2}$} (19node);
 \path[<-] (14node) edge node[left] {$s_{3}$} (20node);
 \path[<-] (14node) edge node[left] {$s_4$} (21node);
 \path[<-] (15node) edge node[right] {$s_{2}$} (21node);
 \path[<-] (16node) edge node[left] {$s_{1}$} (20node);
 \path[<-] (16node) edge node[right] {$s_{4}$} (23node);
 \path[<-] (17node) edge node[left] {$s_{1}$} (21node);
 \path[<-] (17node) edge node[left] {$s_{3}$} (22node);
 \path[<-] (18node) edge node[right] {$s_{2}$} (22node);
 \path[<-] (18node) edge node[left] {$s_{4}$} (24node);
 \node (25) [above of=19, node distance=1.5cm, coordinate] {19};
 \node (26) [above of=20, node distance=1.5cm, coordinate] {20};
 \node (27) [above of=21, node distance=1.5cm, coordinate] {21};
 \node (28) [above of=22, node distance=1.5cm, coordinate] {22};
 \node (29) [above of=23, node distance=1.5cm, coordinate] {23};
 \node (30) [above of=24, node distance=1.5cm, coordinate] {24};
 \node [int] (25node) [right of=25, distance=0cm] {$w_{13}^{52}$}; 
 \node [int] (26node) [right of=26, distance=0cm] {$w_{13}^{53}$};
 \node [int] (27node) [right of=27, distance=0cm] {$w_{13}^{54}$};
 \node [int] (28node) [right of=28, distance=0cm] {$w_{13}^{55}$};
 \node [int] (29node) [right of=29, distance=0cm] 
 {$w_{13}^{56}$};
 \node [int] (30node) [right of=30, distance=0cm] {$w_{13}^{57}$};
 \path[<-] (19node) edge node[left] {$s_{3}$} (25node);
 \path[<-] (20node) edge node[left] {$s_{2}$} (25node);
 \path[<-] (20node) edge node[left] {$s_4$} (29node);
 \path[<-] (21node) edge node[left] {$s_{3}$} (26node);
 \path[<-] (22node) edge node[left] {$s_{1}$} (26node);
 \path[<-] (22node) edge node[left] {} (30node);
 \path[<-] (23node) edge node[left] {$s_{3}$} (27node);
 \path[<-] (23node) edge node[right] {$s_1$} (29node);
 \path[<-] (24node) edge node[left] {$s_{3}$} (28node);
 \path[<-] (24node) edge node[left] {$s_{2}$} (30node);
 \node (31) [above of=25, node distance=1.5cm, coordinate] {25};
 \node (32) [above of=26, node distance=1.5cm, coordinate] {26};
 \node (33) [above of=27, node distance=1.5cm, coordinate] {27};
 \node (34) [above of=28, node distance=1.5cm, coordinate] {28};
 \node (35) [above of=29, node distance=1.5cm, coordinate] {29};
 \node (36) [above of=30, node distance=1.5cm, coordinate] {30};
 \node (37) [above of=30, right of=30, node distance=1.5cm, coordinate] {30};
 \node [int] (31node) [right of=31, distance=0cm] {$w_{13}^{45}$}; 
 \node [int] (32node) [right of=32, distance=0cm] {$w_{13}^{46}$};
 \node [int] (33node) [right of=33, distance=0cm] {$w_{13}^{47}$};
 \node [int] (34node) [right of=34, distance=0cm] {$w_{13}^{48}$};
 \node [int] (35node) [right of=35, distance=0cm] 
 {$w_{13}^{49}$};
 \node [int] (36node) [right of=36, distance=0cm] {$w_{13}^{50}$};
 \node [int] (37node) [right of=37, distance=0cm] {$w_{13}^{51}$};
 \path[<-] (25node) edge node[left] {} (35node);
 \path[<-] (26node) edge node[left] {$s_{4}$} (36node);
 \path[<-] (27node) edge node[right] {$s_1$} (32node);
 \path[<-] (27node) edge node[right] {$s_{4}$} (37node);
 \path[<-] (28node) edge node[left] {$s_{2}$} (31node);
 \path[<-] (28node) edge node[left] {$s_{4}$} (34node);
 \path[<-] (29node) edge node[left] {$s_3$} (32node);
 \path[<-] (29node) edge node[left] {$s_2$} (35node);
 \path[<-] (30node) edge node[left] {$s_3$} (33node);
 \path[<-] (30node) edge node[left] {$s_{1}$} (36node);
 \node (38) [above of=31, left of=31, node distance=1.5cm, coordinate] {31};
 \node (39) [above of=31, node distance=1.5cm, coordinate] {31};
 \node (40) [above of=32, node distance=1.5cm, coordinate] {32};
 \node (41) [above of=33, node distance=1.5cm, coordinate] {33};
 \node (42) [above of=34, node distance=1.5cm, coordinate] {34};
 \node (43) [above of=35, node distance=1.5cm, coordinate] {35};
 \node (44) [above of=36, node distance=1.5cm, coordinate] {36};
 \node (45) [above of=37, node distance=1.5cm, coordinate] {37};
 \node [int] (38node) [right of=38, distance=0cm] {$w_{13}^{37}$}; 
 \node [int] (39node) [right of=39, distance=0cm] {$w_{13}^{38}$};
 \node [int] (40node) [right of=40, distance=0cm] {$w_{13}^{39}$};
 \node [int] (41node) [right of=41, distance=0cm] {$w_{13}^{40}$};
 \node [int] (42node) [right of=42, distance=0cm] 
 {$w_{13}^{41}$};
 \node [int] (43node) [right of=43, distance=0cm] {$w_{13}^{42}$};
 \node [int] (44node) [right of=44, distance=0cm] {$w_{13}^{43}$};
 \node [int] (45node) [right of=45, distance=0cm] {$w_{13}^{44}$};
 \path[<-] (31node) edge node[left] {$s_{1}$} (38node);
 \path[<-] (31node) edge node[right] {$s_{3}$} (40node);
 \path[<-] (31node) edge node[right] {$s_4$} (41node);
 \path[<-] (32node) edge node[left] {$s_{2}$} (39node);
 \path[<-] (32node) edge node[right] {$s_{4}$} (44node);
 \path[<-] (33node) edge node[right] {$s_{2}$} (40node);
 \path[<-] (33node) edge node[right] {$s_4$} (45node);
 \path[<-] (33node) edge node[left] {$s_1$} (43node);
 \path[<-] (34node) edge node[left] {$s_{2}$} (41node);
 \path[<-] (35node) edge node[right] {$s_{3}$} (42node);
 \path[<-] (36node) edge node[right] {$s_3$} (43node);
 \path[<-] (37node) edge node[left] {$s_{1}$} (44node);
 \path[<-] (37node) edge node[left] {$s_{3}$} (45node);
 \node (46) [above of=38, node distance=1.5cm, coordinate] {38};
 \node (47) [above of=39, node distance=1.5cm, coordinate] {39};
 \node (48) [above of=40, node distance=1.5cm, coordinate] {40};
 \node (49) [above of=41, node distance=1.5cm, coordinate] {41};
 \node (50) [above of=42, node distance=1.5cm, coordinate] {42};
 \node (51) [above of=43, node distance=1.5cm, coordinate] {43};
 \node (52) [above of=44, node distance=1.5cm, coordinate] {44};
 \node (53) [above of=45, node distance=1.5cm, coordinate] {45};
 \node (54) [above of=45, right of=45, node distance=1.5cm, coordinate] {45};
 \node [int] (46node) [right of=46, distance=0cm] {$w_{13}^{28}$}; 
 \node [int] (47node) [right of=47, distance=0cm] {$w_{13}^{29}$};
 \node [int] (48node) [right of=48, distance=0cm] {$w_{13}^{30}$};
 \node [int] (49node) [right of=49, distance=0cm] {$w_{13}^{31}$};
 \node [int] (50node) [right of=50, distance=0cm] 
 {$w_{13}^{32}$};
 \node [int] (51node) [right of=51, distance=0cm] {$w_{13}^{33}$};
 \node [int] (52node) [right of=52, distance=0cm] {$w_{13}^{34}$};
 \node [int] (53node) [right of=53, distance=0cm] {$w_{13}^{35}$};
 \node [int] (54node) [right of=54, distance=0cm] {$w_{13}^{36}$};
 \node (55) [above of=47, node distance=1.5cm, coordinate] {47};
 \node (56) [above of=48, node distance=1.5cm, coordinate] {48};
 \node (57) [above of=49, node distance=1.5cm, coordinate] {49};
 \node (58) [above of=50, node distance=1.5cm, coordinate] {50};
 \node (59) [above of=51, node distance=1.5cm, coordinate] {51};
 \node (60) [above of=52, node distance=1.5cm, coordinate] {52};
 \node (61) [above of=56, node distance=1.5cm, coordinate] {56};
 \node (62) [above of=57, node distance=1.5cm, coordinate] {57};
 \node (63) [above of=58, node distance=1.5cm, coordinate] {58};
 \node (64) [above of=59, node distance=1.5cm, coordinate] {59};
 \node (65) [above of=60, node distance=1.5cm, coordinate] {60};
 \node (66) [above of=61, node distance=1.5cm, coordinate] {61};
 \node (67) [above of=62, node distance=1.5cm, coordinate] {62};
 \node (68) [above of=63, node distance=1.5cm, coordinate] {63};
 \node (69) [above of=64, node distance=1.5cm, coordinate] {64};
 \node (70) [above of=66, node distance=1.5cm, coordinate] {66};
 \node (71) [above of=67, node distance=1.5cm, coordinate] {67};
 \node (72) [above of=68, node distance=1.5cm, coordinate] {68};
 \node (73) [above of=69, node distance=1.5cm, coordinate] {69};
 \node (74) [above of=70, node distance=1.5cm, coordinate] {70};
 \node (75) [above of=71, node distance=1.5cm, coordinate] {71};
 \node (76) [above of=72, node distance=1.5cm, coordinate] {72};
 \node (77) [above of=73, node distance=1.5cm, coordinate] {73};
 \node (78) [above of=74, node distance=1.5cm, coordinate] {74};
 \node (79) [above of=75, node distance=1.5cm, coordinate] {75};
 \node (80) [above of=76, node distance=1.5cm, coordinate] {76};
 \node (81) [above of=79, node distance=1.5cm, coordinate] {79};
 \node [int] (55node) [right of=55, distance=0cm] {$w_{13}^{22}$}; 
 \node [int] (56node) [right of=56, distance=0cm] {$w_{13}^{23}$};
 \node [int] (57node) [right of=57, distance=0cm] {$w_{13}^{24}$};
 \node [int] (58node) [right of=58, distance=0cm] {$w_{13}^{25}$};
 \node [int] (59node) [right of=59, distance=0cm] 
 {$w_{13}^{26}$};
 \node [int] (60node) [right of=60, distance=0cm] {$w_{13}^{27}$};
 \node [int] (61node) [right of=61, distance=0cm] {$w_{13}^{17}$};
 \node [int] (62node) [right of=62, distance=0cm] {$w_{13}^{18}$};
 \node [int] (63node) [right of=63, distance=0cm] {$w_{13}^{19}$};
 \node [int] (64node) [right of=64, distance=0cm] {$w_{13}^{20}$}; 
 \node [int] (65node) [right of=65, distance=0cm] {$w_{13}^{21}$};
 \node [int] (66node) [right of=66, distance=0cm] {$w_{13}^{13}$};
 \node [int] (67node) [right of=67, distance=0cm] {$w_{13}^{14}$};
 \node [int] (68node) [right of=68, distance=0cm] 
 {$w_{13}^{15}$};
 \node [int] (69node) [right of=69, distance=0cm] {$w_{13}^{16}$};
 \node [int] (70node) [right of=70, distance=0cm] {$w_{13}^{9}$};
 \node [int] (71node) [right of=71, distance=0cm] {$w_{13}^{10}$};
 \node [int] (72node) [right of=72, distance=0cm] {$w_{13}^{11}$};
 \node [int] (73node) [right of=73, distance=0cm] {$w_{13}^{12}$}; 
 \node [int] (74node) [right of=74, distance=0cm] {$w_{13}^{5}$};
 \node [int] (75node) [right of=75, distance=0cm] {$w_{13}^{6}$};
 \node [int] (76node) [right of=76, distance=0cm] {$w_{13}^{7}$};
 \node [int] (77node) [right of=77, distance=0cm] 
 {$w_{13}^{8}$};
 \node [int] (78node) [right of=78, distance=0cm] {$w_{13}^{2}$};
 \node [int] (79node) [right of=79, distance=0cm] {$w_{13}^{3}$};
 \node [int] (80node) [right of=80, distance=0cm] {$w_{13}^{4}$};
 \node [int] (81node) [right of=81, distance=0cm] {$w_{13}^{1}$};
 \path[<-] (38node) edge node[left] {$s_{3}$} (46node);
 \path[<-] (38node) edge node[left] {$s_{4}$} (47node);
 \path[<-] (39node) edge node[left] {$s_{3}$} (48node);
 \path[<-] (39node) edge node[left] {} (51node);
 \path[<-] (40node) edge node[left] {$s_{1}$} (46node);
 \path[<-] (40node) edge node[right] {$s_{4}$} (52node);
 \path[<-] (41node) edge node[left] {$s_{1}$} (47node);
 \path[<-] (41node) edge node[left] {$s_3$} (50node);
 \path[<-] (42node) edge node[left] {$s_{2}$} (48node);
 \path[<-] (42node) edge node[right] {$s_{4}$} (53node);
 \path[<-] (43node) edge node[left] {$s_{2}$} (49node);
 \path[<-] (43node) edge node[right] {$s_{4}$} (54node);
 \path[<-] (44node) edge node[left] {$s_{2}$} (51node);
 \path[<-] (44node) edge node[right] {$s_3$} (54node);
 \path[<-] (45node) edge node[left] {$s_{2}$} (52node);
 \path[<-] (45node) edge node[right] {$s_{1}$} (54node);
 \path[<-] (46node) edge node[left] {$s_{2}$} (55node);
 \path[<-] (46node) edge node[right] {$s_4$} (58node);
 \path[<-] (47node) edge node[left] {$s_{3}$} (56node);
 \path[<-] (48node) edge node[left] {$s_{4}$} (59node);
 \path[<-] (49node) edge node[left] {$s_{1}$} (55node);
 \path[<-] (49node) edge node[left] {$s_4$} (60node);
 \path[<-] (50node) edge node[left] {$s_1$} (56node);
 \path[<-] (51node) edge node[left] {$s_{3}$} (57node);
 \path[<-] (52node) edge node[left] {$s_{1}$} (58node);
 \path[<-] (53node) edge node[left] {$s_{2}$} (59node);
 \path[<-] (54node) edge node[left] {$s_2$} (60node);
 \path[<-] (55node) edge node[left] {$s_{4}$} (64node);
 \path[<-] (56node) edge node[left] {$s_2$} (61node);
 \path[<-] (57node) edge node[left] {$s_{2}$} (63node);
 \path[<-] (57node) edge node[right] {$s_4$} (65node);
 \path[<-] (58node) edge node[left] {$s_{2}$} (64node);
 \path[<-] (59node) edge node[left] {$s_3$} (62node);
 \path[<-] (60node) edge node[right] {$s_{3}$} (63node);
 \path[<-] (60node) edge node[right] {$s_1$} (64node);
 \path[<-] (62node) edge node[right] {$s_4$} (68node);
 \path[<-] (63node) edge node[left] {$s_{1}$} (66node);
 \path[<-] (63node) edge node[right] {$s_4$} (69node);
 \path[<-] (64node) edge node[left] {$s_3$} (66node);
 \path[<-] (65node) edge node[left] {$s_{3}$} (67node);
 \path[<-] (65node) edge node[left] {$s_2$} (69node);
 \path[<-] (66node) edge node[right] {$s_4$} (73node);
 \path[<-] (67node) edge node[left] {$s_{2}$} (70node);
 \path[<-] (67node) edge node[left] {$s_4$} (71node);
 \path[<-] (68node) edge node[left] {$s_3$} (71node);
 \path[<-] (69node) edge node[left] {$s_{3}$} (72node);
 \path[<-] (69node) edge node[left] {$s_1$} (73node);
 \path[<-] (70node) edge node[left] {$s_1$} (74node);
 \path[<-] (70node) edge node[left] {$s_{4}$} (75node);
 \path[<-] (70node) edge node[right] {$s_3$} (76node);
 \path[<-] (71node) edge node[left] {$s_2$} (75node);
 \path[<-] (72node) edge node[left] {$s_{2}$} (76node);
 \path[<-] (72node) edge node[left] {$s_1$} (77node);
 \path[<-] (73node) edge node[left] {$s_3$} (77node);
 \path[<-] (74node) edge node[left] {$s_4$} (78node);
 \path[<-] (74node) edge node[right] {$s_{3}$} (79node);
 \path[<-] (75node) edge node[left] {$s_1$} (78node);
 \path[<-] (76node) edge node[left] {$s_1$} (79node);
 \path[<-] (77node) edge node[left] {$s_{2}$} (80node);
 \path[<-] (79node) edge node[left] {$s_2$} (81node);
 \path[<-] (80node) edge node[left] {$s_1$} (81node);
\end{tikzpicture}

\item[(2)] The following graph shows all the 10 elements in $\mathcal F^{(\chi_{14},J_{14})}$. Here, $w_{14}^1=s_1s_2s_3s_4s_1s_2s_3s_4s_2s_3s_4w_o$, and corresponds to $W[13857]$ in SageMath.

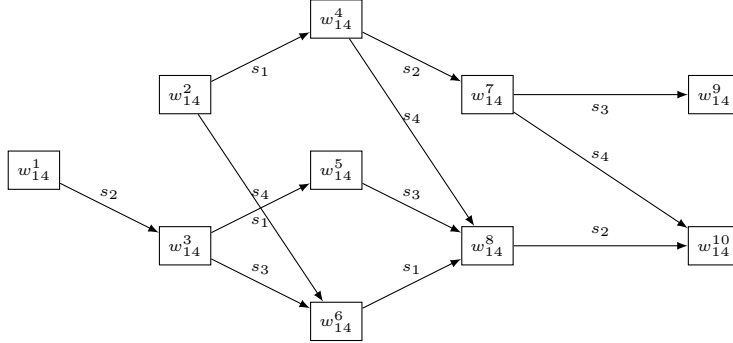
\tikzset{int/.style={draw,minimum size=1em}}
\begin{tikzpicture}[auto,>=latex]
 \node [int] (1) {\tiny{$w_{14}^1$}};
 \node (2) [above of=1, right of=1, node distance=1cm, coordinate] {1};
 \node (3) [below of=1, right of=1, node distance=1cm, coordinate] {1};
 \node [int] (2node) [right of=2, distance=0cm] {\tiny{$w_{14}^2$}}; 
 \node [int] (3node) [right of=3, distance=0cm] {\tiny{$w_{14}^3$}};
 \path[->] (1) edge node[above] {\tiny{$s_{2}$}} (3node);
 \node (4) [above of=2node, right of=2node, node distance=1cm, coordinate] {2};
 \node (5) [below of=2node, right of=2node, node distance=1cm, coordinate] {2};
 \node (6) [below of=3node, right of=3node, node distance=1cm, coordinate] {3};
 \node [int] (4node) [right of=4, distance=0cm] {\tiny{$w_{14}^4$}}; 
 \node [int] (5node) [right of=5, distance=0cm] {\tiny{$w_{14}^5$}};
 \node [int] (6node) [right of=6, distance=0cm] {\tiny{$w_{14}^6$}};
 \path[->] (2node) edge node[below] {\tiny{$s_{1}$}} (4node);
 \path[->] (2node) edge node[above] {\tiny{$s_{4}$}} (6node);
 \path[->] (3node) edge node[below] {\tiny{$s_{1}$}} (5node);
 \path[->] (3node) edge node[above] {\tiny{$s_{3}$}} (6node);
 \node (7) [below of=4node, right of=4node, node distance=1cm, coordinate] {4};
 \node (8) [below of=5node, right of=5node, node distance=1cm, coordinate] {5};
 \node [int] (7node) [right of=7, distance=0cm] {\tiny{$w_{14}^7$}}; 
 \node [int] (8node) [right of=8, distance=0cm] {\tiny{$w_{14}^8$}};
 \path[->] (4node) edge node[below] {\tiny{$s_{2}$}} (7node);
 \path[->] (4node) edge node[above] {\tiny{$s_{4}$}} (8node);
 \path[->] (5node) edge node[above] {\tiny{$s_{3}$}} (8node);
 \path[->] (6node) edge node[above] {\tiny{$s_{1}$}} (8node);
 \node (9) [right of=7node, node distance=2cm, coordinate] {7};
 \node (10) [right of=8node, node distance=2cm, coordinate] {8};
 \node [int] (9node) [right of=9, distance=0cm] {\tiny{$w_{14}^9$}}; 
 \node [int] (10node) [right of=10, distance=0cm] {\tiny{$w_{14}^{10}$}};
 \path[->] (7node) edge node[below] {\tiny{$s_{3}$}} (9node);
 \path[->] (7node) edge node[above] {\tiny{$s_{4}$}} (10node);
 \path[->] (8node) edge node[above] {\tiny{$s_{2}$}} (10node);
\end{tikzpicture}
\item[(3)] The following graph shows all the 9 elements in $\mathcal F^{(\chi_{15}, J_{15})}$. Here $w_{15}^1=s_1s_2s_3s_4s_3s_4s_2s_3s_4s_1s_2w_o$, and corresponds to $W[13812]$ in the SageMath.\\
\tikzset{int/.style={draw,minimum size=1em}}
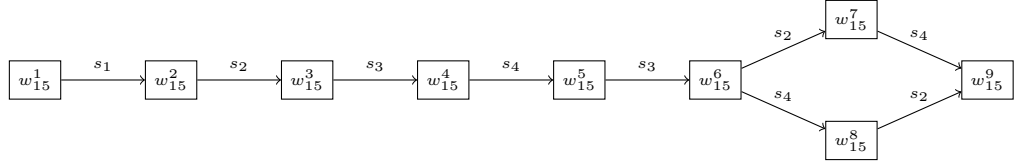
\begin{tikzpicture}
 \node [int] (1) {\tiny{$w_{15}^1$}};
 \node (2) [right of=1, node distance=0.8cm, coordinate] {1};
 \node [int] (2node) [right of=2, distance=0cm] {\tiny{$w_{15}^2$}};
 \node (3) [right of=2node, node distance=0.8cm, coordinate] {2}; 
 \node [int] (3node) [right of=3, distance=0cm] {\tiny{$w_{15}^3$}};
 \path[->] (1) edge node[above] {\tiny{$s_{1}$}} (2node);
 \path[->] (2node) edge node[above] {\tiny{$s_{2}$}} (3node);
 \node (4) [right of=3node, node distance=0.8cm, coordinate] {3};
 \node [int] (4node) [right of=4, distance=0cm] {\tiny{$w_{15}^4$}};
 \node (5) [right of=4node, node distance=0.8cm, coordinate] {4};
 \node [int] (5node) [right of=5, distance=0cm] {\tiny{$w_{15}^5$}};
 \node (6) [right of=5node, node distance=0.8cm, coordinate] {5};
 \node [int] (6node) [right of=6, distance=0cm] {\tiny{$w_{15}^6$}};
 \path[->] (3node) edge node[above] {\tiny{$s_{3}$}} (4node);
 \path[->] (4node) edge node[above] {\tiny{$s_{4}$}} (5node);
 \path[->] (5node) edge node[above] {\tiny{$s_{3}$}} (6node);
 \node (7) [above of=6node, right of=6node, node distance=0.8cm, coordinate] {6};
 \node (8) [below of=6node, right of=6node, node distance=0.8cm, coordinate] {6};
 \node [int] (7node) [right of=7, distance=0cm] {\tiny{$w_{15}^7$}}; 
 \node [int] (8node) [right of=8, distance=0cm] {\tiny{$w_{15}^8$}};
 \path[->] (6node) edge node[above] {\tiny{$s_{2}$}} (7node);
 \path[->] (6node) edge node[above] {\tiny{$s_{4}$}} (8node);
 \node (9) [below of=7 node, right of=7node, node distance=0.8cm, coordinate] {7};
 \node [int] (9node) [right of=9, distance=0cm] {\tiny{$w_{15}^9$}}; 
 \path[->] (7node) edge node[above] {\tiny{$s_{4}$}} (9node);
 \path[->] (8node) edge node[above] {\tiny{$s_{2}$}} (9node);
\end{tikzpicture}
\end{enumerate}

\section{Discrete series at central characters $\chi$ with $|Z(\chi)|=2$} 

\subsection{The character $\chi_{16}$}\label{s character 16}

Let $\chi_{16}$ be the character determined by $Z(\chi_{16})=\left\{ \alpha_1, \alpha_4 \right\}$ and $P(\chi_{16})=\{\beta_{16}^1,\beta_{16}^2,\cdots,\beta_{16}^8\}$, where 
\[ \beta_{16}^1=\sqrt{2}(0,1/2,-a,-b)=-\alpha_1-\alpha_2-(2a+1)\alpha_3-(2a+1)\alpha_4,\quad \beta_{16}^2=s_4(\beta_{16}^1)=\beta_{16}^1+\alpha_4, \]
\[ \beta_{16}^3=s_1(\beta_{16}^1)=\beta_{16}^1+\alpha_1, \quad \beta_{16}^4=s_4s_1(\beta_{16}^1)=\beta_{16}^1+\alpha_1+\alpha_4, \]
\[ \beta_{16}^5=\sqrt{2}(0,a,b,-1/2)=-2a\alpha_1-2a\alpha_2-2a\alpha_3-2a\alpha_4, \quad\beta_{16}^6=s_4(\beta_{16}^5)=\beta_{16}^5+2b\alpha_4, \]
\[\beta_{16}^7=s_1(\beta_{16}^5)=\beta_{16}^5+2a\alpha_1, \quad \beta_{16}^8=s_4s_1(\beta_{16}^5)=\beta_{16}^5+2a\alpha_1+2b\alpha_4.\]
Explicitly,
\begin{align} \label{eqn char 16}
\chi_{16}& = \beta_{16}^2+\beta_{16}^3+\frac{3}{2}\beta_{16}^6+\frac{3}{2}\beta_{16}^7  = \sqrt{2}(-\frac{b}{2}, -b, \frac{b}{2}-\frac{1}{4}, -a+\frac{1}{4}).
\end{align}
Let $J_{16}=\left\{ \beta_{16}^1, \beta_{16}^2, \beta_{16}^3, \beta_{16}^5, \beta_{16}^6, \beta_{16}^7\right\}$.
\textbf{}

\subsection{$\mathbb H^{(\chi_{16}, J_{16})}$}

\begin{theorem} \label{thm chi16 ds}
$(\chi_{16}, J_{16})$ is a skew local region and the sign of $\mathcal{F}^{(\chi_{16},J_{16})}$ is equal to $(-,-,-,+,-,-,-,+)$. Moreover, $\mathbb{H}^{(\chi_{16}, J_{16})}$ is a discrete series, and its dimension is equal to 35.
\end{theorem}
\begin{proof} First, we can check $w=s_3s_4s_3s_4s_2s_3\in\mathcal{F}^{(\chi_{16},J_{16})}$. Thus $(\chi_{16},J_{16})$ is a local region. We also remark that the detailed weight structure is given in (\ref{diagram local region 16}) in an appendix.

We now show $(\chi_{16},J_{16})$ is skew. We first check the first condition of Definition \ref{def skew}. Let $w \in \mathcal F^{(\chi_{16}, J_{16})}$. Note that 
\[    \mathrm{sign}(w\beta_{16}^2) =-\mathrm{sign}(w\beta_{16}^4)=-\mathrm{sign}(s_{w\alpha_1}(w\beta_{16}^2)), 
\]
\[   \mathrm{sign}(w\beta_{16}^3) =-\mathrm{sign}(w\beta_{16}^4)=-\mathrm{sign}(s_{w\alpha_4}(w\beta_{16}^3))
\]
By Lemma \ref{lem rule out orthogonal condition}, $w(\alpha_1)$ and $w(\alpha_4)$ are not simple roots. This implies that $Z(w\chi_{16})=\left\{ w\alpha_1, w\alpha_4\right\}$ does not contain a simple root, and so the first condition of Definition \ref{def skew} is satisfied.

We now check the second condition of Definition \ref{def skew}. We have: 
\begin{enumerate}
\item $w\beta_{16}^2, w\beta_{16}^6, w\alpha_1$ are linearly independent;
\item $\mathrm{sign}(w\beta_{16}^2)=-\mathrm{sign}(w\beta_{16}^4)=-\mathrm{sign}(s_{w\alpha_1}(w\beta_{16}^2))$; 
\item $\mathrm{sign}(w\beta_{16}^6)=-\mathrm{sign}(w\beta_{16}^8)=-\mathrm{sign}(s_{w\alpha_1}(w\beta_{16}^6)$.
\end{enumerate}
Thus, Lemma \ref{lem no rank 2 subroot system} implies that $w\alpha_1$ is not in any rank $2$ subroot system of $R$. Now, one considers $w\beta_{16}^3, w\beta_{16}^7, w\alpha_4$ and applies a similar argument to conclude that $w\alpha_4$ is not in any rank $2$ subroot system of $R$. This implies that for any rank $2$ subroot system $R'$ of $R$, 
\[  \left\{ \alpha \in R': \langle \alpha, w\chi_{16}\rangle =0 \right\}=\emptyset . \]
This verifies the second condition of Definition \ref{def skew}.

It remains to check $\mathbb H^{(\chi_{16}, J_{16})}$ is a discrete series. Let $w \in \mathcal F^{(\chi_{16}, J_{16})}$.
Since $\mathcal F^{(\chi_{16}, J_{16})}$ has the sign $(-,-,-,+,-,-,-,+)$, Proposition \ref{prop for ds criteria} and the expression (\ref{eqn char 16}) imply that $w\chi_{16}(\varpi_{i})<0$ for all $i=1,2,3,4$. This implies that $\mathbb H^{(\chi_{16}, J_{16})}$ is a discrete series.

\end{proof}

\section{Discrete series at central characters $\chi$ with $|Z(\chi)|=4$} \label{ss ds z>2}

\subsection{The character $\chi_{17}$} \label{s character 17}
Let $\chi_{17}$ be the character determined by
\[ Z(\chi_{17})=\{\alpha_{1},\alpha_2,\alpha_{4},\alpha_1+\alpha_2\} \quad \mbox{and}\quad  P(\chi_{17})=\{\beta_{17}^1,\beta_{17}^2,\cdots,\beta_{17}^{11},\beta_{17}^{12}\} \]
with 
\[ \beta_{17}^1=\sqrt{2}(-1/2,1/2,1/2,-1/2)=-2a\alpha_1-(2a+1)\alpha_2-(2a+1)\alpha_3-(2a+1)\alpha_4,\quad \beta_{17}^2=s_4(\beta_{17}^1)=\beta_{17}^1+\alpha_{4},
\]
\[ \beta_{17}^3=s_1(\beta_{17}^1)=\beta_{17}^1+2b\alpha_{1}, \quad \beta_{17}^4=s_4s_1(\beta_{17}^1)=\beta_{17}^1+2b\alpha_{1}+\alpha_{4},\]
\[ \beta_{17}^5=s_2(\beta_{17}^1)=\beta_{17}^1+\alpha_{2}, \quad \beta_{17}^6=s_4s_2(\beta_{17}^1)=\beta_{17}^1+\alpha_{4}+\alpha_{2},\]
\[\beta_{17}^7=s_2s_1(\beta_{17}^1)=\beta_{17}^1+2b\alpha_1+2a\alpha_{2}, \quad \beta_{17}^8=s_4s_2s_1(\beta_{17}^1)=\beta_{17}^1+2b\alpha_1+2a\alpha_{2}+\alpha_4,\]
\[ \beta_{17}^9=s_1s_2(\beta_{17}^1)=\beta_{17}^1+2a\alpha_{1}+\alpha_2, \quad \beta_{17}^{10}=s_4s_1s_2(\beta_{17}^1)=\beta_{17}^1+2a\alpha_1+\alpha_2+\alpha_4,\]
\[ \beta_{17}^{11}=s_1s_2s_1(\beta_{17}^1)=\beta_{17}^1+2a\alpha_1+2a\alpha_{2},  \quad \beta_{17}^{12}=s_4s_1s_2s_1(\beta_{17}^1)=\beta_{17}^1+2a\alpha_{1}+2a\alpha_{2} +\alpha_4. \]
Explicitly, we can write 
\begin{align} \label{eqn char 17 express}
\chi_{17}=(2\beta_{17}^2+(2a-1)\beta_{17}^5+(6a+1)\beta_{17}^6+(4a+4)\beta_{17}^9)/(2+4a). 
\end{align}
Let $J_{17}=\left\{ \beta_{17}^1, \beta_{17}^2, \beta_{17}^3, \beta_{17}^4, \beta_{17}^5, \beta_{17}^6, \beta_{17}^7, \beta_{17}^9\right\}$.


\subsection{$\mathbb{H}^{(\chi_{17},J_{17})}$}

\begin{theorem} \label{thm 17 skew}
$(\chi_{17}, J_{17})$ is a skew local region and the sign of $\mathcal{F}^{(\chi_{17},J_{17})}$ is equal to $(-,-,-,-,-,-,-,+,-,+,+,+)$. Moreover, $\mathbb{H}^{(\chi_{17},J_{17})}$ is a discrete series, and has dimension equal to $35$.
\end{theorem}
\begin{proof} First, we can check $w_{17}^1=s_3s_4s_3s_4s_1s_2s_3s_4s_3\in\mathcal{F}^{(\chi_{17},J_{17})}$. Hence, $(\chi_{17},J_{17})$ is a local region. A complete list of elements in $\mathcal F^{(\chi_{17}, J_{17})}$ is shown in (\ref{diagram local region 17}).

We are now going to show $(\chi_{17}, J_{17})$ is skew. For the first condition of Definition \ref{def skew}, one again uses Lemma \ref{lem rule out orthogonal condition}. The argument is quite similar to the one in Theorem \ref{thm chi16 ds}. For $w\in \mathcal F^{(\chi_{17}, J_{17})}$, we consider the following pairs:
\begin{enumerate}
\item $w(\beta_{17}^1)$ and $w(\beta_{17}^{11})$ to show $w(\alpha_1+\alpha_2)$ is not a simple root;
\item $w(\beta_{17}^7)$ and $w(\beta_{17}^{11})$ to show $w(\alpha_1)$ is not a simple root;
\item $w(\beta_{17}^9)$ and $w(\beta_{17}^{11})$ to show $w(\alpha_2)$ is not a simple root;
\item $w(\beta_{17}^9)$ and $w(\beta_{17}^{10})$ to show $w(\alpha_4)$ is not a simple root.
\end{enumerate}

We now consider the second condition of Definition  \ref{def skew}. The argument is again similar to the one in Theorem \ref{thm chi16 ds}. For $w \in \mathcal F^{(\chi_{17}, J_{17})}$, we consider the following cases:
\begin{enumerate}
\item $w(\beta_{17}^7), w(\beta_{17}^9), w(\alpha_4)$ to show $w(\alpha_4)$ is not in any rank $2$ subroot system of $R$;
\item $w(\beta_{17}^4), w(\beta_{17}^9), w(\alpha_2)$ to show $w(\alpha_2)$ is not in any rank $2$ subroot system of $R$;
\item $w(\beta^6_{17}), w(\beta^7_{17}), w(\alpha_1)$ to show $w(\alpha_1)$ is not in any rank $2$ subroot system of $R$;
\item $w(\beta_{17}^1), w(\beta_{17}^2), w(\alpha_1+\alpha_2)$ to show that $w(\alpha_1+\alpha_2)$ is not in any rank $2$ subroot system of $R$.
\end{enumerate}

Finally, we need to check that $\mathbb{H}^{(\chi_{17},J_{17})}$ is a discrete series. Now, similar to Theorem \ref{thm chi16 ds}, we apply Proposition \ref{prop for ds criteria} and (\ref{eqn char 17 express}) to conclude that $w\chi_{17}(\varpi_i)<0$ for $i=1,2,3,4$.

\end{proof}

\part{Construction of minimally induced modules for discrete series}

\section{Inductive construction of discrete series} \label{s inductive construct ds}



In this section, it is slightly more convenient for us to consider $\mathbb H$ to be the graded Hecke algebra associated to a root system $R$, and consider $V=V_R$. Note that $V$ is a real representation  of $W$ (i.e. admitting a $W$-invariant real structure) and so is isomorphic to $V^{\vee}$, and we shall denote the isomorphism by $\eta$.  Define the $W$-invariant bilinear form $(.,.)$ on $V$ given by
\[  (v_1, v_2) =\eta(v_1)(v_2) .
\]
Define $||v||=\sqrt{(v,v)}$, for $v \in V_0$.

\subsection{$\star$-dual module}

\begin{definition} \label{def star operation}
Define $\star: \mathbb{H}\rightarrow\mathbb{H}$ to be the linear anti-involution determined by $$p^{\star}=-t_{w_{o}}w_{o}(p)t_{w_{o}}^{-1},\quad \forall p\in S(V);$$
$$t_{w}^{\star}=t_{w}^{-1},\quad \forall w\in W .$$
\end{definition}
The above notion of $\star$ is closely related to the smooth dual for $p$-adic groups, see e.g. \cite[Corollary 5.5]{BM93}. The notion also applies for the parabolic subalgebra $\mathbb H_I$ of $\mathbb H$, and in such case, we sometimes denote by $\star_I$ to emphasis the underlying algebra.

To avoid notational complication, we shall assume that $-w_o$ acts an identity on $R$ in the remaining of this section. This happens for type $H_4$. Recall that we are also assuming that $V=V_R$, and an important consequence of this is that for any irreducible $\mathbb H$-module $X$, $X^{\star}\cong X$ (by \cite[Lemma 4.5]{Ch16} and \cite[Proposition 4.3.1]{BC13}). 


\subsection{Inductive construction of discrete series}


Recall that the notion in $\widetilde{\mathbb H}_I$ is defined in Section \ref{ss parabolic subalg}.

\begin{theorem} \cite[Proposition 7.3.88.]{Ch14} \label{thm construct minimal ds}
Let $DS$ be a discrete series of $\mathbb H$. Then there exists $I \subset \Delta$ with $|I|=|\Delta|-1$, a discrete series $DS'$ of $\widetilde{\mathbb H}_I$ and a character $\omega \in V_I^{\bot,\vee}$ such that the following properties hold:
\begin{enumerate}
\item $\omega$ satisfies $\omega(\alpha) <0$ for the unique $\alpha \in \Delta- I$;
\item there is a surjective $\mathbb H$-map $\mathbb H\otimes_{\mathbb H_I} (DS'\otimes \mathbb C_{\omega}) \twoheadrightarrow DS$; and
\item $(\mathbb H\otimes_{\mathbb H_I}(DS'\otimes \mathbb C_{\omega}))^{\star}$ is a standard $\mathbb H$-module.
\end{enumerate}
\end{theorem}

\begin{proof}

Let $n=|\Delta|$ and $\mathcal P_{n-1}$ be the set of all subsets of $\Delta$ of cardinality $n-1$. Let $\mathcal W_{DS}$ be the set of all weights of $DS$. Define a map
\[  \Phi: \mathcal P_{n-1} \times \mathcal W_{DS} \rightarrow \mathbb R 
\]
given by
\[  \Phi(I, \gamma) = -\frac{(\mathrm{Re}(\gamma), \varpi_{\alpha})}{||\varpi_{\alpha}||},
\]
where $\alpha$ is the unique element in $\Delta-I$. 

Now, we pick a pair $(I_0, \gamma_0) \in \mathcal P_{n-1} \times \mathcal W_{DS}$ such that $\Phi(I_0, \gamma_0)$ attains the value:
\begin{align} \label{eqn min ds}
       \mathrm{min}\left\{ \Phi(I,\gamma): (I, \gamma) \in \mathcal P_{n-1}\times \mathcal W_{DS} \right\} .
\end{align}
Let $\alpha_0$ be the unique element in $\Delta-I_0$. Now one picks an irreducible submodule $X$ in $DS|_{\mathbb H_{I_0}}$ whose weights $\gamma$ satisfy the property that 
\[ (\mathrm{Re}(\gamma), \varpi_{\alpha_0})= (\mathrm{Re}(\gamma_0), \varpi_{\alpha_0}) .
\]

One checks that minimality property forces that $X$ is an essentially discrete series. See \cite[Proposition 7.3.88.]{Ch14} for a detailed computation. Then, $X\cong DS'\otimes \mathbb C_{\omega}$, where $DS'$ is some discrete series of $\widetilde{\mathbb H}_{I_0}$ and
\[  \omega  = -\Phi(I_0, \gamma_0) \cdot \frac{\varpi_{\alpha_0}}{||\varpi_{\alpha_0}||} .
\]
Now Frobenius reciprocity gives a surjection from $\mathbb H\otimes_{\mathbb H_{I_0}}(DS' \otimes \mathbb C_{\omega})$ to $DS$.

This gives (2) and (1) follows from our construction. (3) follows from the fact that
\[   (\mathbb H\otimes_{\mathbb H_I}(DS'\otimes \mathbb C_{\omega}))^{\star} \cong \mathbb H\otimes_{\mathbb H_I}((DS'\otimes \mathbb C_{\omega})^{\star})\cong \mathbb H\otimes_{\mathbb H_I} (DS'{}^{\color{red}\star}\otimes \mathbb C_{-\omega})
{\color{red}.}\]

\end{proof}

\begin{definition}
For a discrete series $DS$ of $\mathbb H$, in the notations in the statement of Proposition \ref{thm construct minimal ds}, we say that $\mathbb H\otimes_{\mathbb H_I}(DS' \otimes \mathbb C_{\omega})$ is a {\it minimally induced module} associated to $DS$.
\end{definition}

\subsection{Criteria for existence of some discrete series}

We now state some sufficient conditions for ensuring the existence of some discrete series. Recall that the Langlands parameter is defined in Theorem \ref{thm langlands class gha}. For a Langlands parameter $(U,I)$, denote by $L(U,I)$ the unique irreducible quotient of $\mathrm{Ind}_{\mathbb H_I}^{\mathbb H}U$.

\begin{lemma} \label{lem general strategy}
Let $(U, I)$ be a Langlands parameter for $\mathbb H$. Suppose the followings hold:
\begin{enumerate}
\item $\mathrm{Ind}^{\mathbb H}_{\mathbb H_I}U$ is reducible;
\item the only non-tempered irreducible composition factor in $\mathrm{Ind}_{\mathbb H_I}^{\mathbb H}U$ is $L(U, I)$;
\item any tempered composition factor in $\mathrm{Ind}^{\mathbb H}_{\mathbb H_I}U$ is a discrete series.
\end{enumerate}
Let
\[ p =\mathrm{dim}~ \mathrm{Hom}_{\mathbb H}((\mathrm{Ind}_{\mathbb H_I}^{\mathbb H}U)^{\star}, (\mathrm{Ind}_{\mathbb H_I}^{\mathbb H}U))
\]
Let $K$ be the socle of $\mathrm{Ind}_{\mathbb H_I}^{\mathbb H}U$. Then, the followings hold:
\begin{enumerate}
\item[(i)] $p \geq 1$ and $\mathrm{Ind}^{\mathbb H}_{\mathbb H_I}U$ admits a short exact sequence:
\[  0 \rightarrow K  \rightarrow \mathrm{Ind}^{\mathbb H}_{\mathbb H_I}U \rightarrow L(U, I) \rightarrow 0
\]
such that $K$ is a semisimple $\mathbb H$-module with all simple composition factors being a discrete series.
\item[(ii)] if $p \leq 3$, then $K$ contains precisely $p$ mutually non-isomorphic discrete series.
\end{enumerate}
\end{lemma}

\begin{proof}
The conditions (2) and (3) guarantee that any simple composition factor in $\mathrm{Ind}^{\mathbb H}_{\mathbb H_{\color{red}I}}U$ other than $L(U, I)$ is a discrete series. Now, \cite[Theorem 1.2]{Ch16} asserts that there is no first extension between discrete series of $\mathbb H$, and so $K$ is a direct sum of discrete series. This shows the first assertion.

Now we show the second assertion. If a discrete series, say $X$, appears in $K$ with multiplicity greater than $2$, then the semisimplicity of $K$ implies the following sequence of maps:
\[   \mathrm{Hom}_{\mathbb H}(X\oplus X, X \oplus X) \hookrightarrow \mathrm{Hom}_{\mathbb H}(K, K) \hookrightarrow \mathrm{Hom}_{\mathbb H}((\mathrm{Ind}_{\mathbb H_I}^{\mathbb H}U)^{\star},  \mathrm{Ind}^{\mathbb H}_{\mathbb H_I}U) .
\]
Since the first Hom has dimension equal to $4$, the 
 the last space also has at least dimension $4$. However,  this contradicts to the given condition on $p$.
\end{proof}

\section{Splitting of induced maps} \label{s splitting i}
\label{s splitting induced map}

In Sections \ref{s splitting induced map} and \ref{ss geo lem}, we shall introduce more tools and we shall consider general graded Hecke algebra $\mathbb H$, and use notations in Section \ref{sec notation prelim}  .

\subsection{Induced map} \label{ss induced map split}

Let $I \subset \Delta$. Let $U$ be an $\mathbb H_I$-module. Define 
\[  \iota_{U}: U \hookrightarrow \mathbb H\otimes_{\mathbb H_I} U,
\]
given by $\iota_U(u)=1\otimes u$, for $u\in U$. We shall regard $\iota_U$ as an $\mathbb H_I$-module morphism.

\begin{definition}
We use the above notation. We say that $\iota_{U}$ {\it splits} if there exists an $\mathbb H_I$-module morphism $\tau_{U}:\mathbb H\otimes_{\mathbb H_I}U \rightarrow U$ such that $\tau_{U}\circ \iota_{U}=\mathrm{Id}_{U}$.
\end{definition}

It is a natural question to ask when the embedding $\iota_{U}$ splits. The classical example in this regard is when $\mathbb H\otimes_{\mathbb H_I}U$ is a standard module. This result is basically shown in the proof of \cite[Theorem 2.1]{Ev96} and \cite[Theorem 2.4]{KR02}. We give another class of examples:

\begin{proposition} \label{prop splitting calibrated}
Let $I \subset \Delta$. Let $U$ be an irreducible $\mathbb H_I$-module. Suppose $\mathbb H\otimes_{\mathbb H_I}U$ has an irreducible quotient isomorphic to a calibrated module. Then $\iota_{U}$ splits.
\end{proposition}

\begin{proof}
Let $X$ be an irreducible calibrated quotient of $\mathbb H \otimes_{\mathbb H_I}U$.  Let $\mathrm{pr}$ be the surjection of $\mathbb H \otimes_{\mathbb H_I}U$ onto $X$. By Lemma \ref{lem restrict calibrated}, $X|_{\mathbb H_I}$ is semisimple. Hence, $\mathrm{pr}(1\otimes U)$ is a direct summand in $X|_{\mathbb H_I}$ (and is non-zero by Frobenius reciprocity). Then, one has a surjection $s: X|_{\mathbb H_I} \rightarrow \mathrm{pr}(1\otimes U)$. By identifying $\mathrm{pr}(1\otimes U)$ with $U$, we have a nonzro map $s\circ \mathrm{pr}$ from $\mathbb H\otimes_{\mathbb H_I}U$ to $U$. Now $(s\circ \mathrm{pr})\circ \iota_U$ is non-zero and so by scaling if necessary, the map $s\circ \mathrm{pr}$ gives the map $\tau_{U}$ such that $\tau_U\circ \iota_U=\mathrm{id}_U$.

\end{proof}

\section{The geometric lemma and the second adjointness} \label{ss geo lem}

\subsection{The geoemtric lemma}
We shall need to analyse the structure of restricting some induced modules. The standard tool for doing that is the geometric lemma and we shall recall now. For $I, I' \subset \Delta$, let $W^{I,I'}$ be the set of minimal representatives in $W_I\setminus W/W_{I'}$. 

\begin{proposition} (Geometric lemma)
Let $I, I' \subset \Delta$. Let $U$ be an $\mathbb H_{I'}$-module. Then $(\mathbb H\otimes_{\mathbb H_{I'}} U)|_{\mathbb H_{I}}$ admits a filtration of the form:
\begin{align} \label{eqn weights geometric layer}
\mathbb H_{I} \otimes_{\mathbb H_{I \cap w(I')}} (U |_{\mathbb H_{w^{-1}(I)\cap I'}})^{w^{-1}} ,  
\end{align}
where $w$ runs for elements in $W^{I,I'}$. Moreover, the filtration is compatible with the Bruhat ordering. Here, $(M)^{w^{-1}}$ means to give an $\mathbb H_{I\cap w(I')}$ module structure on $\mathbb H_{w^{-1}(I) \cap I'}$-module $M$ via the formula: for any $x \in M$, 
\[ p.x=w^{-1}(p).x, \mbox{ for all $p \in V$},
\]
\[  t_u.x =t_{w^{-1}uw}.x, \mbox{ for all $u \in W_{I\cap w(I')}$}.
\]
\end{proposition}

\begin{proof}
This result is certainly well-known. We only sketch the proof for the convenience of the reader. For $w \in W^{I, I'}$, let 
\[  X_{\leq w} = \mathrm{span}\left\{ t_{u}t_{w'}\otimes v \right\}, \quad  \mbox{(resp. $X_{<w} = \mathrm{span}\left\{ t_{u}t_{w'}\otimes v\right\}$)},
\]
where $u$ runs for all elements in $W_I$, $w'$ runs for all elements in $W^{I,I'}$ with $w' \leq w$ (resp. $w'<w$) and $v$ runs for all elements in $U$. One then checks that $X_{<w}$ and $X_{\leq w}$ are invariant under $\mathbb H_{I}$-action, and 
\[  X_{\leq w}/X_{<w} \cong \mathbb H_{I} \otimes_{\mathbb H_{I \cap w(I')}} (U|_{\mathbb H_{w^{-1}(I)\cap I'}})^{w^{-1}}  .
\]
\end{proof}

We also need the description of the weights of the layers in (\ref{eqn weights geometric layer}), which follows from \cite[Theorem 6.4]{BM93}:

\begin{lemma} 
The weights of the module in (\ref{eqn weights geometric layer}) are of the form:
\[ w'w(\gamma)  \]
for the weights $\gamma$ of $U$ and $w'\in W_I$.
\end{lemma}

\subsection{The second adjointness}
Let $I \subset \Delta$. Let $w^I_o$ be the longest element in $W/W_I$. Let $I'=w_{o}^I(I)$, which is also a subset of $\Delta$. Note that $w_o^IW_I(w_o^I)^{-1}=W_{I'}$. Define $\theta^I: \mathbb H_{I} \rightarrow \mathbb H_{I'}$ determined by: for $v \in V \subset S(V)$ and $w \in W_I$,
\[  \theta^I(v)= w_{o}^I(v), \quad \theta^I(w)=(w_o^I)w(w_o^I)^{-1}.
\]
When $I=\Delta$, we shall simply write $\theta$ for $\theta^I$. For an $\mathbb H_{I'}$-module $X$, the induced $\mathbb H_{I}$-module, denoted by $\theta^I(X)$, has the same underlying space as $X$ and the action given by: for $x \in X$,
\[   h\cdot_{\theta^I(X)} x=\theta^I(h)\cdot_{X} x ,
\]
for all $h \in \mathbb H_I$.

It is straightforward to check that for a finite-dimensional $\mathbb H$-module $X$,
\[(X^{\star}|_{\mathbb H_I})^{\star_I} \cong \theta^I(X|_{\mathbb H_{I'}}) .\]


\begin{proposition} \label{prop second adjoint} (see \cite{BC13}) (Second adjointness)
Let $X$ be a finite-dimensional $\mathbb H$-module and let $Y$ be a finite-dimensional $\mathbb H_I$-module. Then
\[  \mathrm{Hom}_{\mathbb H}(X, \mathbb H\otimes_{\mathbb H_I}Y) \cong \mathrm{Hom}_{\mathbb H_{I'}}( X , (\theta^I)^{-1}(Y)).\]
\end{proposition}

\begin{proof}
\begin{align*}
 \mathrm{Hom}_{\mathbb H}(X, \mathbb H\otimes_{\mathbb H_I}Y) & \cong \mathrm{Hom}_{\mathbb H}(\mathbb H\otimes_{\mathbb H_I}Y^{\star}, X^{\star}) \\
 & \cong \mathrm{Hom}_{\mathbb H_I}(Y^{\star}, X^{\star}|_{\mathbb H_I}) \\
& \cong \mathrm{Hom}_{\mathbb H_I}((X^{\star}|_{\mathbb H_I})^{\star_I}, Y) \\
& \cong \mathrm{Hom}_{\mathbb H_I}(\theta^I(X|_{\mathbb H_{I'}}), Y) \\
& \cong \mathrm{Hom}_{\mathbb H_{I'}}(X, (\theta^I)^{-1}(Y)),
\end{align*}
where the first isomorphism follows from taking $\star$-duals and using $(\mathbb H \otimes_{\mathbb H_I}Y)^{\star}\cong \mathbb H\otimes_{\mathbb H_I}Y^{\star}$, the second isomorphism follows from Frobenius reciprocity, the third isomorphism follows from taking ${\star}_I$-duals, the forth isomorphism follows from the above discussion, and the last one follows from applying the $\theta^I$-dual.
\end{proof}

\begin{remark} \label{rmk action on twisted maps}
We now mention two cases that we shall use. For $H_4$, let $I=\left\{ \alpha_1, \alpha_3, \alpha_4\right\}$ or $\left\{ \alpha_2, \alpha_3, \alpha_4\right\}$. In such cases, the induced map of $\theta^{I}(I)=I$. In the latter case, $\theta^I$ acts trivially on $V_I$ and acts by $-1$ on $V_I^{\bot}$, and $\theta^I(w)=w$ for $w \in W_I$. In the former case, $\theta^I$ fixes $\alpha_1$, but switches $\alpha_3$ and $\alpha_4$, and acts by $-1$ on $V_I^{\bot}$. As a consequence, in the latter case, for any irreducible $\mathbb H_I$-module $\widetilde{U}\otimes \mathbb C_{\chi}$ (where $\widetilde{U}$ is an irreducible $\widetilde{\mathbb H}_I$-module and $\chi \in V_I^{\bot, \vee}$),
\begin{align} \label{eqn isomorphism latter}
(\theta^I)^{-1}(\widetilde{U}\otimes \mathbb C_{\chi}) \cong \widetilde{U}\otimes \mathbb C_{-\chi} \cong (\widetilde{U}\otimes \mathbb C_{\chi})^{\star_I}.
\end{align}
In the former case, if we pick $\widetilde{U}$ to be a discrete series of $\widetilde{\mathbb H}_I$, then the equation (\ref{eqn isomorphism latter}) still holds by \cite[Corollary 8.3]{Ch16}.
\end{remark}

\section{Existence of discrete series at $W\chi_{17}$} \label{s exist ds 17}




In Sections \ref{s exist ds 17} and \ref{s exist ds 16}, we shall work on the graded Hecke algebra of type $H_4$. We shall work on the central character $W\chi_{17}$ before $W\chi_{16}$ since the minimally induced module for $W\chi_{17}$ is slightly simpler.

\subsection{Some structure of an induced module} \label{ss minimal induce chi 17}
Let $I_{17}=\left\{ \alpha_1, \alpha_3, \alpha_4 \right\}$, and let $\mathrm{St}_{I_{17}}$ be the Steinberg module of $\widetilde{\mathbb H}_{I_{17}}$. Let
\[\varpi'=(-b+\frac{1}{4})\varpi_2. \]
Let $U'=\mathrm{St}_{I_{17}}\otimes \mathbb C_{\varpi'}$ and let $ X' = \mathbb H \otimes_{\mathbb H_{I_{17}}}U' $.

Let
\[  \chi'=-\frac{1}{4}\alpha_1-(b+1)\alpha_3-(b+1)\alpha_4=\sqrt{2}(\frac{1}{8},\frac{a}{4}, -a, -\frac{b}{4} ) .
\]
Note that
\[  ||\chi'||^2=a+\frac{5}{8}, \quad ||\varpi'||^2=\frac{7}{8}-a.
\]

Let $v \in 1 \otimes U'$ be a non-zero vector. The weight of $v$ is
\[ \chi'+\varpi'=-a\alpha_1-(2a-\frac{1}{2})\alpha_2-(3a+\frac{1}{2})\alpha_3-(2a+1)\alpha_4= \sqrt{2}(-\frac{b}{2}+\frac{1}{4}, \frac{b}{2}, -a, -\frac{1}{4}).
\]
Moreover, $\chi'+\varpi'$ and $\chi_{17}$ are in the same $W$-orbit. Let $w_{17}^*$ be the unique element in $W$ such that $w_{17}^*(\chi_{17})=\chi'+\varpi'$ and $R(w_{17}^*)\cap Z(\chi_{17})=\emptyset$. The corresponding element in SageMath is $W[3101]$.

\begin{proposition} \label{prop structure of a minimal ds}
The $\mathbb H$-module $X'$ has precisely one non-tempered composition factor. Moreover, the non-tempered composition factor has the Langlands parameter $(U'{}^{\star}, I_{17})$.
\end{proposition}
\begin{proof}

    First, $Q:=L(U'{}^{\star}, I_{17})$ is an non-tempered quotient module of $X'{}^{\star}$. Thus, it suffices to show that $Q$ is the only non-tempered composition factor of $X'$.

Suppose there exists another non-tempered composition factor, say $L$, of $X'$ with Langlands parameter $(U=\widetilde{U}\otimes \mathbb C_{\varpi_L},I)$ where $\widetilde{U}$ is a tempered $\widetilde{\mathbb{H}}_I$-module and $\varpi_L \in V_I^{\vee, \bot}$. Let $\chi_L$ be a weight of $\widetilde{U}$.

Then, $\chi_L+\varpi_L$ is a weight of $X'$ and so $\chi_L+\varpi_L$ is in the $W$-orbit of $\chi_{17}$. Thus,
\begin{align} \label{eqn lenght compar langlands}
||\chi'||^2+||\varpi'||^2= || \chi_{17}||^2=||\chi_L+\varpi_L||=||\chi_L||^2+||\varpi_L||^2.
\end{align} \label{eqn chi 17 orth decomp 2}
Moreover, by tha Langlands classification, we also have:
\begin{align} \label{eqn langland inequality}
\varpi_L < -\varpi' .
\end{align}


We first observe, by Proposition \ref{prop temp induced from ds} and the real central character condition for $\chi_L$, that $||\chi_L||$ is equal to the length of some discrete series central character of a lower rank graded Hecke algebra. We shall appeal to some classification result and lengths of discrete series central characters in an appendix (See Section \ref{s appendix low rank ds}) to show that no such $L$ exists. We divide into two cases.

\begin{itemize}
\item Case (1): Suppose  $||\chi_L||\geq ||\chi'||$. Then, from Section \ref{s appendix low rank ds}, there are only two possible situations:
\begin{enumerate}
\item[(i)] $I$ is of type $H_3$ and $\widetilde{U}$ is the Steinberg module: In this case, $\chi_L=-\frac{1}{4}((10a+5)\alpha_2+(20a+8)\alpha_3+(18a+6)\alpha_4)$, and so we have 
\[ ||\chi_L||=\sqrt{\frac{48a+19}{8}}>\sqrt{\frac{3}{2}}=||\chi_{17}||.\]
\item[(ii)] $I$ is of type $H_3$ and $\chi_L=-\frac{1}{2}((2a+2)\alpha_2+(4a+3)\alpha_3+6a\alpha_4)$. We have 
\[ ||\chi_L||=\sqrt{2}>\sqrt{\frac{3}{2}}=||\chi_{17}||. \] 
\end{enumerate}
In both situations, we have 
\[ ||\chi_L||^2+||\varpi_L||^2\geq||\chi_L||^2>||\chi_{17}||^2 \]
which contradicts (\ref{eqn lenght compar langlands}).
\item Case (2): Suppose $||\chi_L||<||\chi'||$. In this case, we have $||\varpi_L||>||-\varpi'||$ by (\ref{eqn lenght compar langlands}).  On the other hand, $\varpi_L<-\varpi'$ implies $-\varpi'=\varpi_L+\gamma$ for some $\gamma$ as a positive sum of some simple roots. Thus, $\langle \varpi_L, \gamma \rangle \geq 0$. Then 
\[ ||-\varpi'||^2=||\varpi_L||^2+2\langle \varpi_L, \gamma \rangle+||\gamma||^2 > ||\varpi_L||^2
\]
giving a contradiction.
\end{itemize}

The above show that it is impossible to have a non-tempered composition factor in $X'$, other than $L(U'^{{\color{red}\star}},I)$.
\end{proof}

\begin{lemma} \label{lem temp=ds 17}
All the tempered modules of $\mathbb H$ with the central character $W\chi_{17}$ are discrete series.
\end{lemma}
\begin{proof}
Let $L$ be a tempered $\mathbb H$-module with the central character $W\chi_{17}$, but not a discrete series. By Proposition \ref{prop temp induced from ds}, $L$ is a composition factor of $\mathbb H\otimes_{\mathbb H_I}(\widetilde{U}\otimes \mathbb C_{\mu})$ for some  $I\neq\Delta$, $\widetilde{\mathbb H}_I$-discrete series $\widetilde{U}$ and $\mu \in V_I^{\bot, \vee}$ satisfying $\mathrm{Re}(\mu(\alpha))=0$ for all $\alpha \notin I$. However, this is impossible by comparing the lengths of discrete series central character in Section \ref{s appendix low rank ds}.
\end{proof}

\begin{lemma} \label{lem mini ds 17}
The module $X'$ is a minimally induced module associated to the discrete series $\mathbb H^{(\chi_{17}, J_{17})}$, where $\chi_{17}$ and $J_{17}$ are defined in Section \ref{s character 17}.
\end{lemma}

\begin{proof}
This follows from an extensive computation of the weights from the skew local region $(\chi_{17}, J_{17})$. Alternatively, we use the notations in an appendix (Section \ref{ss explicit p wchi17}), we note that
\[ \chi'+\varpi'=\frac{1}{4a+2} (2w_{17}^*(\beta_{17}^2)+(2a-1)w_{17}^*(\beta_{17}^5)+(6a+1)w_{17}^*(\beta_{17}^6)+(4a+4)w_{17}^*(\beta_{17}^9)).
\]
Moreover, the computation in Section \ref{ss explicit p wchi17} shows that $w_{17}^*$ is in the local region $\mathcal F^{(\chi_{17}, J_{17})}$.

Suppose $(I_{17},\chi'+\varpi')$ does not attain the minimum value of the set (\ref{eqn min ds}) for $\mathbb H^{(\chi_{17}, J_{17})}$. Then let $\gamma$ be a weight of $\mathbb H^{(\chi_{17}, J_{17})}$ and $I\subset \Delta$ with $|I|=|\Delta|-1$ such that the data attains the minimum value of  (\ref{eqn min ds}). Let
\[  C_1= -\frac{(\mathrm{Re}(\gamma), \varpi_{\beta})}{||\varpi_{\beta}||} , \quad C_2=-\frac{(\varpi', \varpi_{2})}{||\varpi_{2}||}=||\varpi'|| ,
\]
where $\beta$ is the unique element in $\Delta-I$. We have $C_1<C_2$ by the minimality of choice of $(\gamma, I)$.
On the other hand,
\[  \gamma=\lambda{\color{red}-}C_1\frac{\varpi_{\beta}}{||\varpi_{\beta}||}, \
\]
for some $\lambda \in V_I^{\vee}$. With the equality $||\gamma||=||\chi'+\varpi'||$, we deduce that
\[
||\lambda|| >||\chi'|| .
\]
By Theorem \ref{thm construct minimal ds} and its proof, $W_{I}\lambda$ is a discrete series central character of $\widetilde{\mathbb H}_I$. Thus, the analysis in Case 1 of the proof of Proposition \ref{prop structure of a minimal ds} shows that it is impossible. In other words, $(I_{17}, \chi'+\varpi')$ attains the minimum value of (\ref{eqn min ds}). This shows the lemma.
\end{proof}

\subsection{Stabilizer set}

\begin{lemma} \label{lem stabilizer group 17}
The set $\{w \in W^{I_{17}, I_{17}} |w(\chi'+\varpi')=\chi'+\varpi'  \}$ precisely consists of the following three elements: the identity and 
\[ us_2s_4u^{-1}, \mbox{where $u=s_2s_3s_4s_3s_4s_2s_3s_4s_1s_2s_3s_4s_2s_3s_1$,} \]
\[ vs_2s_4v^{-1}, \mbox{where $v=s_2s_3s_4s_1s_2s_3s_4s_2s_3s_4s_1s_2s_3s_4s_1s_2s_3s_4s_1s_2s_3s_1$} . \]
The three elements above are also in the reduced expressions.


\end{lemma}

\begin{proof}
 The above lemma is a direct computation on SageMath. The full stabilizer group of $\chi'+\varpi'$ is also given in an appendix (Section \ref{ss tabilizer 17}). Moreover, one can directly check (e.g. from the computations of the stabilizer group in Section \ref{ss tabilizer 17}) that 
 \[ us_2s_4u^{-1}(\alpha_i)>0, \quad vs_2s_4v^{-1}(\alpha_i)>0\]
for all $i=1,3,4$. Hence, the three elements are minimal representatives in $W/W_{I_{17}}$. Moreover, the elements are involutions and so they are also in $W^{I_{17}, I_{17}}$.
 
\end{proof}

\subsection{Restriction of $X'$}

\begin{lemma} \label{lem induced moduels h4 reduce}
$X'$ is reducible.
\end{lemma}

\begin{proof}
We have known that $\mathbb H^{(\chi_{17}, J_{17})}$ and $L(U'{}^{\star}, I_{17})$ are two composition factors of $X'$ (see, in particular, Lemma \ref{lem mini ds 17}). Thus $X'$ must have at least $2$ composition factors and this implies the lemma.
\end{proof}

\begin{corollary} \label{cor induced modules h4 17}
The natural embedding $\iota_{U'}$ from $U'$ to $X'$ splits.
\end{corollary}

\begin{proof}
This follows from Proposition \ref{prop splitting calibrated} and Lemma \ref{lem mini ds 17}.
\end{proof}

\begin{lemma} \label{lem ext vanishing}
Let $I'\subset I_{17}$. Let $Y$ be a twisted Steinerg $\mathbb H_{I'}$-module i.e. $Y|_{W_{I'}}$ is the $1$-dimensional sign module of $W_{I'}$. Let $X$ be a twisted Steinberg $\mathbb{H}_{I_{17}}$-module. Then, if 
\[  \mathrm{Hom}_{\mathbb H_{I_{17}}}(\mathbb H_{I_{17}}\otimes_{\mathbb H_{I'}} Y, X) =0 ,
\]
then, for all $i \geq 0$,
\[  \mathrm{Ext}^i_{\mathbb H_{I_{17}}}(\mathbb H_{I_{17}}\otimes_{\mathbb H_{I'}} Y, X) =0 .
\]
\end{lemma}

\begin{proof}
By Frobenius reciprocity, we have 
\[ \mathrm{Ext}^i_{\mathbb H_{I_{17}}}(\mathbb H_{I_{17}}\otimes_{\mathbb H_{I'}}Y, X) \cong \mathrm{Ext}^i_{\mathbb H_{I'}}(Y, X|_{\mathbb H_{I'}}) .\]
Since both $Y$ and $ X|_{\mathbb H_{I'}}$ are twisted Steinberg modules, the zeroness for the Hom implies that $Y$ and $X|_{\mathbb H_{I'}}$ have different central character of $\mathbb H_{I'}$. This then in turn implies that 
\[ \mathrm{Ext}^i_{\mathbb H_{I_{17}}}(\mathbb H_{I_{17}}\otimes_{\mathbb H_{I'}} Y, X) =0 
\]
for all $i\geq 0$.
\end{proof}

\begin{corollary} \label{cor coscole of induce h4}
The module $U'$ appears with multiplicity $3$ in the cosocle of $X'|_{\mathbb H_{I_{17}}}$.
\end{corollary}

\begin{proof}
Let $\iota_{U'}: U' \rightarrow X'$ be the natural embedding in Section \ref{ss induced map split}, and let $Y=\mathrm{coker}(\iota_{U'})$. Now, by Corollary \ref{cor induced modules h4 17},
\[   X'|_{\mathbb H_{I_{17}}} \cong Y \oplus U' .
\]
Thus, it suffices to show the following claim:

\noindent
{\it Claim:}
\[  \mathrm{dim}~\mathrm{Hom}_{\mathbb H_{I_{17}}}(Y, U') =2 . \]

\noindent
{\it Proof of Claim:} We shall denote the three elements in Lemma \ref{lem stabilizer group 17}, by $w_1=1, w_2, w_3$ respectively. Since
\[  l(w_2)= 32 , l(w_3)=46 , l(w_3w_2^{-1})=22
\]
$$(w_3w_2^{-1}=s_3s_4s_2s_3s_4s_1s_2s_3s_4s_2s_3s_4s_2s_3s_4s_1s_2s_3s_4s_3s_2s_1),$$
$w_2$ and $w_3$ are incomparable under the Bruhat ordering.

For $w \in W^{I_{17},I_{17}}$ and $w\neq 1$, let 
\[  Z_w =\mathbb H_{I_{17}}\otimes_{\mathbb H_{I_{17}\cap w(I_{17})}} ((U')|_{\mathbb H_{w^{-1}(I_{17})\cap I_{17}}})^{w^{-1}} .
\]
The geometric lemma in Section \ref{ss geo lem} implies that $Y$ has a submodule $Y'$ such that
\begin{enumerate}
\item $Y'$ admits a short exact sequence:
\[  0 \rightarrow Y'' \rightarrow Y' \rightarrow Z_{w_2}\oplus Z_{w_3} \rightarrow 0,
\]
where $Y''$ has a filtration with layers for all $Z_w$ with 
$w \not\geq w_2$ and $w \not\geq w_3$ (and $w\neq 1$);
\item $Y/Y'$ admits a filtration with layers for all $Z_w$ with either $w >w_2$ or $w > w_3$.
\end{enumerate}
It is clear by Lemmas \ref{lem stabilizer group 17} and \ref{lem ext vanishing} that 
\begin{align} \label{eqn 2 dim submodule}
 \mathrm{dim}~ \mathrm{Hom}_{\mathbb H_{I_{17}}}(Y', U') =2 .
\end{align}

By Lemma \ref{lem stabilizer group 17}, for $w \in W^{I_{17},I_{17}}$ satisfying $w>w_2$ or $w>w_3$,
\[  \mathrm{dim}~\mathrm{Hom}_{\mathbb H_{I_{17}}}(Z_w, U')= 0
\]
and so, by Lemma \ref{lem ext vanishing}, $\mathrm{dim}~\mathrm{Ext}^i_{\mathbb H_{I_{17}}}(Z_w, U')=0$ for all $i$. Thus, a standard long exact sequence argument gives that for all $i$,
\begin{align} \label{eqn zero dim for quotient}
  \mathrm{dim}~\mathrm{Ext}^i_{\mathbb H_{I_{17}}}(Y/Y', U') =0 .
\end{align}
Thus, the long exact sequence:
\[  0 \rightarrow \mathrm{Hom}_{\mathbb H_{I_{17}}}(Y/Y', U') \rightarrow \mathrm{Hom}_{\mathbb H_{I_{17}}}(Y, U') \rightarrow \mathrm{Hom}_{\mathbb H_{I_{17}}}(Y', U') \rightarrow \mathrm{Ext}^1_{\mathbb H_{I_{17}}}(Y/Y', U') .
\]
Since the first and last spaces have zero dimensional by (\ref{eqn zero dim for quotient}), 
\[  \mathrm{Hom}_{\mathbb H_{I_{17}}}(Y, U') \cong \mathrm{Hom}_{\mathbb H_{I_{17}}}(Y', U') 
\]
Now the claim follows from (\ref{eqn 2 dim submodule}), as desired.
\end{proof}

\begin{corollary} \label{cor after cosocle}
$\mathrm{dim}~\mathrm{Hom}_{\mathbb H}(X', X'{}^{\star})=3$.
\end{corollary}

\begin{proof}
We have
\[ \mathrm{Hom}_{\mathbb H}(X', X'{}^{\star}) \cong \mathrm{Hom}_{\mathbb H}(X', \mathrm{Ind}_{\mathbb H_{I_{17}}}^{\mathbb H}(U'{}^{\star})) \cong \mathrm{Hom}_{\mathbb H_{I_{17}}}(X'|_{\mathbb H_{I_{17}}}, U') ,
\]
where the second isomorphism follows from Proposition \ref{prop second adjoint} and Remark \ref{rmk action on twisted maps}. Now the corollary follows from Corollary \ref{cor coscole of induce h4}.
\end{proof}

\subsection{Existence of discrete series for $\chi_{17}$}

\begin{proposition} \label{prop existence ds 17}
There exist 3 mutually non-isomorphic discrete series in the Jordan-H\"older series of $X'$. 
\end{proposition}

\begin{proof}
Proposition \ref{prop structure of a minimal ds} checks the condition (2) in Lemma \ref{lem general strategy}. Lemmas \ref{lem temp=ds 17} and \ref{lem induced moduels h4 reduce} check the conditions (3) and (1) in  Lemma  \ref{lem general strategy} respectively. With Corollary \ref{cor after cosocle}, Lemma \ref{lem general strategy} implies that there exist $3$ mutually non-isomorphic discrete series in $X'$.
\end{proof}

\section{Existence of discrete series at $W\chi_{16}$} \label{s exist ds 16}

\subsection{Minimally induced module for $\mathbb H^{(\chi_{16},J_{16})}$} \label{ss skew 16 ds}

Let $I_{16}=\left\{ \alpha_2,\alpha_3, \alpha_4\right\}$. Let
\[  \varpi''=-b\varpi_{1}.
\]
Let $\widetilde{M}$ be the irreducible  discrete series of $ \widetilde{\mathbb H}_{I_{16}}$ with dimension equal to $4$. Moreover, $\widetilde{M}$ precisely has the following four weights:
\[  \chi'':=-(1+a)\alpha_2-(2a+3/2)\alpha_3-3a\alpha_4 =-\sqrt{2} (\frac{3-2a}{4}, \frac{a+1}{2}, \frac{1}{4},0 ) {\color{red},}
\]
\[  s_4\chi''= -(1+a)\alpha_2 -(2a+3/2)\alpha_3-(1+2a)\alpha_4 =-\frac{\sqrt{2}}{2}(1, 2a, 2b, 0),
\]
\[  s_3s_4\chi''= -(1+a)\alpha_2 -(1/2+3a)\alpha_3-(1+2a)\alpha_4=-\frac{\sqrt{2}}{2}(1+b, 3a-1, \frac{1}{2}, 0),
\]
\[   s_2s_3s_4\chi''= -(2a-1/2)\alpha_2 -(1/8+3a)\alpha_3-(1+2a)\alpha_4=-\frac{\sqrt{2}}{2}(2b, 1, 2a, 0).\]
Let $w_{16}^*$ be the element in $W$ such that $w_{16}^*(\chi_{16})=\chi''+\varpi''$ and $R(w_{16}^*)\cap Z(\chi_{16})=\emptyset$. The corresponding element in SageMath is $W[1040]$. 

Note that
\[  ||\chi''||^2=2, \quad ||\varpi''||^2 = \frac{1}{2}.
\]
Let $U''=\widetilde{M}\otimes \mathbb C_{\varpi''}$ and let  $X''=\mathbb H\otimes_{\mathbb H_{I_{16}}}U''$.




\subsection{Structure of $X''$}

\begin{proposition} \label{prop 16 minimal induce ds}
There is precisely one non-temepred composition factor in $X''$. Moreover, the Langlands parameter of the non-tempered composition factor is $(U''^{\star}, I_{16})$.
\end{proposition}
\begin{proof}
The proof is similar to the one of Proposition \ref{prop structure of a minimal ds}. We need to do some computation in this particular case to show there is only one non-tempered composition factor. 

Note that $X''{}^{\star}$ is a standard module of $\mathbb H$ (see Theorem \ref{thm construct minimal ds}(3)), and the Langlands quotient is $L(U''^{\star}, I_{16})$. Suppose $L$ is another non-tempered composition factor in $X''$. Similar to the proof of Proposition \ref{prop structure of a minimal ds}, let $(\widetilde{U}\otimes \mathbb C_{\varpi_L}, I)$ be the Langlands parameter of $L$, where $\widetilde{U}$ is a tempered $\widetilde{\mathbb H}_I$-module and $\varpi_L \in V_I^{\bot, \vee}$. Let $\chi_L$ be a weight of $\widetilde{U}$.


We again consider two cases:

\begin{itemize}
\item Case (1):  $||\chi_L||\geq ||\chi''||$. By the classification of discrete series central character in Section \ref{s appendix low rank ds}, there is only one possible situation:
\begin{itemize}
\item $I$ is of type $H_3$ and $\widetilde{U}$ is the Steinberg module. But, we then have $||\chi_L||=\sqrt{\frac{48a+19}{8}}>\sqrt{\frac{5}{2}}=||\chi_{16}||$. Then, arguing as in the Case (1) of the proof of Proposition \ref{prop structure of a minimal ds}, one shows this case is impossible.
\end{itemize}
\item Case (2): If $||\chi_L||<||\chi''||$. Then one applies the same argument in the Case (2) of the proof of Proposition \ref{prop structure of a minimal ds} to show it is impossible.
\end{itemize}
\end{proof}

\begin{lemma}
All the tempered modules of $\mathbb H$ with the central character $W\chi_{16}$ are discrete series.
\end{lemma}

\begin{proof}
This is similar to Lemma \ref{lem temp=ds 17}. One compares the lengths of the discrete series central characters in Section \ref{s appendix low rank ds} with $||\chi_{16}||$.
\end{proof}

\begin{lemma} \label{lem mini ds 16}
The module $X''$ is a minimally induced module associated to the discrete series $\mathbb H^{(\chi_{16}, J_{16})}$, where $\chi_{16}$ and $J_{16}$ are defined in Section \ref{s character 16}.
\end{lemma}

\begin{proof}
It is similar to the one in Lemma \ref{lem mini ds 17}. Indeed, {\color{red}}
\[   w_{16}^*\chi_{16} = w_{16}^*\beta_{16}^2+w_{16}^*\beta_{16}^3+\frac{3}{2}w_{16}^*\beta_{16}^6+\frac{3}{2}w_{16}^*\beta_{16}^7 
\]
and so $w_{16}^*$ is in $\mathcal F^{(\chi_{16}, J_{16})}$.  Now, the analysis in Proposition \ref{prop 16 minimal induce ds} also gives that $X''$ is the minimally induced module associated to the discrete series $\mathbb H^{(\chi_{16}, J_{16})}$. 
\end{proof}

\subsection{Stabilizer subgroup}

\begin{proposition} \label{prop stabilizer chi16} 
Let $W^{I_{16},I_{16}}$ be the set of minimal representatives in $W_{I_{16}}\setminus W/W_{I_{16}}$. Then $\{w \in W^{I_{16},I_{16}} : w(\chi''+\varpi'')=\chi''+\varpi'' \}$ consists of precisely the following two elements: the identity and 
\[ xs_2s_4x^{-1}, \mbox{where $x=s_1s_2s_3s_4s_3s_4s_2s_3s_4s_1s_2s_3$.} \]


\end{proposition}

\begin{proof}
This is done by a direct check in SageMath. On the other hand, one checks (e.g. use computations in Section \ref{ss stabilizer chi16 others}) that:
\begin{align} \label{eqn u^* action on I}
     xs_2s_4x^{-1}\alpha_2=\alpha_3, \quad xs_2s_4x^{-1}\alpha_3=\alpha_2, \quad xs_2s_4x^{-1}\alpha_4=\alpha_1+\alpha_2+(2a+1)\alpha_3+(2a+1)\alpha_4,
\end{align}
and so $xs_2s_4x^{-1}$ is a minimal representative in $W/W_{I_{16}}$. Since $xs_2s_4x^{-1}$ is an involution, $xs_2s_4x^{-1}$ is in $W^{I_{16},I_{16}}$.
\end{proof}

\subsection{Restriction on $U''$} \label{ss str rest U''}
Let $u^{*}$ be the unique non-trivial element in $\{w \in W^{I_{16},I_{16}} : w(\chi''+\varpi'')=\chi''+\varpi'' \}$. 

 It follows from the equations (\ref{eqn u^* action on I}) that
\[ I':= I_{16} \cap u^*{}^{-1}(I_{16}) =\left\{ \alpha_2, \alpha_3 \right\} .\]
Note that $U''|_{\mathbb H_{I'}}$ is a direct summand of two irreducible $\mathbb H_{I'}$-modules, with one to be one-dimensional twisted Steinberg module with the weight $u^{*}\chi_{16}$ and another to be three-dimensional of the weights $s_4w_{16}^{*}\chi_{16}, s_3s_4w_{16}^{*}\chi_{16}, s_2s_3s_4w_{16}^{*}\chi_{16}$.

\subsection{Some vanishing Ext-groups}

A complication for $X''$ is that $U''$ is not a twisted Steinberg module, and so we cannot use a simple trick like Lemma \ref{lem ext vanishing}. We show other useful lemma:

\begin{lemma} \label{lem ext vanish 16}
Let $w \in W^{I_{16},I_{16}}$. Let $I_w=I_{16} \cap w(I_{16})$. Let $\gamma$ be a weight of $U''$. Then, if $\mathrm{Ext}^i_{\mathbb H_{I_w}}((U'')^{w^{-1}}, U'')\neq 0$ for some $i$, then there exist $w_1, w_2 \in W_{I_{16}}$ such that $w_2ww_1(\gamma)=\gamma$.
\end{lemma}

\begin{proof}
The non-vanishing Ext-groups imply that there exists a weight $\gamma'$ of $U''$ such that $w(\gamma')$ and $\gamma$ are in the same $W_{I_{w}}$-orbit and so also in the same $W_{I_{16}}$-orbit. Thus, $w(\gamma')$, $\gamma'$ and $\gamma$ are in the same $W_{I_{16}}$-orbit. Now, one can find $w_1, w_2 \in W_{I_{16}}$ such that $w_1(\gamma)=\gamma'$ and $w_2w(\gamma')=\gamma$. This then implies the lemma.
\end{proof}

\begin{lemma} \label{lem ext vanishing for larger layer}
 Let $w' \in W^{I_{16},I_{16}}$ such that $w'>u^*$. Let $I_{w'}=I_{16}\cap w'(I_{16})$. Then, for all $i \geq 0$,
\[   \mathrm{Ext}^i_{\mathbb H_{I_{16}}}(\mathbb H_{I_{16}}\otimes_{\mathbb H_{I_{w'}}} (U'')^{w'{}^{-1}}, U'') =0. 
\]
\end{lemma}

\begin{proof}
Suppose the Ext-group is non-zero. By Lemma \ref{lem ext vanish 16}, there exist $x, y \in W_{I_{16}}$ such that $yw'x(\chi''+\varpi'')=\chi''+\varpi''$. But then, $l(yw'x)>l(u^*)$. However, the complete list of the stabilizer group of $\chi''+\varpi''$ in Section \ref{ss stabilizer 16} shows that it is impossible.

\end{proof}

\subsection{Existence of discrete series}

\begin{lemma} \label{lem dim greater than 1}
Recall that $u^{*}$ and $I'$ are defined in Section \ref{ss str rest U''}.
\[ \mathrm{dim}~\mathrm{Hom}_{\mathbb H_{16}}(\mathbb H_{I_{16}} \otimes_{\mathbb H_{I'}}(U'')^{u^*{}^{-1}}, U'')  \geq 1 .
\]
\end{lemma}

\begin{proof}
The structure of $U''|_{\mathbb H_{I'}}$ is described in Section \ref{ss str rest U''}, and in particular, it has a one-dimensional direct summand, denoted by $Y$. By Proposition \ref{prop stabilizer chi16}, the twist $Y^{u^{\star}{}^{-1}}$ is isomorphic to $Y$, as $\mathbb H_{I'}$-modules, which then contributes to one dimension in $\mathrm{Hom}_{\mathbb H_{I'}}((U'')^{u^{\star}{}^{-1}}, U'')$. Now the lemma follows from Frobenius reciprocity.
\end{proof}

\begin{lemma} \label{lem multiplicity at least 2 socle 16}
The module $U''$ appears with multiplicity at least $2$ in the cosocle of $X''|_{\mathbb H_{I_{16}}}$.
\end{lemma}

\begin{proof}

Let $\iota_{U''}: U'' \rightarrow X''$ be the natural embedding in Section \ref{ss induced map split}. Then, we have a splitting:
\begin{align} \label{eqn split decompose x''}
X''|_{\mathbb H_{I_{16}}} \cong U'' \oplus Y''  .
\end{align}

Now, for $w \in W^{I_{16},I_{16}}$ and $w\neq 1$, let
\[ Z_w = \mathbb H_{I_{16}} \otimes_{\mathbb H_{I_{16}\cap w(I_{16})}} (U''|_{\mathbb H_{w^{-1}(I_{16})\cap I_{16}}})^{w^{-1}}  .
\]
By the geometric lemma, $Y''$ admits a filtration whose geometric lemma admits successive subquotients of the form $Z_w$. 

The lemma will follow from (\ref{eqn split decompose x''}) and the following claim.

\noindent
{\it Claim:} 
\[  \mathrm{dim}~\mathrm{Hom}_{\mathbb H_{I_{16}}}(Y'', U'') \geq 1
\]

\noindent
{\it Proof of Claim:} 
Let $\widetilde{Y}$ be the quotient of $Y''$ that precisely contains all layers $Z_{w}$ with $w \geq u^*$. Then it suffices to show 
\[  \mathrm{dim} \mathrm{Hom}_{\mathbb H_{I_{16}}}(\widetilde{Y}, U'') \geq1  .
\]
Now, $\widetilde{Y}$ admits a short exact sequence
\[  0 \rightarrow Z_{u^*} \rightarrow \widetilde{Y} \rightarrow \widetilde{Z} \rightarrow 0 ,
\]
where $\widetilde{Z}$ admits filtrations with successive subquotients isomorphic to $Z_w$ with $w>u^*$. Then, by Lemma \ref{lem ext vanishing for larger layer}, 
\[  \mathrm{Ext}^i_{\mathbb H_{I_{16}}}(\widetilde{Z}, U'')=0
\]
for all $i$. Then, a standard argument of homological algebra implies that 
\[  \mathrm{Hom}_{\mathbb H_{I_{16}}}(Z_{u^*}, U'')\cong \mathrm{Hom}_{\mathbb H_{I_{16}}}(\widetilde{Y}, U'') .
\]
Now, it follows from Lemma \ref{lem dim greater than 1} that the former space has dimension at least one and so the latter space has dimension at least one. This proves the claim. \\

\end{proof}

\begin{proposition} \label{prop existence ds 16}
There exist at least two non-isomorphic discrete series in $X''$.
\end{proposition}

\begin{proof}
Let 
\[  p=\mathrm{dim}~\mathrm{Hom}_{\mathbb H}(X'', X''{}^{\star})=\mathrm{dim}~\mathrm{Hom}_{\mathbb H_{I_{16}}}(X''|_{\mathbb H_{I_{16}}} , U'' ) ,
\]
where the second equality follows from the second adjointness and Remark \ref{rmk action on twisted maps}. 

By Lemma \ref{lem multiplicity at least 2 socle 16}, $p \geq 2$. Combining with Lemma \ref{lem general strategy}, the cosocle of $X''$ contains at least $2$ discrete series. By Lemma \ref{lem mini ds 16}, it remains to show that at least one such discrete series is not isomorphic to $\mathbb H^{(\chi_{16}, J_{16})}$.

Suppose all discrete series in the cosocle of $X''$ are isomorphic to $\mathbb H^{(\chi_{16}, J_{16})}$, and let $r$ be the multiplicity of $\mathbb H^{(\chi_{16}, J_{16})}$ in the cosocle of $X''$. Since
\[  \mathrm{Hom}_{\mathbb H}(\mathrm{cosoc}(X''), \mathrm{cosoc}(X'')) \cong \mathrm{Hom}_{\mathbb H}(X'', X''{}^{\star}),
\]
(by Lemma \ref{lem general strategy}(i)), we have 
\[  \mathrm{dim}~\mathrm{Hom}_{\mathbb H}(X'', X''{}^{\star}) =r^2.
\]
On the other hand, 
\begin{align} \label{eqn fb ds 16}
\mathrm{Hom}_{\mathbb H}(X'', X''{}^{\star}) \cong \mathrm{Hom}_{\mathbb H}(X'', \mathrm{cosoc}(X'')) &\cong \mathrm{Hom}_{\mathbb H}(X'', (\mathbb H^{(\chi_{16}, J_{16})})^{\oplus r}) \\ 
& \cong \mathrm{Hom}_{\mathbb H}(X'', \mathbb H^{(\chi_{16}, J_{16})})^{\oplus r} \\
&\cong \mathrm{Hom}_{\mathbb H_{I_{16}}}(U'', \mathbb H^{(\chi_{16}, J_{16})})^{\oplus r} ,
\end{align}
where the first isomorphism follows from a standard homological argument and the last isomorphism follows from Frobenius reciprocity.

By the structure of a calibrated module in Corollary \ref{cor property of cal rep}(2), 
\[  \mathrm{dim}~\mathrm{Hom}_{\mathbb H_{I_{16}}}(U'',\mathbb H^{(\chi_{16}, J_{16})}) \leq 1\]
and so is equal to $1$ by Lemma \ref{lem mini ds 16}, 
and so (\ref{eqn fb ds 16}) gives that
\[ \mathrm{dim}~\mathrm{Hom}_{\mathbb H}(X'', X''{}^{\star})=r .\]
Since $r\neq r^2$ (as $r\geq 2$), we obtain a contradiction.
\end{proof}

\part{Classification of discrete series and some more properties}

\section{Classification theorem of discrete series for $\mathbb H$ of type $H_4$} \label{s classify ds}

\subsection{Upper bound on the number of discrete series}

\begin{proposition} \label{prop upper bound ds}
There are at most $20$ isomorphism classes of discrete series for the graded Hecke algebra $\mathbb H$ of type $H_4$.
\end{proposition}

\begin{proof}
It follows from \cite[Corollary 8.2]{Ch16} that the number of ismorphism classes of discrete series for $\mathbb H$ is less than or equal to the number of elliptic conjugacy classes in the Weyl group of type $H_4$. Now, the proposition follows from that the number of elliptic conjugacy classes  in the Weyl group of type $H_4$ is $20$ (see \cite[Remark 5.3]{Ch13}).
\end{proof}

\subsection{Classification Theorem}
\begin{theorem} \label{thm classification}
Let $\mathbb H$ be the graded Hecke algebra associated to the root system of type $H_4$ (see Section \ref{ss gha}). Then,
\begin{enumerate}
\item There are precisely 20 isomorphism classes of discrete series of $\mathbb H$.
\item  The set of discrete series central character of $\mathbb H$ is precisely the set of Heckman-Opdam central characters of $\mathbb H$.
\item There are precisely one isomorphism class of discrete series of $\mathbb H$ at the central character $W\chi_i$, where $i=1, \ldots, 15$, and precisely two isomorphism classes of discrete series at the central character $W\chi_{16}$ and precisely three isomorphism classes of discrete series at the central character $W\chi_{17}$.
\end{enumerate}
\end{theorem}

\begin{proof}
Proposition \ref{prop upper bound ds} gives an upper bound on the number of isomorphism classes of discrete series. We constructed 17 discrete series in Theorems \ref{thm regular ds}, \ref{thm ds z=1}, \ref{thm chi16 ds}, \ref{thm 17 skew}, and show existence of three additional discrete series in Propositions \ref{prop existence ds 17} and \ref{prop existence ds 16}. These give the three assertions.
\end{proof}

\subsection{Combine with previous cases}

\begin{corollary} \label{cor dscc-hocc}
Let $R$ be a root system and let $c$ be a positive parameter function. Let $\mathbb H$ be the graded Hecke algebra associated to $R$ and the paremeter function $c$. Then, the set of discrete series central characters of $\mathbb H$ coincides with the set of Heckman-Opdam central characters of $\mathbb H$.
\end{corollary}

\begin{proof}
When $R$ is of crystallographic type, it follows from \cite[Theorem 8.3]{KL87} for crystallographic types of equal parameter (and the Lusztig reduction to graded Hecke algebra) and from \cite[Theorem 7.1]{OS10} (also see \cite[, Theorem 3.29]{Op04}) for other unequal parameter cases. When $R$ is of non-crystallographic type, it follows from \cite{KR02} for type $I_2(n)$, \cite{Kri99} for type $H_3$ and Theorem \ref{thm classification} for type $H_4$.
\end{proof}

\section{Anti-spherical discrete series} \label{s sph ds}

In this section, we shall first consider general graded Hecke algebra $\mathbb H$ defined in Definition \ref{def graded hecke alg} . We assume that the parameter function $c_{\alpha}> 0$ for all $\alpha \in \Delta$.




\subsection{Anti-dominant character and local region}

\begin{definition}
We say $\chi \in V^{\vee}$ is \textit{anti-dominant} if $\chi(\alpha)\leq0$ for all $\alpha \in R^+$.
\end{definition}

We first recall the following standard result, see e.g. \cite{Hu90}, or see \cite[Lemma 1.2]{DL18} for a more general statement:

\begin{proposition}
Let $\chi \in V^{\vee}$. Then there exists a unique anti-dominant character in the $W$-orbit $W\chi$.     
\end{proposition}



\begin{proposition} \label{prop local region antid c}
Let $\chi \in V^{\vee}$. Let $\chi^{*}$ be the anti-dominant character in $W\chi$ and let $ w^*$ be the unique element in $W$ such that $R(w^*)\cap Z(\chi)=\emptyset$ and $\chi^{*}=w^*\chi$. Then $w^*\in\mathcal{F}^{(\chi,P(\chi))}$ and $w^*(\beta)\in R^-$ for all $\beta \in P(\chi)$. 
\end{proposition}
\begin{proof} 
For any $\beta \in P(\chi)$, 
\[ \langle w^*\beta, w^*\chi \rangle =\langle \beta, \chi \rangle =c_{\beta} . \]
Since $w^*\chi$ is anti-dominant and $k_{\beta}>0$, we must have that $w^*\beta$ is in $R^-$. This implies that $w^* \in \mathcal F^{(\chi, P(\chi))}$.
 \end{proof}




\subsection{Some properties of anti-sphericity}

\begin{definition} An $\mathbb{H}$-module $M$ is \textit{anti-spherical} if it contains the sign representation of $W$ when restricted to $W$.
\end{definition}

\begin{lemma} 
Let $X$ be a finite-dimensional $\mathbb{H}$-module.  Then $X$ is anti-spherical if and only if $X^{\star}$ is anti-spherical.
\end{lemma}

\begin{proof}
It follows from Definition \ref{def star operation} that $X^{\star}|_W$ is the dual $W$-representation of $X|_W$. This implies that $X$ contains the sign $W$-representation if and only if $X^{\star}$ contains the sign $W$-representation. Now the lemma follows.
\end{proof}

\begin{lemma} \label{lem induce spherical}
Let $X$ be a finite-dimensional $\mathbb H_I$-module. Then $\mathbb H\otimes_{\mathbb H_I}X$ is anti-spherical if and only if $X$ is anti-spherical.

\end{lemma}

\begin{proof}
Note that, as $W$-representations, $(\mathbb H\otimes_{\mathbb H_I}X)|_W \cong  \mathbb C[W]\otimes_{\mathbb C[W_I]}(X|_{W_I})$. On the other hand, by Frobenius reciprocity,
\[  \mathrm{Hom}_{W}(\mathbb C[W]\otimes_{\mathbb C[W_I]}(X|_{W_I}), \mathrm{sign}) \cong \mathrm{Hom}_{W_I}(X|_{W_I}, \mathrm{sign}) .\]
Now the lemma follows from the above two isomorphisms.
\end{proof}

\subsection{Iwahori-Matsumoto involution}
\begin{definition}
    The \textit{Iwahori-Matsumoto involution} $\mathrm{IM}$ is a linear involution $\mathrm{IM}:\mathbb{H}_I\rightarrow\mathbb{H}_I$ characterized by $$\mathrm{IM}(v)=-v, \quad\forall v\in V,$$
    $$\mathrm{IM}(t_{w})=(-1)^{l(w)}t_{w},\quad \forall w\in W_I.$$

    For an $\mathbb{H}_I$-module $X$, $\mathrm{IM}(X)$ is an $\mathbb{H}_I$-module isomorphic to $X$ as vector spaces with $\mathbb{H}_I$-action given by: $$\pi_{\mathrm{IM}(X)}(h)x=\pi_{X}(\mathrm{IM}(h))x,\forall x\in X,$$ where $\pi_{\mathrm{IM}(X)}$ (resp.$\pi_{X})$ are the maps defining the $\mathbb{H}_I$-action on $\mathrm{IM}(X)$ (resp. $X$). 
\end{definition}

The following lemma is well-known and the proof is straightforward:

\begin{lemma} \label{lem trans im}
Let $I\subset \Delta$. Let $U$ be an $\mathbb H_I$-module. Then 
\[  \mathrm{IM}(\mathbb H\otimes_{\mathbb H_I}U) \cong \mathbb H\otimes_{\mathbb H_I}\mathrm{IM}(U) .
\]
\end{lemma}

\begin{proposition} \label{prop im of ds is non sph}
Suppose $R$ is of rank at least $1$. For any essentially discrete series $X$ of $\mathbb{H}$, $\mathrm{IM}(X)$ is not anti-spherical.
\end{proposition}
\begin{proof}
Recall that $\Delta$ is the set of simple roots of $R$. When $|\Delta|=1$, it is a standard straightforward verification that $X|_W$ is the $1$-dimensional module with only the sign $W$-representation. Hence, $\mathrm{IM}(X)|_W$ is the trivial $W$-representation. This shows the proposition in this case.

We now proceed by induction. Assume $|\Delta|=k \geq 2$. By Theorem \ref{thm construct minimal ds}, there exists a subset $I\subset \Delta $ such that
\begin{enumerate}
\item $|I|=k-1$; 
\item there exists a discrete series $ds$ of $\widetilde{\mathbb H}_I$ and a character $\varpi \in V_I^{\bot, \vee}$ such that $\mathbb H\otimes_{\mathbb H_I} (ds\otimes \mathbb C_{\varpi})$ surjects onto $X$.
\end{enumerate}

Now, by the induction assumption, $\mathrm{IM}(ds\otimes \mathbb C_{\varpi})$ is not anti-spherical and so, by Lemmas \ref{lem induce spherical} and \ref{lem trans im}, $\mathrm{IM}(\mathbb H\otimes_{\mathbb H_I}(ds\otimes \mathbb C_{\varpi}))$ is not anti-spherical. 

On the other hand, the surjection from $\mathbb H\otimes_{\mathbb H_I}(ds \otimes \mathbb C_{\varpi})$ onto $X$ induces a surjection from $\mathrm{IM}(\mathbb H\otimes_{\mathbb H_I}(ds\otimes \mathbb C_{\varpi})$ onto $\mathrm{IM}(X)$. Since $\mathrm{IM}(\mathbb H\otimes_{\mathbb H_I}(ds\otimes \mathbb C_{\varpi}))$ is not anti-spherical, we then also have that $\mathrm{IM}(X)$ is not anti-spherical, as desired.


\end{proof}

\subsection{Principal series}
\begin{definition} 
For $\chi \in V^{\vee}$, the \textit{principal series} $I(\chi)$ is the $\mathbb{H}$-module defined by $$I(\chi)=\mathrm{Ind}_{S(V)}^{\mathbb{H}}\mathbb{C}_{\chi}=\mathbb{H}\otimes_{S(V)}\mathbb{C}_{\chi}$$ with $\mathbb{H}$ acting by left multiplication.
\end{definition}




\begin{lemma} \label{lem central char unique sp}
Let $\chi \in V^{\vee}$. Then there is precisely one irreducible anti-spherical $\mathbb H$-module with central character $W\chi$.
\end{lemma}

\begin{proof}
Note that the principal series $I(\chi)$ is isomorphic to the regular representation $\mathbb C[W]$, as $W$-representations. Thus the sign representation appears with multiplicity one in $I(\chi)$. Thus, there is precisely one simple composition factor in $I(\chi)$ which is anti-spherical.

On the other hand, it follows from \cite[Proposition 2.8(a) and (c)]{KR02} (also see \cite[Theorem 2.10]{Ka82}) that each irreducible module with the central character $W\chi$ is a factor in the Jordan-H\"older series of $I(\chi)$. Combining with the previous paragraph, we obtain the lemma.
\end{proof}

\subsection{Non-antispherical composition factors of $I(\chi)$}

\begin{lemma} \label{lem sph implies anti weight}
Let $X$ be an irreducible $\mathbb{H}$-module. If $X$ does not contain an anti-dominant weight, then $X$ is not anti-spherical.   
\end{lemma}
\begin{proof}
When $|\Delta|=1$, this case is well-known and straightforward. We remark that one needs to use the parameter function $c$ to be positive here.

Now, we consider $|\Delta|=k\geq 2$. Then, by the Langlands classification and Proposition \ref{prop temp induced from ds}, $\mathrm{IM}(X)$ is an irreducible subquotient of $\mathbb H \otimes_{\mathbb{H}_I}(ds\otimes\mathbb{C}_{\nu})$, where $I \subset \Delta$, $ds$ is a discrete series for $\widetilde{\mathbb{H}}_I$ and $\nu \in V_I^{\bot, \vee}$ which satisfies $\nu(\alpha)\geq0$, for all $\alpha \notin I$. 

We consider the following two cases:
\begin{itemize}
\item Case 1: $|I|= 0$.  Then $\mathrm{IM}(X)$ contains a dominant weight and so $X$ contains an anti-dominant weight. This gives a contradiction, and so this case indeed cannot happen.
\item Case 2: $0<|I| \leq k$. Then $X$ is an irreducible quotient of $\mathrm{IM}(\mathbb H\otimes_{\mathbb H_I}(ds\otimes \mathbb C_{\nu}))$ and so, by Lemma  \ref{lem trans im}, is an irreducible quotient of $\mathbb{H} \otimes_{\mathbb H_I}(\mathrm{IM}(ds)\otimes\mathbb{C}_{-\nu})$. 

Since $\mathrm{IM}(ds)$ is not anti-spherical by Proposition \ref{prop im of ds is non sph}, $\mathbb H\otimes_{\mathbb H_I}(\mathrm{IM}(ds)\otimes\mathbb{C}_{-\nu})$ is not anti-spherical according to Lemma \ref{lem induce spherical}, and so its irreducible quotient $X$ cannot be anti-spherical.
\end{itemize}
\end{proof}

\begin{lemma} \label{lem unique anti wt module}
Let $\chi \in V^{\vee}$. There exists a unique irreducible $\mathbb H$-module which has the central character $W\chi$ and has an anti-dominant weight.
\end{lemma}

\begin{proof}
By applying the (covariant and exact) Iwahori-Matsumoto involution functor, we can prove the statement for a dominant weight instead of an anti-dominant weight. Let $\chi^*$ be the dominant character in $W\chi$. By Frobenius reciprocity, it suffices to show that $I(\chi^*)$ has a unique irreducible quotient. To this end, we write $\chi^*=\sum_{\alpha \in \Delta} d_{\alpha} \varpi_{\alpha}^{\vee}$ for some $d_{\alpha}$ satisfying $\mathrm{Re}(d_{\alpha})\geq 0$.  Let $I'=\left\{ \alpha \in \Delta : \mathrm{Re}(d_{\alpha})=0 \right\}$. By the transitivity of inductions, we have:
\[ I(\chi^*)= \mathbb H\otimes_{S(V)}\mathbb C_{\chi^*} \cong \mathbb H\otimes_{\mathbb H_{I'}}(\mathbb H_{I'}\otimes_{S(V)}\mathbb C_{\chi^*}) .
\]

Now note $\mathbb H_{I'}\otimes_{S(V)}\mathbb C_{\chi^*}$ is irreducible by \cite[Proposition 2.9]{KR02} (we need the assumption that $c_{\alpha}\neq 0$ for all $\alpha$ here). Moreover,
\begin{align*} 
    \mathbb H_{I'}\otimes_{S(V_{I'})}\mathbb C_{\chi^*}\cong (\widetilde{\mathbb H}_{I'}\otimes_{S(V_{I'})}\mathbb C_{\chi_1^*})\otimes \mathbb C_{\chi_2^*} , 
    \end{align*}
where $\chi_1^*=\sum_{\alpha \in I'}d_{\alpha}\varpi_{\alpha}^{\vee}$ and $\chi_2'=\sum_{\alpha\notin I'}d_{\alpha} \varpi_{\alpha}^{\vee}$. Hence, $(\mathbb H_{I'}\otimes_{S(V)}\mathbb C_{\chi^*}, I')$ is a Langlands parameter. Thus, $ \mathbb H\otimes_{S(V)}\mathbb C_{\chi^*}$ is a standard module and so has a unique irreducible quotient, as desired.
\end{proof}

\begin{theorem} \label{thm spherical anti-dominant}
Let $\mathbb H$ be a graded Hecke algebra of a positive parameter function defined in Definition \ref{def graded hecke alg}. Then an irreducible $\mathbb H$-module $X$ is anti-spherical if and only if $X$ has the anti-dominant weight.     
\end{theorem}
\begin{proof}
The only if direction is Lemma \ref{lem sph implies anti weight}. For the if direction, by Lemma \ref{lem unique anti wt module}, for each central character $W\chi$, there is only one irreducible module with the central character $W\chi$ that has an anti-dominant weight. Now the if part follows from the only if direction and Lemma \ref{lem central char unique sp}.
\end{proof} 

We also remark that Theorem \ref{thm spherical anti-dominant} and the proof of Lemma \ref{lem unique anti wt module} essentially determine the Iwahori-Matsumoto involution of an anti-spherical discrete series:


\subsection{Specializing to $\mathbb H_4$ case}

We now obtain the anti-sphericity of calibrated discrete series of $\mathbb H$.

\begin{corollary}
For $i=1,2,4,5,7,8,9, 12$, $\mathbb H^{(\chi_i, J_i)}$ are anti-spherical discrete series, and for $i=3,6,10,11,13,14,15,16,17$, $\mathbb H^{(\chi_i, J_i)}$ are not anti-spherical.    
\end{corollary}

\section{Branching laws for some discrete series} \label{s ext branching ds}

\subsection{Anti-spherical regular discrete series}

An irreducible module $X$ of $\mathbb H$ is said to be \textit{regular} if $X$ has a regular central character.

\begin{proposition} \label{thm vanishing ext branching}
Let $R$ be a root system and let $\Delta$ be a fixed set of simple roots of $R$. Let $\mathbb H$ be the graded Hecke algebra associated to $R$ and a positive parameter function $c$. Let $X$ be a regular anti-spherical discrete series of $\mathbb H$. Let $I \subset \Delta$. Let $Y$ be an anti-spherical discrete series of $\widetilde{\mathbb H}_I$. Then, for all $i \geq 1$
\[  \mathrm{Ext}^i_{\widetilde{\mathbb H}_I}(X|_{\widetilde{\mathbb H}_I}, Y) =0 .
\]
\end{proposition}

\begin{proof}
Let $X$ be a regular anti-spherical discrete series of $\mathbb H$. By Lemma \ref{lem restrict calibrated}, $X|_{\widetilde{\mathbb H}_I}$ is a direct sum of irreducible $\widetilde{\mathbb H}_I$-modules of regular central characters. Note that any direct summand in $X|_{\widetilde{\mathbb H}_I}$ is anti-spherical. Otherwise, by Frobenius reciprocity, $X$ is a simple quotient of a module parabolically induced from a non-anti-spherical module, and so $X$ cannot be anti-spherical by Lemma \ref{lem induce spherical}.

Suppose, for some $i \geq 1$,
\[   \mathrm{Ext}_{\widetilde{\mathbb H}_I}^i(X|_{\widetilde{\mathbb H}_I}, Y)\neq 0 .
\]
Then, there is a direct summand $\widetilde{Y}$ in $X|_{\widetilde{\mathbb H}_I}$ such that 
\[  \mathrm{Ext}_{\widetilde{\mathbb H}_I}^i(\widetilde{Y}, Y)\neq 0 .
\]
However, there is only one anti-spherical module at each central character. Then $\widetilde{Y}\cong Y$. This contradicts the higher vanishing Ext-groups for discrete series \cite{Ch16}.
\end{proof}

We also remark that the branching problem for graded Hecke algebra does not correspond precisely with the branching problem for $p$-adic groups. For example, in type $A$ case, the restriction of the Steinberg module for graded Hecke algebra of type $A_n$ is still finite-dimensional and so cannot be projective, in contrast with \cite{Ch21, CS21}.

\subsection{Dropping the anti-spherical condition}

It is interesting to note that the vanishing Ext-branching in Theorem \ref{thm vanishing ext branching} fails in general if one drops the anti-spherical condition. We now consider $\mathbb H$ to be the graded Hecke algebra of type $H_4$. We use the discrete series at $\chi_3$, and let $X=\mathbb H^{(\chi_3, J_3)}$. The weights of $X$ are explicitly given by:
\begin{itemize}
\item $-s_4s_3\chi_3=-\frac{\sqrt{2}}{4(1+2a)}(2a,8a+2,1,30a+9)\\
=-[(9a+3)\alpha_1+(18a+\frac{11}{2})\alpha_2+(27a+8)\alpha_3+(21a+\frac{13}{2})\alpha_4],$
\item $-s_3s_4s_3\chi_3=-\frac{\sqrt{2}}{4(1+2a)}(6a+2,4a+2,-1,30a+9)\\=-[(9a+3)\alpha_1+(18a+\frac{11}{2})\alpha_2+(25a+\frac{17}{2})\alpha_3+(21a+\frac{13}{2})\alpha_4],$
\item $-s_2s_3s_4s_3\chi_3=-\frac{\sqrt{2}}{4(1+2a)}(1,2a,8a+2,30a+9)\\=-[(9a+3)\alpha_1+(16a+6)\alpha_2+(25a+\frac{17}{2})\alpha_3+(21a+\frac{13}{2})\alpha_4],$
\item $-s_1s_2s_3s_4s_3\chi_3=-\frac{\sqrt{2}}{4(1+2a)}(8a+3,14a+4,8a+2,26a+7)\\=-[(7a+3)\alpha_1+(16a+6)\alpha_2+(25a+\frac{17}{2})\alpha_3+(21a+\frac{13}{2})\alpha_4].$\\
\end{itemize}

Note that 
\[   \langle \alpha_4, -s_4s_3\chi_3 \rangle =\frac{1}{2}
\]
and so $X|_{\widetilde{\mathbb H}_{\left\{ \alpha_4\right\}}}$ has an irreducible $1$-dimensional non-tempered composition factor, denoted by $U$. On the other hand, for the $1$-dimensional discrete series, denoted by $\mathrm{St}$, it follows from \cite[Theorem 7.1]{Ch17} that
\[  \mathrm{Ext}_{\widetilde{\mathbb H}_{\left\{\alpha_4\right\}}}^1(X|_{\widetilde{\mathbb H}_{\left\{ \alpha_4\right\}}}, \mathrm{St}) \neq 0 .
\]

\subsection{Conjecture}
One may ask how about dropping the regularity condition:

\begin{conjecture}
We use the notations in Theorem \ref{thm vanishing ext branching}. Let $X'$ be an anti-spherical discrete series of $\mathbb H$. Let $I \subset \Delta$. Let $Y'$ be an anti-spherical discrete series of $\widetilde{\mathbb H}_{I}$. Then, for $i \geq 1$,
\[  \mathrm{Ext}^i_{\widetilde{\mathbb H}_I}(X'|_{\widetilde{\mathbb H}_I}, Y') = 0 .
\]
\end{conjecture}

One complication is that the restriction is not a semisimple module in general.

\subsection{Proof of Theorem \ref{thm ext branching for discrete series}}

We end this section with explaining Theorem \ref{thm ext branching for discrete series} -- an Ext-vanishing extension of Theorem \ref{thm construct minimal ds}. Let $DS$ be a discrete series of $\mathbb H$. Using the notations in the proof of Theorem \ref{thm construct minimal ds}, one can pick such discrete series submodule $DS'$ in $DS|_{\widetilde{\mathbb H}_{I_0}}$. One now has
\[   DS|_{\mathbb H_{I_0}}= X_0\oplus Y_0 ,
\]
where $X_0$ (resp. $Y_0$) contains all irreducible composition factors with the central character equal (resp. not equal) to that of $DS'\otimes \mathbb C_{\omega}$.

Indeed, the minimality forces that any irreducible composition factor in $DS|_{\mathbb H_{I_0}}$ with the same central character as $DS'\otimes \mathbb C_{\omega}$ is also an essentially discrete series. Thus, from the vanishing higher Ext-groups for discrete series \cite{Ch16}, we have: for $i \geq 1$,
\[  \mathrm{Ext}^i_{\widetilde{\mathbb H}_{I_0}}(X_0, DS')=0.
\]

It remains to show that
\[ \mathrm{Ext}^i_{\widetilde{\mathbb H}_{I_0}}(Y_0, DS')=0 .
\]
Suppose not to obtain contradiction. Then we find an irreducible $\mathbb H_{I_0}$-composition factor $Y$ in $Y_0$ such that 
$\mathrm{Ext}^i_{\widetilde{\mathbb H}_{I_0}}(Y, DS')\neq 0$. We can write $Y$ as $Y'\otimes \mathbb C_{k \varpi_{\alpha_0}}$ for some irreducible $\widetilde{\mathbb H}_{I_0}$-module $Y'$, $k \in \mathbb R$ and $\alpha_0\in \Delta-I_0$. Since $Y'$ has the same central character as $DS'$ and $Y$ has the same central character as $DS$, we can only have $k\varpi_{\alpha_0}=\pm \omega$. However, if $k\varpi_{\alpha_0}=\omega$, it contradicts the choice of $Y_0$, and if $k\varpi_{\alpha_0}=-\omega$, it contradicts that $DS$ is a discrete series.


\section{Springer Correspondence} \label{s springer correspond}

The classical Springer representation is defined as the top degree of the $\ell$-adic cohomology on some Springer fibers. For the non-crystallgraphic types, it seems to be no known geometry for such theory. On the other hand, Lusztig establishes how a total Springer representation can be lifted to a tempered module of a graded Hecke algebra. Now, the Springer theory can be modeled in the framework of graded Hecke algebras. 

For a nilpotent element in a complex semisimple Lie algebra, one can attach to a so-called $\mathfrak{sl}_2$-triple and in particular attach to a semisimple element in the $\mathfrak{sl}_2$-triple. From this viewpoint, one can use the central characters of $\mathbb H$ which support a tempered module with a real central character to replace nilpotent orbits. The distinguished nilpotent orbits should correspond to the Heckman-Opdam central characters. The cohomology of Springer fibers will be replaced by the $W$-structure of those tempered modules. One expects there is a natural bijection between the set of irreducible representations of $W$ and the set of tempered $\mathbb H$-modules of real central characters.

We shall use the symbols $\sigma_j$ $(j=1, \ldots, 34)$ to denote the irreducible representations of $W(H_4)$ corresponding to the characters in \cite[Table 4]{Gr74}\footnote{The character $\chi_8$ in \cite[Table 4]{Gr74} does not satisfy the orthogonality of characters and the values of $\chi_8$ on the conjugacy classes $K_8$ and $K_9$ should be $4\alpha$ rather than $-4\alpha$.}. In particular, $\sigma_1$ is the trivial representation, $\sigma_3$ is the reflection representation and $\sigma_{34}$ is the $48$-dimensional representation. We write down an expected/conjectural Springer representation attached to each Heckman-Opdam central character in Table \ref{table sp H4}. In the table, the second column represents the dimension of the corresponding $\sigma_i$, and for the columns under $\chi_i$, the entries represent the expected value $\mathrm{dim}~\mathrm{Hom}_W(DS_i|_W, \sigma_j)$ (all unfilled entries are zero) for the calibrated discrete series $DS_i$ with the central character $W\chi_i$. The entries under the columns $\chi_{16}', \chi_{17}', \chi_{17}''$ have similar meaning, but the values are for the non-calibrated modules at the central characters $W\chi_{16}$ and $W\chi_{17}$. From here we also expect that there exist antispherical discrete series at the central characters $W\chi_{16}$ and $W\chi_{17}$. Moreover, the expected dimension of the non-calibrated discrete series at $W\chi_{16}$ is 155, and the expected dimensions of the two non-calibrated discrete series at $W\chi_{17}$ are 101 and 331.

{\small
\begin{table} 
\begin{tabular}{|c|c|c|c|c|c|c|c|c|c|c|c|c|c|c|c|c|c|c|c|c|c|}
\hline
      & $\mathrm{dim}$  & $\chi_1$ & $\chi_3$ & $\chi_4$ & $\chi_{15}$ & $\chi_2$ & $\chi_6$ & $\chi_7$ & $\chi_5$ & $\chi_{8}$ & $\chi_{16}' $ & $\chi_{9}$ & $\chi_{10}$ & $\chi_{14}$ & $\chi_{16}$ & $\chi_{13}$ & $\chi_{17}'$ & $\chi_{17}''$ & $\chi_{17}$ & $\chi_{11}$ & $\chi_{12}$ \\
\hline
$\sigma_2$ & 1 & 1        &     0    &  1      & 0  & 1 &  0 & 1 & 1 & 1 & 1 & 1 & 0 & 0& 0 & 0 & 1 &  0 & 0 &  0 & 1 \\
\hline
$\sigma_6$ & 4&        &  1       &  0     &  0 &  0 &  1  & 0  & 0 & 0 & 0 & 0 & 1 & 0 & 0 & 0 & 0  & 0  & 0 & 1 & 0  \\
\hline
$\sigma_{4}$ & 4   &         &          &  1       &  0   & 1 &   0   & 1 & 1 & 1 & 1 & 1 & 0 & 1 & 1 & 1 & 1 & 0 & 0 & 0 & 1\\
\hline
$\sigma_{14}$  & 9 &  &  &   &   1     &  0 &  0 & 0  & 0 & 0 & 0 &  1 & 0 &0 & 0 & 0 & 1 & 1 & 0 & 1 &  1 \\
\hline
$\sigma_{12}$   & 9 &   &   &   &       & 1 & 0 & 1  &  1&  1 & 1 & 1 & 0 & 0 & 1 & 1 &  1  & 0 & 1 & 0 & 1\\
\hline
$\sigma_{21}$   &  16 &  &   &   &      & & 1 & 0  &  0  & 0 & 1 & 1 &  1  & 0 & 0 & 0 &   1  & 0  & 0 & 1 & 1\\
\hline
$\sigma_{19}$   &  16 &  &   &   &      & & & 1 & 1 & 1& 1 & 1 & 0 & 0 & 0 & 1 &  1  & 0    & 0  & 0 &  2   \\
\hline
$\sigma_{28}$   &  25&  &   &   &       & & &  & 1 & 1& 1 & 1 & 0 & 0 &  0 & 0 & 2 & 0 & 0  & 1   &  2\\
\hline
$\sigma_{32}$   & 36 &  &   &   &     & & &  & & 1 & 1 & 2 & 1 & 0 & 0 & 0 & 2  & 1  & 0   & 2 & 2  \\
\hline
$\sigma_{26}$   & 24&  &    &   &     & & &  & & 1 & 1& 1 & 0 & 0 & 0 & 0 &  1  & 0 & 0   & 0  & 1 \\
\hline
$\sigma_{24}$   & 24 &   &    &   &     & & &  & & & 1 & 1 & 0 & 0 & 0 &  0 &  1  & 0   & 0  & 0 & 1 \\
\hline
$\sigma_{33}$   & 40 &   &    &   &     & & &  & & & & 1 & 0  & 0 & 0 &  0 &  1  & 1  & 0  & 1  & 1\\
\hline
$\sigma_{30}$  & 30 &  &   &      &   & & &  & & & & &  1 & 0 & 0 & 0 & 0 &0  & 0 & 1  & 0    \\
\hline
$\sigma_7$      & 6 &  &   &    &    & & &  & & & & & & 1 &  1 & 1 &   0   & 0    & 0  & 0 & 0 \\
\hline
$\sigma_{16}$   & 16 &  &   &    &    & & &  & & & & & & & 1 & 1 &  0  & 0   & 1   & 0  & 0 \\
\hline
$\sigma_{29}$   & 30 &  &   &    &    & & &  & & & & & & & & 1 &   0   & 0  & 0  & 0  & 1 \\
\hline
$\sigma_{34}$  & 48 &  &    &   &    & & &  & & & & & & & &  &1   & 0   & 0   & 1 & 1 \\
\hline
$\sigma_{22}$  & 18 &  &    &   &    & & &  & & & & & & & & & 1    & 0   & 0  & 0 & 1\\
\hline
$\sigma_{17}$  &  16 &  &    &   &    & & &  & & & & & & & & &         & 1   & 0 & 1   & 0  \\
\hline
$\sigma_{15}$  &  10 &  &    &   &    & & &  & & & & & & & & &   & &1  & 0    & 0   \\
\hline
$\sigma_{25}$  & 24 &  &    &   &     & & &  & & & & & & & & &  &   & & 1  &  0  \\
\hline
$\sigma_{23}$  & 24 &  &   &   &      & & &  & & & & & & & & &  & & & & 1           \\
\hline
$\sigma_{9}$  & 8 &  &    &   &      & & &  & & & & & & & & &    & & &  & 1  \\  
\hline
\end{tabular}
 \caption{Springer correspondence for discrete series of type $H_4$ } \label{table sp H4}
\end{table}
}

The above correspondence agrees with the dimensions computed in previous sections, and also satisfies the following orthogonality relations: for two discrete series $DS, DS'$ of $\mathbb H$,
\[ \sum_{i=0}^4(-1)^i\mathrm{dim}~\mathrm{Ext}^i(DS, DS')=\sum_{i=0}^4(-1)^i\mathrm{dim}~\mathrm{Hom}_W(DS\otimes \wedge^iV, DS') =\left\{ \begin{array}{cc} 1 & \mbox{if $DS\cong DS'$} \\ 0 & \mbox{if $DS\not\cong DS'$}. \end{array} \right.
\]
This information is an important guidance on finding above correspondence while one needs a bit more to get the uniqueness. For example, there are two discrete series with $35$ dimensions, and one may use the minimally induced modules constructed in Sections \ref{ss minimal induce chi 17} and \ref{ss skew 16 ds} to decide the choice. We plan to give a detail account on this in a forthcoming article \cite{CH25+}.

According to \cite{LA82} (also see \cite[Page 412, Table C.2]{GP00}), the discrete series with the central characters $W\chi_1, W\chi_4$, $W\chi_2,  W\chi_6, W\chi_{5}, W\chi_{8}$ correspond to special representations under the Springer correspondence. $W\chi_8$ corresponds to the superspecial representation in the sense of \cite{Lu25}.


\part{Appendices}

\section{Appendix: Regular discrete series characters} \label{s rdsc}

\begin{itemize}
\item  $\chi_1=(17a+6)\alpha_1+(34a+\frac{23}{2})\alpha_2+(51a+\frac{33}{2})\alpha_3+(42a+13)\alpha_4$
\item  $\chi_2=(7a+3)\beta^1_2+(8a+2)\beta^2_2+(7a+\frac{9}{2})\beta^3_2+(2a+\frac{5}{2})\beta^4_2,$
\item   $\chi_3=(9a+3)\beta^1_3+(18a+\frac{11}{2})\beta^2_3+(15a+6)\beta^3_3-(6a+\frac{1}{2})\beta^4_3,$
 \item   $\chi_4=(10a+4)\beta^1_4+(20a+\frac{15}{2})\beta^2_4+(4a+\frac{7}{2})\beta^3_4+(14a+6)\beta^4_4$
\item $\chi_5=(5a+2)\beta^1_5+5\beta^2_5+(-5a+\frac{9}{2})\beta^3_5+\frac{7}{2}\beta^4_5$
\item $\chi_6=(7a+1)\beta^1_6+(-7a+\frac{9}{2})\beta^2_6+(-3a+7)\beta^3_6+(3a+\frac{11}{2})\beta^4_6$
 \item $\chi_7=(5a+3)\beta^1_7+(-5a+\frac{11}{2})\beta^2_7+\frac{15}{2}\beta^3_7+4\beta^4_7$
 \item   $\chi_8=(-14a+13)\beta^1_8+(-10a+9)\beta^2_8+(-20a+\frac{35}{2})\beta^3_8+(-4a+\frac{11}{2})\beta^4_8$
\item   $\chi_9=(a+1)\beta^1_9+\frac{1}{4}\beta^2_9+(-a+\frac{3}{2})\beta^3_9+\frac{3}{2}\beta^4_9$
 \item   $\chi_{10}=(-7a+8)\beta^1_{10}+(-2a+\frac{7}{2})\beta^2_{10}+(-7a+\frac{13}{2})\beta^3_{10}+(-8a+6)\beta^4_{10}$
\item    $\chi_{11}=(-18a+\frac{29}{2})\beta^1_{11}+(6a-\frac{7}{2})\beta^2_{11}+(-15a+\frac{27}{2})\beta^3_{11}+(-9a+\frac{15}{2})\beta^4_{11}$
 \item    $\chi_{12}=(-51a+42)\beta^1_{12}+(-17a+\frac{29}{2})\beta^2_{12}+(-42a+34)\beta^3_{12}+(-34a+\frac{57}{2})\beta^4_{12}$
\end{itemize}

\section{Appendix: Weight structure for regular discrete series}

The weight structure of regular discrete series is shown as follows. We only show those with dimensions less than or equal to 55. \\

\noindent
(1) The local region $\mathcal F^{(\chi_1, J_1)}=\left\{ w_o\right\}$. Then $\mathbb H^{(\chi_1, J_1)}$ is the so-called (1-dimensional) Steinberg module with the weight $-\chi_1$.\\
(2) Elements in $\mathcal F^{(\chi_2, J_2)}$:
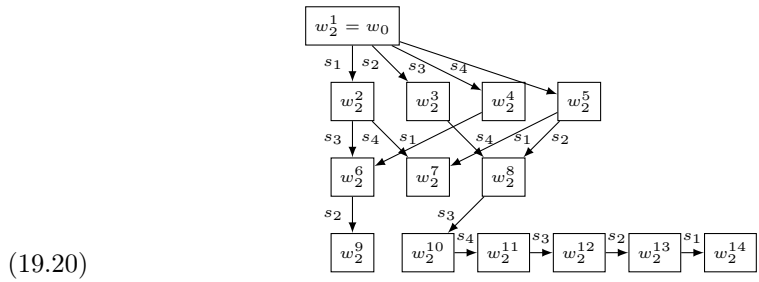
\begin{align}
\tikzset{int/.style={draw,minimum size=1em}}
\begin{tikzpicture}[auto,>=latex]
\node [int] (1node) [right of=1, distance=0cm]{\tiny{$w_2^1=w_0$}};
 \node (2) [below of=1, node distance=1cm, coordinate] {1};
 \node (3) [right of=2, node distance=1cm, coordinate] {2};
 \node (4) [right of=3, node distance=1cm, coordinate] {3};
 \node (5) [right of=4, node distance=1cm, coordinate] {4};
 \node (6) [below of=2, node distance=1cm, coordinate] {2};
 \node (7) [below of=3, node distance=1cm, coordinate] {3};
 \node (8) [below of=4, node distance=1cm, coordinate] {4};
 \node (9) [below of=6, node distance=1cm, coordinate] {6};
 \node (10) [right of=9, node distance=1cm, coordinate] {9};
 \node (11) [right of=10, node distance=1cm, coordinate] {10};
 \node (12) [right of=11, node distance=1cm, coordinate] {11};
 \node (13) [right of=12, node distance=1cm, coordinate] {12};
 \node (14) [right of=13, node distance=1cm, coordinate] {13};
 \node [int] (2node) [right of=2, distance=0cm] {\tiny{$w_2^2$}}; 
 \node [int] (3node) [right of=3, distance=0cm] {\tiny{$w_2^3$}};
 \node [int] (4node) [right of=4, distance=0cm] {\tiny{$w_2^4$}};
 \node [int] (5node) [right of=5, distance=0cm] {\tiny{$w_2^5$}};
 \node [int] (6node) [right of=6, distance=0cm] {\tiny{$w_2^6$}};
 \node [int] (7node) [right of=7, distance=0cm] {\tiny{$w_2^7$}}; 
 \node [int] (8node) [right of=8, distance=0cm] {\tiny{$w_2^8$}};
 \node [int] (9node) [right of=9, distance=0cm] {\tiny{$w_2^9$}};
 \node [int] (10node) [right of=10, distance=0cm] {\tiny{$w_2^{10}$}};
 \node [int] (11node) [right of=11, distance=0cm] {\tiny{$w_2^{11}$}};
 \node [int] (12node) [right of=12, distance=0cm] {\tiny{$w_2^{12}$}}; 
 \node [int] (13node) [right of=13, distance=0cm] {\tiny{$w_2^{13}$}};
 \node [int] (14node) [right of=14, distance=0cm] {\tiny{$w_2^{14}$}};
 \path[->] (1node) edge node[left] {\tiny{$s_1$}} (2node);
 \path[->] (1node) edge node[left] {\tiny{$s_2$}} (3node);
 \path[->] (1node) edge node[left] {\tiny{$s_3$}} (4node);
 \path[->] (1node) edge node[left] {\tiny{$s_4$}} (5node);
 \path[->] (2node) edge node[left] {\tiny{$s_3$}} (6node);
 \path[->] (2node) edge node[left] {\tiny{$s_4$}} (7node);
 \path[->] (3node) edge node[right] {\tiny{$s_4$}} (8node);
 \path[->] (4node) edge node[left] {\tiny{$s_1$}} (6node);
 \path[->] (5node) edge node[right] {\tiny{$s_1$}} (7node);
 \path[->] (5node) edge node[right] {\tiny{$s_2$}} (8node);
 \path[->] (6node) edge node[left] {\tiny{$s_2$}} (9node);
 \path[->] (8node) edge node[left] {\tiny{$s_3$}} (10node);
 \path[->] (10node) edge node[above] {\tiny{$s_4$}} (11node);
 \path[->] (11node) edge node[above] {\tiny{$s_3$}} (12node);
 \path[->] (12node) edge node[above] {\tiny{$s_2$}} (13node);
 \path[->] (13node) edge node[above] {\tiny{$s_1$}} (14node);
\end{tikzpicture}
\end{align}\\
(3) Elements in $\mathcal F^{(\chi_3, J_3)}$:
\tikzset{int/.style={draw,minimum size=1em}}
\begin{tikzpicture}
 \node [int] (1) {\tiny{$s_4s_3w_o$}};
 \node (2) [right of=1, node distance=0.8cm, coordinate] {1};
 \node [int] (2node) [right of=2, distance=0cm] {\tiny{$s_3s_4s_3w_o$}};
 \node (3) [right of=2node, node distance=1cm, coordinate] {2}; 
 \node [int] (3node) [right of=3, distance=0cm] {\tiny{$s_2s_3s_4s_3w_o$}};
 \path[->] (1) edge node[above] {\tiny{$s_{3}$}} (2node);
 \path[->] (2node) edge node[above] {\tiny{$s_{2}$}} (3node);
 \node (4) [right of=3node, node distance=1.2cm, coordinate] {3};
 \node [int] (4node) [right of=4, distance=0cm] {\tiny{$s_1s_2s_3s_4s_3w_o$}};
 \path[->] (3node) edge node[above] {\tiny{$s_{1}$}} (4node);
\end{tikzpicture}\\
(4) Elements in $\mathcal F^{(\chi_4, J_4)}$:
\tikzset{int/.style={draw,minimum size=1em}}
\begin{tikzpicture}
 \node [int] (1) {\tiny{$s_4w_o$}};
 \node (2) [right of=1, node distance=0.8cm, coordinate] {1};
 \node [int] (2node) [right of=2, distance=0cm] {\tiny{$w_o$}};
 \node (3) [right of=2node, node distance=1cm, coordinate] {2}; 
 \node [int] (3node) [right of=3, distance=0cm] {\tiny{$s_3w_o$}};
 \path[->] (1) edge node[above] {\tiny{$s_{4}$}} (2node);
 \path[->] (2node) edge node[above] {\tiny{$s_{3}$}} (3node);
 \node (4) [right of=3node, node distance=1.2cm, coordinate] {3};
 \node [int] (4node) [right of=4, distance=0cm] {\tiny{$s_2s_3w_o$}};
 \path[->] (3node) edge node[above] {\tiny{$s_{2}$}} (4node);
 \node (5) [right of=4node, node distance=1.2cm, coordinate] {4};
 \node [int] (5node) [right of=5, distance=0cm] {\tiny{$s_1s_2s_3w_o$}};
 \path[->] (4node) edge node[above] {\tiny{$s_{1}$}} (5node);
\end{tikzpicture}\\
(5) Elements in $\mathcal F^{(\chi_5, J_5)}$:
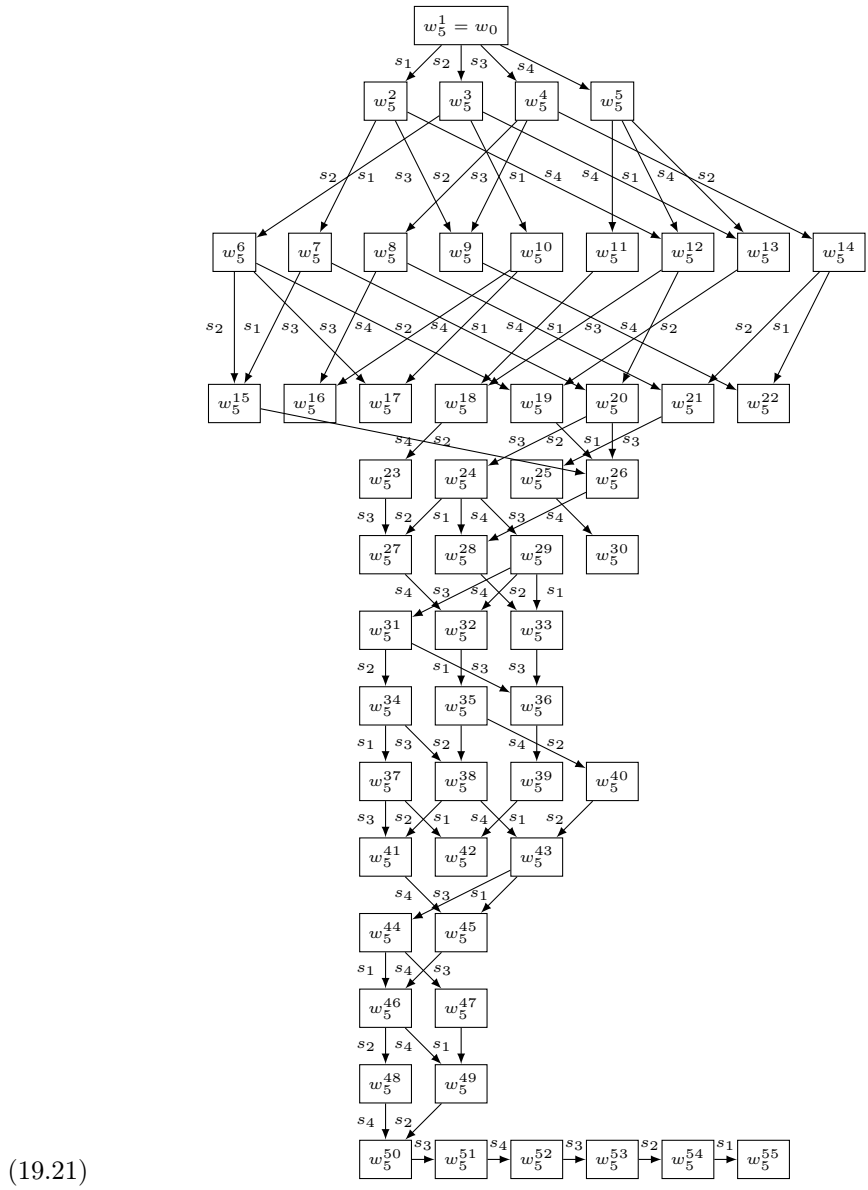
\begin{align}
\tikzset{int/.style={draw,minimum size=1em}}
\begin{tikzpicture}[auto,>=latex]
\node [int] (1node) [right of=1, distance=0cm]{\tiny{$w_5^1=w_0$}};
 \node (2) [below of=1, left of=1, node distance=1cm, coordinate] {1};
 \node (3) [right of=2, node distance=1cm, coordinate] {2};
 \node (4) [right of=3, node distance=1cm, coordinate] {3};
 \node (5) [right of=4, node distance=1cm, coordinate] {4};
 \node (6) [below of=2, left of=2, node distance=2cm, coordinate] {2};
 \node (7) [right of=6, node distance=1cm, coordinate] {6};
 \node (8) [right of=7, node distance=1cm, coordinate] {7};
 \node (9) [right of=8, node distance=1cm, coordinate] {8};
 \node (10) [right of=9, node distance=1cm, coordinate] {9};
 \node (11) [right of=10, node distance=1cm, coordinate] {10};
 \node (12) [right of=11, node distance=1cm, coordinate] {11};
 \node (13) [right of=12, node distance=1cm, coordinate] {12};
 \node (14) [right of=13, node distance=1cm, coordinate] {13};
 \node (15) [below of=6, node distance=2cm, coordinate] {6};
 \node (16) [right of=15, node distance=1cm, coordinate] {15};
 \node (17) [right of=16, node distance=1cm, coordinate] {16};
 \node (18) [right of=17, node distance=1cm, coordinate] {17};
 \node (19) [right of=18, node distance=1cm, coordinate] {18};
 \node (20) [right of=19, node distance=1cm, coordinate] {19};
 \node (21) [right of=20, node distance=1cm, coordinate] {20};
 \node (22) [right of=21, node distance=1cm, coordinate] {21};
 \node (23) [below of=17, node distance=1cm, coordinate] {17};
 \node (24) [right of=23, node distance=1cm, coordinate] {23};
 \node (25) [right of=24, node distance=1cm, coordinate] {24};
 \node (26) [right of=25, node distance=1cm, coordinate] {25};
 \node (27) [below of=23, node distance=1cm, coordinate] {23};
 \node (28) [right of=27, node distance=1cm, coordinate] {27};
 \node (29) [right of=28, node distance=1cm, coordinate] {28};
 \node (30) [right of=29, node distance=1cm, coordinate] {29};
 \node (31) [below of=27, node distance=1cm, coordinate] {27};
 \node (32) [right of=31, node distance=1cm, coordinate] {31};
 \node (33) [right of=32, node distance=1cm, coordinate] {32};
 \node (34) [below of=31, node distance=1cm, coordinate] {31};
 \node (35) [right of=34, node distance=1cm, coordinate] {34};
 \node (36) [right of=35, node distance=1cm, coordinate] {35};
 \node (37) [below of=34, node distance=1cm, coordinate] {34};
 \node (38) [right of=37, node distance=1cm, coordinate] {37};
 \node (39) [right of=38, node distance=1cm, coordinate] {38};
 \node (40) [right of=39, node distance=1cm, coordinate] {39};
 \node (41) [below of=37, node distance=1cm, coordinate] {37};
 \node (42) [right of=41, node distance=1cm, coordinate] {41};
 \node (43) [right of=42, node distance=1cm, coordinate] {42};
 \node (44) [below of=41, node distance=1cm, coordinate] {41};
 \node (45) [right of=44, node distance=1cm, coordinate] {44};
 \node (46) [below of=44, node distance=1cm, coordinate] {44};
 \node (47) [right of=46, node distance=1cm, coordinate] {46};
 \node (48) [below of=46, node distance=1cm, coordinate] {46};
 \node (49) [right of=48, node distance=1cm, coordinate] {48};
 \node (50) [below of=48, node distance=1cm, coordinate] {48};
 \node (51) [right of=50, node distance=1cm, coordinate] {50};
 \node (52) [right of=51, node distance=1cm, coordinate] {51};
 \node (53) [right of=52, node distance=1cm, coordinate] {52};
 \node (54) [right of=53, node distance=1cm, coordinate] {53};
 \node (55) [right of=54, node distance=1cm, coordinate] {54};
 \node [int] (2node) [right of=2, distance=0cm] {\tiny{$w_5^2$}}; 
 \node [int] (3node) [right of=3, distance=0cm] {\tiny{$w_5^3$}};
 \node [int] (4node) [right of=4, distance=0cm] {\tiny{$w_5^4$}};
 \node [int] (5node) [right of=5, distance=0cm] {\tiny{$w_5^5$}};
 \node [int] (6node) [right of=6, distance=0cm] {\tiny{$w_5^6$}};
 \node [int] (7node) [right of=7, distance=0cm] {\tiny{$w_5^7$}}; 
 \node [int] (8node) [right of=8, distance=0cm] {\tiny{$w_5^8$}};
 \node [int] (9node) [right of=9, distance=0cm] {\tiny{$w_5^9$}};
 \node [int] (10node) [right of=10, distance=0cm] {\tiny{$w_5^{10}$}};
 \node [int] (11node) [right of=11, distance=0cm] {\tiny{$w_5^{11}$}};
 \node [int] (12node) [right of=12, distance=0cm] {\tiny{$w_5^{12}$}}; 
 \node [int] (13node) [right of=13, distance=0cm] {\tiny{$w_5^{13}$}};
 \node [int] (14node) [right of=14, distance=0cm] {\tiny{$w_5^{14}$}};
 \node [int] (15node) [right of=15, distance=0cm] {\tiny{$w_5^{15}$}};
 \node [int] (16node) [right of=16, distance=0cm] {\tiny{$w_5^{16}$}};
 \node [int] (17node) [right of=17, distance=0cm] {\tiny{$w_5^{17}$}}; 
 \node [int] (18node) [right of=18, distance=0cm] {\tiny{$w_5^{18}$}};
 \node [int] (19node) [right of=19, distance=0cm] {\tiny{$w_5^{19}$}};
 \node [int] (20node) [right of=20, distance=0cm] {\tiny{$w_5^{20}$}};
 \node [int] (21node) [right of=21, distance=0cm] {\tiny{$w_5^{21}$}};
 \node [int] (22node) [right of=22, distance=0cm] {\tiny{$w_5^{22}$}}; 
 \node [int] (23node) [right of=23, distance=0cm] {\tiny{$w_5^{23}$}};
 \node [int] (24node) [right of=24, distance=0cm] {\tiny{$w_5^{24}$}};
 \node [int] (25node) [right of=25, distance=0cm] {\tiny{$w_5^{25}$}};
 \node [int] (26node) [right of=26, distance=0cm] {\tiny{$w_5^{26}$}};
 \node [int] (27node) [right of=27, distance=0cm] {\tiny{$w_5^{27}$}}; 
 \node [int] (28node) [right of=28, distance=0cm] {\tiny{$w_5^{28}$}};
 \node [int] (29node) [right of=29, distance=0cm] {\tiny{$w_5^{29}$}};
 \node [int] (30node) [right of=30, distance=0cm] {\tiny{$w_5^{30}$}};
 \node [int] (31node) [right of=31, distance=0cm] {\tiny{$w_5^{31}$}};
 \node [int] (32node) [right of=32, distance=0cm] {\tiny{$w_5^{32}$}}; 
 \node [int] (33node) [right of=33, distance=0cm] {\tiny{$w_5^{33}$}};
 \node [int] (34node) [right of=34, distance=0cm] {\tiny{$w_5^{34}$}};
 \node [int] (35node) [right of=35, distance=0cm] {\tiny{$w_5^{35}$}};
 \node [int] (36node) [right of=36, distance=0cm] {\tiny{$w_5^{36}$}};
 \node [int] (37node) [right of=37, distance=0cm] {\tiny{$w_5^{37}$}}; 
 \node [int] (38node) [right of=38, distance=0cm] {\tiny{$w_5^{38}$}};
 \node [int] (39node) [right of=39, distance=0cm] {\tiny{$w_5^{39}$}};
 \node [int] (40node) [right of=40, distance=0cm] {\tiny{$w_5^{40}$}};
 \node [int] (41node) [right of=41, distance=0cm] {\tiny{$w_5^{41}$}};
 \node [int] (42node) [right of=42, distance=0cm] {\tiny{$w_5^{42}$}}; 
 \node [int] (43node) [right of=43, distance=0cm] {\tiny{$w_5^{43}$}};
 \node [int] (44node) [right of=44, distance=0cm] {\tiny{$w_5^{44}$}};
 \node [int] (45node) [right of=45, distance=0cm] {\tiny{$w_5^{45}$}};
 \node [int] (46node) [right of=46, distance=0cm] {\tiny{$w_5^{46}$}};
 \node [int] (47node) [right of=47, distance=0cm] {\tiny{$w_5^{47}$}}; 
 \node [int] (48node) [right of=48, distance=0cm] {\tiny{$w_5^{48}$}};
 \node [int] (49node) [right of=49, distance=0cm] {\tiny{$w_5^{49}$}};
 \node [int] (50node) [right of=50, distance=0cm] {\tiny{$w_5^{50}$}};
 \node [int] (51node) [right of=51, distance=0cm] {\tiny{$w_5^{51}$}};
 \node [int] (52node) [right of=52, distance=0cm] {\tiny{$w_5^{52}$}}; 
 \node [int] (53node) [right of=53, distance=0cm] {\tiny{$w_5^{53}$}};
 \node [int] (54node) [right of=54, distance=0cm] {\tiny{$w_5^{54}$}};
 \node [int] (55node) [right of=55, distance=0cm] {\tiny{$w_5^{55}$}};
 \path[->] (1node) edge node[left] {\tiny{$s_1$}} (2node);
 \path[->] (1node) edge node[left] {\tiny{$s_2$}} (3node);
 \path[->] (1node) edge node[left] {\tiny{$s_3$}} (4node);
 \path[->] (1node) edge node[left] {\tiny{$s_4$}} (5node);
 \path[->] (2node) edge node[left] {\tiny{$s_2$}} (7node);
 \path[->] (2node) edge node[left] {\tiny{$s_3$}} (9node);
 \path[->] (2node) edge node[right] {\tiny{$s_4$}} (12node);
 \path[->] (3node) edge node[right] {\tiny{$s_1$}} (6node);
 \path[->] (3node) edge node[left] {\tiny{$s_3$}} (10node);
 \path[->] (3node) edge node[left] {\tiny{$s_4$}} (13node);
 \path[->] (4node) edge node[left] {\tiny{$s_2$}} (8node);
 \path[->] (4node) edge node[right] {\tiny{$s_1$}} (9node);
 \path[->] (4node) edge node[left] {\tiny{$s_4$}} (14node);
 \path[->] (5node) edge node[right] {\tiny{}} (11node);
 \path[->] (5node) edge node[left] {\tiny{$s_1$}} (12node);
 \path[->] (5node) edge node[right] {\tiny{$s_2$}} (13node);
 \path[->] (6node) edge node[left] {\tiny{$s_2$}} (15node);
 \path[->] (6node) edge node[left] {\tiny{$s_3$}} (17node);
 \path[->] (6node) edge node[left] {\tiny{$s_4$}} (19node);
 \path[->] (7node) edge node[left] {\tiny{$s_1$}} (15node);
 \path[->] (7node) edge node[left] {\tiny{$s_4$}} (20node);
 \path[->] (8node) edge node[left] {\tiny{$s_3$}} (16node);
 \path[->] (8node) edge node[left] {\tiny{$s_4$}} (21node);
 \path[->] (9node) edge node[right] {\tiny{$s_4$}} (22node);
 \path[->] (10node) edge node[left] {\tiny{$s_2$}} (16node);
 \path[->] (10node) edge node[right] {\tiny{$s_1$}} (17node);
 \path[->] (11node) edge node[right] {\tiny{$s_1$}} (18node);
 \path[->] (12node) edge node[right] {\tiny{$s_3$}} (18node);
 \path[->] (12node) edge node[right] {\tiny{$s_2$}} (20node);
 \path[->] (13node) edge node[left] {\tiny{}} (19node);
 \path[->] (14node) edge node[left] {\tiny{$s_2$}} (21node);
 \path[->] (14node) edge node[left] {\tiny{$s_1$}} (22node);
 \path[->] (15node) edge node[left] {\tiny{$s_4$}} (26node);
 \path[->] (18node) edge node[right] {\tiny{$s_2$}} (23node);
 \path[->] (19node) edge node[left] {\tiny{$s_2$}} (26node);
 \path[->] (20node) edge node[left] {\tiny{$s_3$}} (24node);
 \path[->] (20node) edge node[left] {\tiny{$s_1$}} (26node);
 \path[->] (21node) edge node[right] {\tiny{$s_3$}} (25node);
 \path[->] (23node) edge node[left] {\tiny{$s_3$}} (27node);
 \path[->] (24node) edge node[left] {\tiny{$s_2$}} (27node);
 \path[->] (24node) edge node[left] {\tiny{$s_1$}} (28node);
 \path[->] (24node) edge node[left] {\tiny{$s_4$}} (29node);
 \path[->] (25node) edge node[left] {\tiny{$s_4$}} (30node);
 \path[->] (26node) edge node[left] {\tiny{$s_3$}} (28node);
 \path[->] (27node) edge node[left] {\tiny{$s_4$}} (32node);
 \path[->] (28node) edge node[left] {\tiny{$s_4$}} (33node);
 \path[->] (29node) edge node[left] {\tiny{$s_3$}} (31node);
 \path[->] (29node) edge node[right] {\tiny{$s_2$}} (32node);
 \path[->] (29node) edge node[right] {\tiny{$s_1$}} (33node);
 \path[->] (31node) edge node[left] {\tiny{$s_2$}} (34node);
 \path[->] (31node) edge node[left] {\tiny{$s_1$}} (36node);
 \path[->] (32node) edge node[right] {\tiny{$s_3$}} (35node);
 \path[->] (33node) edge node[left] {\tiny{$s_3$}} (36node);
 \path[->] (34node) edge node[left] {\tiny{$s_1$}} (37node);
 \path[->] (34node) edge node[left] {\tiny{$s_3$}} (38node);
 \path[->] (35node) edge node[left] {\tiny{$s_2$}} (38node);
 \path[->] (35node) edge node[left] {\tiny{$s_4$}} (40node);
 \path[->] (36node) edge node[right] {\tiny{$s_2$}} (39node);
 \path[->] (37node) edge node[left] {\tiny{$s_3$}} (41node);
 \path[->] (37node) edge node[left] {\tiny{$s_2$}} (42node);
 \path[->] (38node) edge node[right] {\tiny{$s_1$}} (41node);
 \path[->] (38node) edge node[left] {\tiny{$s_4$}} (43node);
 \path[->] (39node) edge node[right] {\tiny{$s_1$}} (42node);
 \path[->] (40node) edge node[left] {\tiny{$s_2$}} (43node);
 \path[->] (41node) edge node[left] {\tiny{$s_4$}} (45node);
 \path[->] (43node) edge node[left] {\tiny{$s_3$}} (44node);
 \path[->] (43node) edge node[left] {\tiny{$s_1$}} (45node);
 \path[->] (44node) edge node[left] {\tiny{$s_1$}} (46node);
 \path[->] (44node) edge node[left] {\tiny{$s_4$}} (47node);
 \path[->] (45node) edge node[right] {\tiny{$s_3$}} (46node);
 \path[->] (46node) edge node[left] {\tiny{$s_2$}} (48node);
 \path[->] (46node) edge node[left] {\tiny{$s_4$}} (49node);
 \path[->] (47node) edge node[left] {\tiny{$s_1$}} (49node);
 \path[->] (48node) edge node[left] {\tiny{$s_4$}} (50node);
 \path[->] (49node) edge node[left] {\tiny{$s_2$}} (50node);
 \path[->] (50node) edge node[above] {\tiny{$s_3$}} (51node);
 \path[->] (51node) edge node[above] {\tiny{$s_4$}} (52node);
 \path[->] (52node) edge node[above] {\tiny{$s_3$}} (53node);
 \path[->] (53node) edge node[above] {\tiny{$s_2$}} (54node);
 \path[->] (54node) edge node[above] {\tiny{$s_1$}} (55node);
\end{tikzpicture}
\end{align}\\
(6) Elements in $\mathcal F^{(\chi_6, J_6)}$:
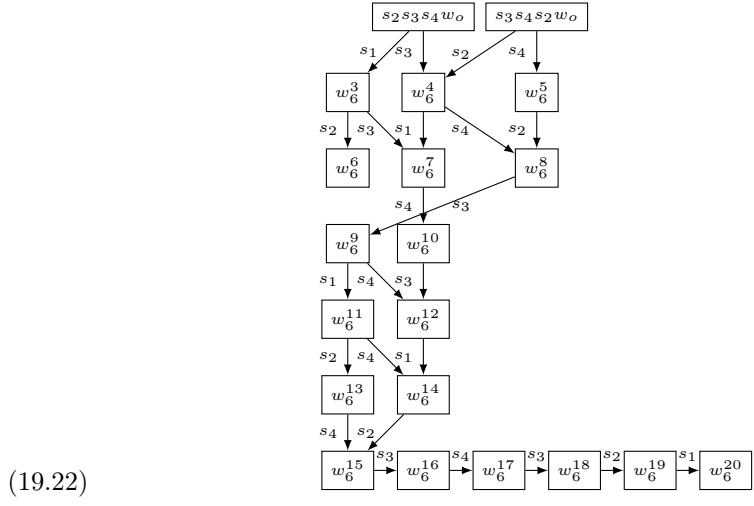
\begin{align}
\tikzset{int/.style={draw,minimum size=1em}}
\begin{tikzpicture}[auto,>=latex]
\node [int] (1node) [right of=1, distance=0cm]{\tiny{$s_2s_3s_4w_o$}};
 \node (2) [right of=1, node distance=1.5cm, coordinate] {1};
 \node (3) [below of=1, left of=1, node distance=1cm, coordinate] {1};
 \node (4) [right of=3, node distance=1cm, coordinate] {3};
 \node (5) [right of=4, node distance=1.5cm, coordinate] {4};
 \node (6) [below of=3, node distance=1cm, coordinate] {3};
 \node (7) [below of=4, node distance=1cm, coordinate] {4};
 \node (8) [below of=5, node distance=1cm, coordinate] {5};
 \node (9) [below of=6, node distance=1cm, coordinate] {6};
 \node (10) [right of=9, node distance=1cm, coordinate] {9};
 \node (11) [below of=9, node distance=1cm, coordinate] {9};
 \node (12) [right of=11, node distance=1cm, coordinate] {11};
 \node (13) [below of=11, node distance=1cm, coordinate] {11};
 \node (14) [right of=13, node distance=1cm, coordinate] {13};
 \node (15) [below of=13, node distance=1cm, coordinate] {13};
 \node (16) [right of=15, node distance=1cm, coordinate] {15};
 \node (17) [right of=16, node distance=1cm, coordinate] {16};
 \node (18) [right of=17, node distance=1cm, coordinate] {17};
 \node (19) [right of=18, node distance=1cm, coordinate] {18};
 \node (20) [right of=19, node distance=1cm, coordinate] {19};
 \node [int] (2node) [right of=2, distance=0cm] {\tiny{$s_3s_4s_2w_o$}}; 
 \node [int] (3node) [right of=3, distance=0cm] {\tiny{$w_6^3$}};
 \node [int] (4node) [right of=4, distance=0cm] {\tiny{$w_6^4$}};
 \node [int] (5node) [right of=5, distance=0cm] {\tiny{$w_6^5$}};
 \node [int] (6node) [right of=6, distance=0cm] {\tiny{$w_6^6$}};
 \node [int] (7node) [right of=7, distance=0cm] {\tiny{$w_6^7$}}; 
 \node [int] (8node) [right of=8, distance=0cm] {\tiny{$w_6^8$}};
 \node [int] (9node) [right of=9, distance=0cm] {\tiny{$w_6^9$}};
 \node [int] (10node) [right of=10, distance=0cm] {\tiny{$w_6^{10}$}};
 \node [int] (11node) [right of=11, distance=0cm] {\tiny{$w_6^{11}$}};
 \node [int] (12node) [right of=12, distance=0cm] {\tiny{$w_6^{12}$}}; 
 \node [int] (13node) [right of=13, distance=0cm] {\tiny{$w_6^{13}$}};
 \node [int] (14node) [right of=14, distance=0cm] {\tiny{$w_6^{14}$}};
 \node [int] (15node) [right of=15, distance=0cm] {\tiny{$w_6^{15}$}};
 \node [int] (16node) [right of=16, distance=0cm] {\tiny{$w_6^{16}$}};
 \node [int] (17node) [right of=17, distance=0cm] {\tiny{$w_6^{17}$}};
 \node [int] (18node) [right of=18, distance=0cm] {\tiny{$w_6^{18}$}}; 
 \node [int] (19node) [right of=19, distance=0cm] {\tiny{$w_6^{19}$}};
 \node [int] (20node) [right of=20, distance=0cm] {\tiny{$w_6^{20}$}};
 \path[->] (1node) edge node[left] {\tiny{$s_1$}} (3node);
 \path[->] (1node) edge node[left] {\tiny{$s_3$}} (4node);
 \path[->] (2node) edge node[left] {\tiny{$s_2$}} (4node);
 \path[->] (2node) edge node[left] {\tiny{$s_4$}} (5node);
 \path[->] (3node) edge node[left] {\tiny{$s_2$}} (6node);
 \path[->] (3node) edge node[left] {\tiny{$s_3$}} (7node);
 \path[->] (4node) edge node[left] {\tiny{$s_1$}} (7node);
 \path[->] (4node) edge node[left] {\tiny{$s_4$}} (8node);
 \path[->] (5node) edge node[left] {\tiny{$s_2$}} (8node);
 \path[->] (7node) edge node[left] {\tiny{$s_4$}} (10node);
 \path[->] (8node) edge node[right] {\tiny{$s_3$}} (9node);
 \path[->] (9node) edge node[left] {\tiny{$s_1$}} (11node);
 \path[->] (9node) edge node[left] {\tiny{$s_4$}} (12node);
 \path[->] (10node) edge node[left] {\tiny{$s_3$}} (12node);
 \path[->] (11node) edge node[left] {\tiny{$s_2$}} (13node);
 \path[->] (11node) edge node[left] {\tiny{$s_4$}} (14node);
 \path[->] (12node) edge node[left] {\tiny{$s_1$}} (14node);
 \path[->] (13node) edge node[left] {\tiny{$s_4$}} (15node);
 \path[->] (14node) edge node[left] {\tiny{$s_2$}} (15node);
 \path[->] (15node) edge node[above] {\tiny{$s_3$}} (16node);
 \path[->] (16node) edge node[above] {\tiny{$s_4$}} (17node);
 \path[->] (17node) edge node[above] {\tiny{$s_3$}} (18node);
 \path[->] (18node) edge node[above] {\tiny{$s_2$}} (19node);
 \path[->] (19node) edge node[above] {\tiny{$s_1$}} (20node);
\end{tikzpicture}
\end{align}\\
(7) Elements in $\mathcal F^{(\chi_7, J_7)}$:
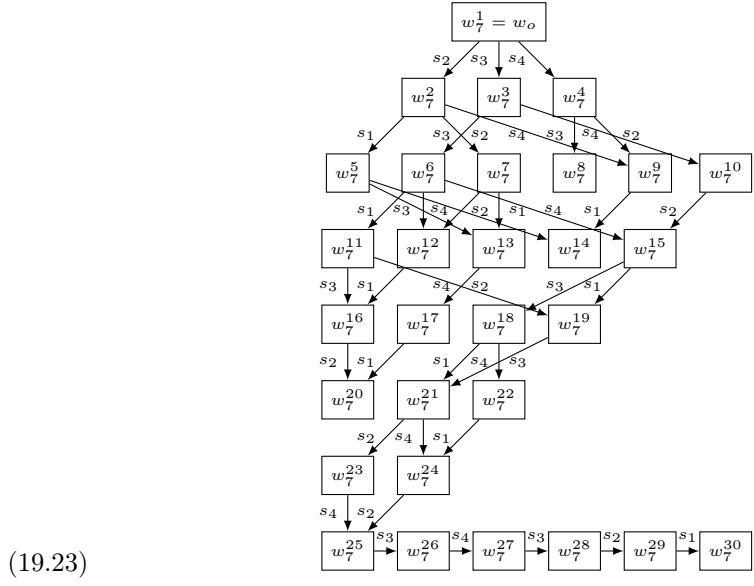
\begin{align}
\tikzset{int/.style={draw,minimum size=1em}}
\begin{tikzpicture}[auto,>=latex]
\node [int] (1node) [right of=1, distance=0cm]{\tiny{$w_7^1=w_o$}};
 \node (2) [below of=1, left of=1, node distance=1cm, coordinate] {1};
 \node (3) [right of=2, node distance=1cm, coordinate] {2};
 \node (4) [right of=3, node distance=1cm, coordinate] {3};
 \node (5) [below of=2, left of=2, node distance=1cm, coordinate] {2};
 \node (6) [right of=5, node distance=1cm, coordinate] {5};
 \node (7) [below of=3, node distance=1cm, coordinate] {3};
 \node (8) [below of=4, node distance=1cm, coordinate] {4};
 \node (9) [right of=8, node distance=1cm, coordinate] {8};
 \node (10) [right of=9, node distance=1cm, coordinate] {9};
 \node (11) [below of=5, node distance=1cm, coordinate] {5};
 \node (12) [right of=11, node distance=1cm, coordinate] {11};
 \node (13) [right of=12, node distance=1cm, coordinate] {12};
 \node (14) [right of=13, node distance=1cm, coordinate] {13};
 \node (15) [right of=14, node distance=1cm, coordinate] {14};
 \node (16) [below of=11, node distance=1cm, coordinate] {11};
 \node (17) [right of=16, node distance=1cm, coordinate] {16};
 \node (18) [right of=17, node distance=1cm, coordinate] {17};
 \node (19) [right of=18, node distance=1cm, coordinate] {18};
 \node (20) [below of=16, node distance=1cm, coordinate] {16};
 \node (21) [right of=20, node distance=1cm, coordinate] {20};
 \node (22) [right of=21, node distance=1cm, coordinate] {21};
 \node (23) [below of=20, node distance=1cm, coordinate] {20};
 \node (24) [right of=23, node distance=1cm, coordinate] {23};
 \node (25) [below of=23, node distance=1cm, coordinate] {23};
 \node (26) [right of=25, node distance=1cm, coordinate] {25};
 \node (27) [right of=26, node distance=1cm, coordinate] {26};
 \node (28) [right of=27, node distance=1cm, coordinate] {27};
 \node (29) [right of=28, node distance=1cm, coordinate] {28};
 \node (30) [right of=29, node distance=1cm, coordinate] {29};
 \node [int] (2node) [right of=2, distance=0cm] {\tiny{$w_7^2$}}; 
 \node [int] (3node) [right of=3, distance=0cm] {\tiny{$w_7^3$}};
 \node [int] (4node) [right of=4, distance=0cm] {\tiny{$w_7^4$}};
 \node [int] (5node) [right of=5, distance=0cm] {\tiny{$w_7^5$}};
 \node [int] (6node) [right of=6, distance=0cm] {\tiny{$w_7^6$}};
 \node [int] (7node) [right of=7, distance=0cm] {\tiny{$w_7^7$}}; 
 \node [int] (8node) [right of=8, distance=0cm] {\tiny{$w_7^8$}};
 \node [int] (9node) [right of=9, distance=0cm] {\tiny{$w_7^9$}};
 \node [int] (10node) [right of=10, distance=0cm] {\tiny{$w_7^{10}$}};
 \node [int] (11node) [right of=11, distance=0cm] {\tiny{$w_7^{11}$}};
 \node [int] (12node) [right of=12, distance=0cm] {\tiny{$w_7^{12}$}}; 
 \node [int] (13node) [right of=13, distance=0cm] {\tiny{$w_7^{13}$}};
 \node [int] (14node) [right of=14, distance=0cm] {\tiny{$w_7^{14}$}};
 \node [int] (15node) [right of=15, distance=0cm] {\tiny{$w_7^{15}$}};
 \node [int] (16node) [right of=16, distance=0cm] {\tiny{$w_7^{16}$}};
 \node [int] (17node) [right of=17, distance=0cm] {\tiny{$w_7^{17}$}};
 \node [int] (18node) [right of=18, distance=0cm] {\tiny{$w_7^{18}$}}; 
 \node [int] (19node) [right of=19, distance=0cm] {\tiny{$w_7^{19}$}};
 \node [int] (20node) [right of=20, distance=0cm] {\tiny{$w_7^{20}$}};
 \node [int] (21node) [right of=21, distance=0cm] {\tiny{$w_7^{21}$}};
 \node [int] (22node) [right of=22, distance=0cm] {\tiny{$w_7^{22}$}}; 
 \node [int] (23node) [right of=23, distance=0cm] {\tiny{$w_7^{23}$}};
 \node [int] (24node) [right of=24, distance=0cm] {\tiny{$w_7^{24}$}};
 \node [int] (25node) [right of=25, distance=0cm] {\tiny{$w_7^{25}$}};
 \node [int] (26node) [right of=26, distance=0cm] {\tiny{$w_7^{26}$}};
 \node [int] (27node) [right of=27, distance=0cm] {\tiny{$w_7^{27}$}};
 \node [int] (28node) [right of=28, distance=0cm] {\tiny{$w_7^{28}$}}; 
 \node [int] (29node) [right of=29, distance=0cm] {\tiny{$w_7^{29}$}};
 \node [int] (30node) [right of=30, distance=0cm] {\tiny{$w_7^{30}$}};
 \path[->] (1node) edge node[left] {\tiny{$s_2$}} (2node);
 \path[->] (1node) edge node[left] {\tiny{$s_3$}} (3node);
 \path[->] (1node) edge node[left] {\tiny{$s_4$}} (4node);
 \path[->] (2node) edge node[left] {\tiny{$s_1$}} (5node);
 \path[->] (2node) edge node[left] {\tiny{$s_3$}} (7node);
 \path[->] (2node) edge node[left] {\tiny{$s_4$}} (9node);
 \path[->] (3node) edge node[right] {\tiny{$s_2$}} (6node);
 \path[->] (3node) edge node[left] {\tiny{$s_4$}} (10node);
 \path[->] (4node) edge node[left] {\tiny{$s_3$}} (8node);
 \path[->] (4node) edge node[right] {\tiny{$s_2$}} (9node);
 \path[->] (5node) edge node[left] {\tiny{$s_3$}} (13node);
 \path[->] (5node) edge node[left] {\tiny{$s_4$}} (14node);
 \path[->] (6node) edge node[left] {\tiny{$s_1$}} (11node);
 \path[->] (6node) edge node[left] {\tiny{}} (12node);
 \path[->] (6node) edge node[right] {\tiny{$s_4$}} (15node);
 \path[->] (7node) edge node[right] {\tiny{$s_2$}} (12node);
 \path[->] (7node) edge node[right] {\tiny{$s_1$}} (13node);
 \path[->] (9node) edge node[left] {\tiny{$s_1$}} (14node);
 \path[->] (10node) edge node[left] {\tiny{$s_2$}} (15node);
 \path[->] (11node) edge node[left] {\tiny{$s_3$}} (16node);
 \path[->] (11node) edge node[left] {\tiny{$s_4$}} (19node);
 \path[->] (12node) edge node[left] {\tiny{$s_1$}} (16node);
 \path[->] (13node) edge node[right] {\tiny{$s_2$}} (17node);
 \path[->] (15node) edge node[left] {\tiny{$s_3$}} (18node);
 \path[->] (15node) edge node[left] {\tiny{$s_1$}} (19node);
 \path[->] (16node) edge node[left] {\tiny{$s_2$}} (20node);
 \path[->] (17node) edge node[left] {\tiny{$s_1$}} (20node);
 \path[->] (18node) edge node[left] {\tiny{$s_1$}} (21node);
 \path[->] (18node) edge node[left] {\tiny{$s_4$}} (22node);
 \path[->] (19node) edge node[right] {\tiny{$s_3$}} (21node);
 \path[->] (21node) edge node[left] {\tiny{$s_2$}} (23node);
 \path[->] (21node) edge node[left] {\tiny{$s_4$}} (24node);
 \path[->] (22node) edge node[left] {\tiny{$s_1$}} (24node);
 \path[->] (23node) edge node[left] {\tiny{$s_4$}} (25node);
 \path[->] (24node) edge node[left] {\tiny{$s_2$}} (25node);
 \path[->] (25node) edge node[above] {\tiny{$s_3$}} (26node);
 \path[->] (26node) edge node[above] {\tiny{$s_4$}} (27node);
 \path[->] (27node) edge node[above] {\tiny{$s_3$}} (28node);
 \path[->] (28node) edge node[above] {\tiny{$s_2$}} (29node);
 \path[->] (29node) edge node[above] {\tiny{$s_1$}} (30node);
\end{tikzpicture}
\end{align}
\newpage





\section{Appendix: Weight structure for calibrated modules with central characters $\chi_{16}$ and $\chi_{17}$}

Let $w_{16}^1=s_3s_4s_3s_4s_2s_3$ and the corresponding element in the SageMath is $W[119]$. 

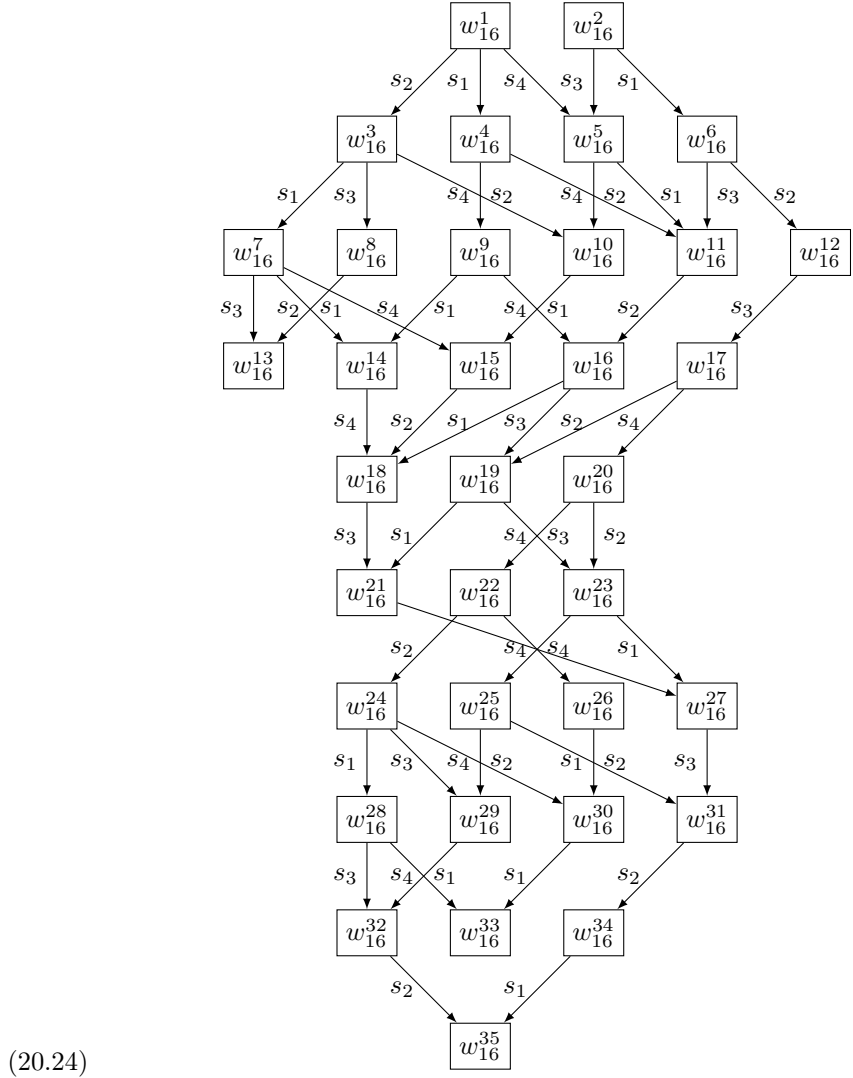
\begin{align} \label{diagram local region 16}
\tikzset{int/.style={draw,minimum size=1em}}
\begin{tikzpicture}[auto,>=latex]
\node [int] (1node) [right of=1, distance=0cm]{$w_{16}^1$};
\node (2) [right of=1, node distance=1.5cm, coordinate] {1};
\node (3) [left of=4, node distance=1.5cm, coordinate] {4};
\node (4) [below of=1, node distance=1.5cm, coordinate] {1};
\node (5) [right of=4, node distance=1.5cm, coordinate] {4};
\node (6) [right of=5, node
distance=1.5cm, coordinate] {5};
\node (7) [left of=3, below of=3, node distance=1.5cm, coordinate] {3};
\node (8) [right of=7, node distance=1.5cm, coordinate] {7};
\node (9) [right of=8, node distance=1.5cm, coordinate] {8};
\node (10) [right of=9, node distance=1.5cm, coordinate] {9};
\node (11) [right of=10, node distance=1.5cm, coordinate] {10};
\node (12) [right of=11, node distance=1.5cm, coordinate] {11};
\node (13) [below of=7, node distance=1.5cm, coordinate] {7};
\node (14) [right of=13, node distance=1.5cm, coordinate] {13};
\node (15) [right of=14, node distance=1.5cm, coordinate] {14};
\node (16) [right of=15, node distance=1.5cm, coordinate] {15};
\node (17) [right of=16, node distance=1.5cm, coordinate] {16};
\node (18) [below of=14, node distance=1.5cm, coordinate] {14};
\node (19) [right of=18, node distance=1.5cm, coordinate] {18};
\node (20) [right of=19, node distance=1.5cm, coordinate] {19};
\node (21) [below of=18, node distance=1.5cm, coordinate] {18};
\node (22) [right of=21, node distance=1.5cm, coordinate] {21};
\node (23) [right of=22, node distance=1.5cm, coordinate] {22};
\node (24) [below of=21, node distance=1.5cm, coordinate] {21};
\node (25) [right of=24, node distance=1.5cm, coordinate] {24};
\node (26) [right of=25, node distance=1.5cm, coordinate] {25};
\node (27) [right of=26, node distance=1.5cm, coordinate] {26};
\node (28) [below of=24, node distance=1.5cm, coordinate] {24};
\node (29) [right of=28, node distance=1.5cm, coordinate] {28};
\node (30) [right of=29, node distance=1.5cm, coordinate] {29};
\node (31) [right of=30, node distance=1.5cm, coordinate] {30};
\node (32) [below of=28, node distance=1.5cm, coordinate] {28};
\node (33) [right of=32, node distance=1.5cm, coordinate] {32};
\node (34) [right of=33, node distance=1.5cm, coordinate] {33};
\node (35) [below of=33, node distance=1.5cm, coordinate] {33};
\node [int] (2node) [right of=2, distance=0cm] {$w_{16}^2$}; 
\node [int] (3node) [right of=3, distance=0cm] {$w_{16}^3$};
\node [int] (4node) [right of=4, distance=0cm] {$w_{16}^4$};
\node [int] (5node) [right of=5, distance=0cm] {$w_{16}^5$};
\node [int] (6node) [right of=6, distance=0cm] {$w_{16}^6$};
\node [int] (7node) [right of=7, distance=0cm] {$w_{16}^7$}; 
\node [int] (8node) [right of=8, distance=0cm] {$w_{16}^8$};
\node [int] (9node) [right of=9, distance=0cm] {$w_{16}^9$};
\node [int] (10node) [right of=10, distance=0cm] {$w_{16}^{10}$};
\node [int] (11node) [right of=11, distance=0cm] {$w_{16}^{11}$};
\node [int] (12node) [right of=12, distance=0cm] {$w_{16}^{12}$}; 
\node [int] (13node) [right of=13, distance=0cm] {$w_{16}^{13}$};
\node [int] (14node) [right of=14, distance=0cm] {$w_{16}^{14}$};
\node [int] (15node) [right of=15, distance=0cm] {$w_{16}^{15}$};
\node [int] (16node) [right of=16, distance=0cm] {$w_{16}^{16}$};
\node [int] (17node) [right of=17, distance=0cm] {$w_{16}^{17}$};
\node [int] (18node) [right of=18, distance=0cm] {$w_{16}^{18}$};
\node [int] (19node) [right of=19, distance=0cm] {$w_{16}^{19}$}; 
\node [int] (20node) [right of=20, distance=0cm] {$w_{16}^{20}$};
\node [int] (21node) [right of=21, distance=0cm] {$w_{16}^{21}$};
\node [int] (22node) [right of=22, distance=0cm] {$w_{16}^{22}$};
\node [int] (23node) [right of=23, distance=0cm] {$w_{16}^{23}$};
\node [int] (24node) [right of=24, distance=0cm] {$w_{16}^{24}$}; 
\node [int] (25node) [right of=25, distance=0cm] {$w_{16}^{25}$};
\node [int] (26node) [right of=26, distance=0cm] {$w_{16}^{26}$};
\node [int] (27node) [right of=27, distance=0cm] {$w_{16}^{27}$};
\node [int] (28node) [right of=28, distance=0cm] {$w_{16}^{28}$};
\node [int] (29node) [right of=29, distance=0cm] {$w_{16}^{29}$}; 
\node [int] (30node) [right of=30, distance=0cm] {$w_{16}^{30}$};
\node [int] (31node) [right of=31, distance=0cm] {$w_{16}^{31}$};
\node [int] (32node) [right of=32, distance=0cm] {$w_{16}^{32}$};
\node [int] (33node) [right of=33, distance=0cm] {$w_{16}^{33}$};
\node [int] (34node) [right of=34, distance=0cm] {$w_{16}^{34}$};
\node [int] (35node) [right of=35, distance=0cm] {$w_{16}^{35}$};
\path[->] (1node) edge node[left] {$s_2$} (3node);
\path[->] (1node) edge node[left] {$s_1$} (4node);
\path[->] (1node) edge node[left] {$s_4$} (5node);
\path[->] (2node) edge node[left] {$s_3$} (5node);
\path[->] (2node) edge node[left] {$s_1$} (6node);
\path[->] (3node) edge node[left] {$s_1$} (7node);
\path[->] (3node) edge node[left] {$s_3$} (8node);
\path[->] (3node) edge node[left] {$s_4$} (10node);
\path[->] (4node) edge node[right] {$s_2$} (9node);
\path[->] (4node) edge node[left] {$s_4$} (11node);
\path[->] (5node) edge node[right] {$s_2$} (10node);
\path[->] (5node) edge node[right] {$s_1$} (11node);
\path[->] (6node) edge node[right] {$s_3$} (11node);
\path[->] (6node) edge node[right] {$s_2$} (12node);
\path[->] (7node) edge node[left] {$s_3$} (13node);
\path[->] (7node) edge node[left] {$s_2$} (14node);
\path[->] (7node) edge node[right] {$s_4$} (15node);
\path[->] (8node) edge node[right] {$s_1$} (13node);
\path[->] (9node) edge node[right] {$s_1$} (14node);
\path[->] (9node) edge node[left] {$s_4$} (16node);
\path[->] (11node) edge node[left] {$s_2$} (16node);
\path[->] (12node) edge node[left] {$s_3$} (17node);
\path[->] (10node) edge node[right] {$s_1$} (15node);
\path[->] (14node) edge node[left] {$s_4$} (18node);
\path[->] (15node) edge node[left] {$s_2$} (18node);
\path[->] (16node) edge node[left] {$s_1$} (18node);
\path[->] (16node) edge node[left] {$s_3$} (19node);
\path[->] (17node) edge node[left] {$s_2$} (19node);
\path[->] (17node) edge node[left] {$s_4$} (20node);
\path[->] (18node) edge node[left] {$s_3$} (21node);
\path[->] (19node) edge node[left] {$s_1$} (21node);
\path[->] (19node) edge node[left] {$s_4$} (23node);
\path[->] (20node) edge node[right] {$s_3$} (22node);
\path[->] (20node) edge node[right] {$s_2$} (23node);
\path[->] (21node) edge node[left] {$s_4$} (27node);
\path[->] (22node) edge node[left] {$s_2$} (24node);
\path[->] (22node) edge node[right] {$s_4$} (26node);
\path[->] (23node) edge node[left] {} (25node);
\path[->] (23node) edge node[left] {$s_1$} (27node);
\path[->] (24node) edge node[left] {$s_1$} (28node);
\path[->] (24node) edge node[left] {$s_3$} (29node);
\path[->] (24node) edge node[left] {$s_4$} (30node);
\path[->] (25node) edge node[right] {$s_2$} (29node);
\path[->] (25node) edge node[left] {$s_1$} (31node);
\path[->] (26node) edge node[right] {$s_2$} (30node);
\path[->] (27node) edge node[left] {$s_3$} (31node);
\path[->] (28node) edge node[left] {$s_3$} (32node);
\path[->] (28node) edge node[left] {$s_4$} (33node);
\path[->] (29node) edge node[right] {$s_1$} (32node);
\path[->] (30node) edge node[left] {$s_1$} (33node);
\path[->] (31node) edge node[left] {$s_2$} (34node);
\path[->] (32node) edge node[left] {$s_2$} (35node);
\path[->] (34node) edge node[left] {$s_1$} (35node);
\end{tikzpicture}
\end{align}

Let $w_{17}^1=s_3s_4s_3s_4s_1s_2s_3s_4s_3$. The corresponding element in SageMath is $W[345]$.

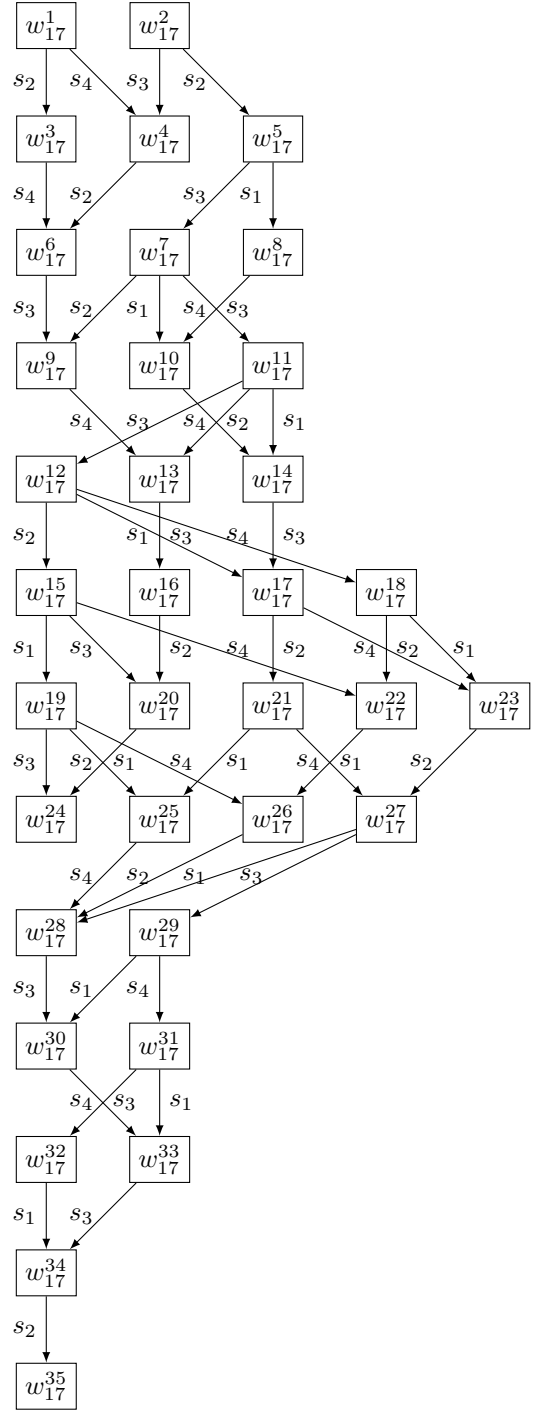
\begin{align} \label{diagram local region 17}
\tikzset{int/.style={draw,minimum size=1em}}
\begin{tikzpicture}[auto,>=latex]
 \node [int] (1node) [right of=1, distance=0cm]{$w_{17}^1$};
 \node (2) [right of=1, node distance=1.5cm, coordinate] {1};
 \node (3) [below of=1, node distance=1.5cm, coordinate] {1};
 \node (4) [right of=3, node distance=1.5cm, coordinate] {3};
 \node (5) [right of=4, node distance=1.5cm, coordinate] {4};
 \node (6) [below of=3, node distance=1.5cm, coordinate] {3};
 \node (7) [right of=6, node distance=1.5cm, coordinate] {6};
 \node (8) [right of=7, node distance=1.5cm, coordinate] {7};
 \node (9) [below of=6, node distance=1.5cm, coordinate] {6};
 \node (10) [right of=9, node distance=1.5cm, coordinate] {9};
 \node (11) [right of=10, node distance=1.5cm, coordinate] {10};
 \node (12) [below of=9, node distance=1.5cm, coordinate] {9};
 \node (13) [right of=12, node distance=1.5cm, coordinate] {12};
 \node (14) [right of=13, node distance=1.5cm, coordinate] {13};
 \node (15) [below of=12, node distance=1.5cm, coordinate] {12};
 \node (16) [right of=15, node distance=1.5cm, coordinate] {15};
 \node (17) [right of=16, node distance=1.5cm, coordinate] {16};
 \node (18) [right of=17, node distance=1.5cm, coordinate] {17};
 \node (19) [below of=15, node distance=1.5cm, coordinate] {15};
 \node (20) [right of=19, node distance=1.5cm, coordinate] {19};
 \node (21) [right of=20, node distance=1.5cm, coordinate] {20};
 \node (22) [right of=21, node distance=1.5cm, coordinate] {21};
 \node (23) [right of=22, node distance=1.5cm, coordinate] {22};
 \node (24) [below of=19, node distance=1.5cm, coordinate] {19};
 \node (25) [right of=24, node distance=1.5cm, coordinate] {24};
 \node (26) [right of=25, node distance=1.5cm, coordinate] {25};
 \node (27) [right of=26, node distance=1.5cm, coordinate] {26};
 \node (28) [below of=24, node distance=1.5cm, coordinate] {24};
 \node (29) [right of=28, node distance=1.5cm, coordinate] {28};
 \node (30) [below of=28, node distance=1.5cm, coordinate] {28};
 \node (31) [right of=30, node distance=1.5cm, coordinate] {30};
 \node (32) [below of=30, node distance=1.5cm, coordinate] {30};
 \node (33) [right of=32, node distance=1.5cm, coordinate] {32};
 \node (34) [below of=32, node distance=1.5cm, coordinate] {32};
 \node (35) [below of=34, node distance=1.5cm, coordinate] {34};
 \node [int] (2node) [right of=2, distance=0cm] {$w_{17}^2$}; 
 \node [int] (3node) [right of=3, distance=0cm] {$w_{17}^3$};
 \node [int] (4node) [right of=4, distance=0cm] {$w_{17}^4$};
 \node [int] (5node) [right of=5, distance=0cm] {$w_{17}^5$};
 \node [int] (6node) [right of=6, distance=0cm] {$w_{17}^6$};
 \node [int] (7node) [right of=7, distance=0cm] {$w_{17}^7$}; 
 \node [int] (8node) [right of=8, distance=0cm] {$w_{17}^8$};
 \node [int] (9node) [right of=9, distance=0cm] {$w_{17}^9$};
 \node [int] (10node) [right of=10, distance=0cm] {$w_{17}^{10}$};
 \node [int] (11node) [right of=11, distance=0cm] {$w_{17}^{11}$};
 \node [int] (12node) [right of=12, distance=0cm] {$w_{17}^{12}$}; 
 \node [int] (13node) [right of=13, distance=0cm] {$w_{17}^{13}$};
 \node [int] (14node) [right of=14, distance=0cm] {$w_{17}^{14}$};
 \node [int] (15node) [right of=15, distance=0cm] {$w_{17}^{15}$};
 \node [int] (16node) [right of=16, distance=0cm] {$w_{17}^{16}$};
 \node [int] (17node) [right of=17, distance=0cm] {$w_{17}^{17}$};
 \node [int] (18node) [right of=18, distance=0cm] {$w_{17}^{18}$};
 \node [int] (19node) [right of=19, distance=0cm] {$w_{17}^{19}$}; 
 \node [int] (20node) [right of=20, distance=0cm] {$w_{17}^{20}$};
 \node [int] (21node) [right of=21, distance=0cm] {$w_{17}^{21}$};
 \node [int] (22node) [right of=22, distance=0cm] {$w_{17}^{22}$};
 \node [int] (23node) [right of=23, distance=0cm] {$w_{17}^{23}$};
 \node [int] (24node) [right of=24, distance=0cm] {$w_{17}^{24}$}; 
 \node [int] (25node) [right of=25, distance=0cm] {$w_{17}^{25}$};
 \node [int] (26node) [right of=26, distance=0cm] {$w_{17}^{26}$};
 \node [int] (27node) [right of=27, distance=0cm] {$w_{17}^{27}$};
 \node [int] (28node) [right of=28, distance=0cm] {$w_{17}^{28}$};
 \node [int] (29node) [right of=29, distance=0cm] {$w_{17}^{29}$}; 
 \node [int] (30node) [right of=30, distance=0cm] {$w_{17}^{30}$};
 \node [int] (31node) [right of=31, distance=0cm] {$w_{17}^{31}$};
 \node [int] (32node) [right of=32, distance=0cm] {$w_{17}^{32}$};
 \node [int] (33node) [right of=33, distance=0cm] {$w_{17}^{33}$};
 \node [int] (34node) [right of=34, distance=0cm] {$w_{17}^{34}$};
 \node [int] (35node) [right of=35, distance=0cm] {$w_{17}^{35}$};
 \path[->] (1node) edge node[left] {$s_2$} (3node);
 \path[->] (1node) edge node[left] {$s_4$} (4node);
 \path[->] (2node) edge node[left] {$s_3$} (4node);
 \path[->] (2node) edge node[left] {$s_2$} (5node);
 \path[->] (3node) edge node[left] {$s_4$} (6node);
 \path[->] (4node) edge node[left] {$s_2$} (6node);
 \path[->] (5node) edge node[left] {$s_3$} (7node);
 \path[->] (5node) edge node[left] {$s_1$} (8node);
 \path[->] (6node) edge node[left] {$s_3$} (9node);
 \path[->] (7node) edge node[left] {$s_2$} (9node);
 \path[->] (7node) edge node[left] {$s_1$} (10node);
 \path[->] (7node) edge node[left] {$s_4$} (11node);
 \path[->] (8node) edge node[right] {$s_3$} (10node);
 \path[->] (9node) edge node[left] {$s_4$} (13node);
 \path[->] (10node) edge node[left] {$s_4$} (14node);
 \path[->] (11node) edge node[left] {$s_3$} (12node);
 \path[->] (11node) edge node[right] {$s_2$} (13node);
 \path[->] (11node) edge node[right] {$s_1$} (14node);
 \path[->] (12node) edge node[left] {$s_2$} (15node);
 \path[->] (12node) edge node[left] {$s_1$} (17node);
 \path[->] (12node) edge node[right] {$s_4$} (18node);
 \path[->] (13node) edge node[right] {$s_3$} (16node);
 \path[->] (14node) edge node[right] {$s_3$} (17node);
 \path[->] (15node) edge node[left] {$s_1$} (19node);
 \path[->] (15node) edge node[left] {$s_3$} (20node);
 \path[->] (15node) edge node[right] {$s_4$} (22node);
 \path[->] (16node) edge node[right] {$s_2$} (20node);
 \path[->] (17node) edge node[right] {$s_2$} (21node);
 \path[->] (17node) edge node[left] {$s_4$} (23node);
 \path[->] (18node) edge node[right] {$s_2$} (22node);
 \path[->] (18node) edge node[right] {$s_1$} (23node);
 \path[->] (19node) edge node[left] {$s_3$} (24node);
 \path[->] (19node) edge node[left] {$s_2$} (25node);
 \path[->] (19node) edge node[right] {$s_4$} (26node);
 \path[->] (20node) edge node[right] {$s_1$} (24node);
 \path[->] (21node) edge node[right] {$s_1$} (25node);
 \path[->] (21node) edge node[left] {$s_4$} (27node);
 \path[->] (22node) edge node[right] {$s_1$} (26node);
 \path[->] (23node) edge node[left] {$s_2$} (27node);
 \path[->] (25node) edge node[left] {$s_4$} (28node);
 \path[->] (26node) edge node[left] {$s_2$} (28node);
 \path[->] (27node) edge node[left] {$s_1$} (28node);
 \path[->] (27node) edge node[left] {$s_3$} (29node);
 \path[->] (28node) edge node[left] {$s_3$} (30node);
 \path[->] (29node) edge node[left] {$s_1$} (30node);
 \path[->] (29node) edge node[left] {$s_4$} (31node);
 \path[->] (30node) edge node[left] {$s_4$} (33node);
 \path[->] (31node) edge node[right] {$s_3$} (32node);
 \path[->] (31node) edge node[right] {$s_1$} (33node);
 \path[->] (32node) edge node[left] {$s_1$} (34node);
 \path[->] (33node) edge node[left] {$s_3$} (34node);
 \path[->] (34node) edge node[left] {$s_2$} (35node);
\end{tikzpicture}
\end{align}

\section{Appendix: Low rank discrete series central characters} \label{s appendix low rank ds}

\subsection{$A_1$}

The root system of type $A_1$ has one simple root denoted by $\alpha$. We assume $\langle \alpha, \alpha \rangle=2$. The only one discrete series central character for $\mathbb H(A_1)$ is $\gamma=\frac{c}{2}\alpha$. We have
\[  ||\gamma||^2=\frac{c^2}{2} .
\]

\subsection{$A_2$}

The root system of type $A_2$ has two simple roots denoted by $\alpha_1$ and $\alpha_2$. We again assume $\langle \alpha_1, \alpha_1\rangle =\langle \alpha_2, \alpha_2\rangle =2$. The only discrete series central character for $\mathbb H(A_2)$ is $c(\alpha_1+\alpha_2)$, with square of the length equal to $2c^2$.



\subsection{$I_2(5)$}

The root system of type $I_2(5)$ has two simple roots denoted by $\alpha_1$ and $\alpha_2$. We again assume $\langle \alpha_1, \alpha_1\rangle=\langle \alpha_2, \alpha_2\rangle =2$. There are two discrete series characters: 
\begin{enumerate}
\item $2(b+1)c(\alpha_1+\alpha_2);$
\item $c(\alpha_1+\alpha_2).$
\end{enumerate}
with repsective square lengths equal to $\frac{c^2}{2b^2}$ and $\frac{c^2}{2a^2}$.

\subsection{$A_3$}

The root system of type $A_3$ has three simple roots denoted by $\alpha_1,\alpha_2,\alpha_3$, corresponding to the Dynkin diagram:
\[
\setlength{\unitlength}{0.8cm}
\begin{picture}(12,1.3)
\thicklines
\put(6,0.8){\circle*{0.2}}
\put(7,0.8){\circle*{0.2}}
\put(8,0.8){\circle*{0.2}}
\put(6,0.8){\line(1,0){1}}
\put(7,0.8){\line(1,0){1}}
\put(5.8,0.4){$\alpha_1$}
\put(6.8,0.4){$\alpha_2$}
\put(7.8,0.4){$\alpha_3$}
\end{picture}.
\]
We again assume $\langle \alpha_1, \alpha_1\rangle =\langle \alpha_2, \alpha_2\rangle =\langle \alpha_3, \alpha_3\rangle= 2$. The only discrete series central character for $\mathbb H(A_3)$ is $c(\frac{3}{2}\alpha_1+\alpha_2+\frac{3}{2}\alpha_3)$ with the square of the length equal to $5c^2$.

\subsection{$H_3$}

The root system of type $H_3$ has three simple roots denoted by $\alpha_1, \alpha_2, \alpha_3$, corresponding to the Dynkin diagram:
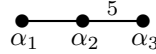
\[
\setlength{\unitlength}{0.8cm}
\begin{picture}(12,1.3)
\thicklines
\put(7.4,0.9){\footnotesize $5$}
\put(6,0.8){\circle*{0.2}}
\put(7,0.8){\circle*{0.2}}
\put(8,0.8){\circle*{0.2}}
\put(6,0.8){\line(1,0){1}}
\put(7,0.8){\line(1,0){1}}
\put(5.8,0.4){$\alpha_1$}
\put(6.8,0.4){$\alpha_2$}
\put(7.8,0.4){$\alpha_3$}
\end{picture}.
\]
We again assume $\langle \alpha_1, \alpha_1 \rangle= \langle \alpha_2, \alpha_2\rangle =\langle \alpha_3, \alpha_3 \rangle =2$. The four discrete series central characters are:
\begin{enumerate}
\item $\frac{c}{2}((10a+5)\alpha_1+(20a+8)\alpha_2+(18a+6)\alpha_3);$
\item $c((2a+2)\alpha_1+(4a+3)\alpha_2+6a\alpha_3);$
\item $\frac{c}{2}((5-2a)\alpha_1+(8-4a)\alpha_2+6a\alpha_3);$
\item $\frac{c}{8+12a}((20a+7)\alpha_1+(40a+12)\alpha_2+(34a+11)\alpha_3)$
\end{enumerate}
with respective squares of lengths $\frac{(48a+19)c^2}{2}$, $8c^2$, $\frac{(43-48a)c^2}{2}$, $\frac{11c^2}{2}$.

\subsection{$A_1\times A_1$, $A_2\times A_1$ and $I_2(5)\times A_1$}

The root systems of types $A_1\times A_1$, $A_2\times A_1$ and $I_2(5)\times A_1$ also appear as a subroot system of $H_4$. The discrete series central characters of those types is a sum of discrete series central characters above. For example, if we consider three simple roots denoted by $\alpha_1,\alpha_2,\alpha_3$, corresponding to the Dynkin diagram:
\[
\setlength{\unitlength}{0.8cm}
\begin{picture}(12,1.3)
\thicklines
\put(6,0.8){\circle*{0.2}}
\put(7,0.8){\circle*{0.2}}
\put(8,0.8){\circle*{0.2}}
\put(6,0.8){\line(1,0){1}}
\put(5.8,0.4){$\alpha_1$}
\put(6.8,0.4){$\alpha_2$}
\put(7.8,0.4){$\alpha_3$}
\end{picture}
\]
The only discrete series central character  for $\mathbb H(A_2\times\mathbb{Z}_2)$ is $c(\alpha_1+\alpha_2+\frac{1}{2}\alpha_3)$.





\section{Appendix: The sets $P(\chi)$ and $Z(\chi)$ for $\chi_{16}$ and $\chi_{17}$}
\subsection{The set $P(w_{17}^*(\chi_{17}))$} \label{ss explicit p wchi17}

We use the notations in Section \ref{ss minimal induce chi 17}. Let $w_{17}^*$ be the element in $W$ such that $w_{17}^*(\chi_{17})=\chi'+\varpi'$ and $R(w_{17}^*)\cap Z(\chi_{17})=\emptyset$. In SageMath code, it corresponds to $W[3101]$. The set $P(w_{17}^*(\chi_{17}))$ contains the following vectors:
\begin{enumerate}
\item $w_{17}^*(\beta_{17}^1)=-4a\alpha_1-(6a+1)\alpha_2-(8a+2)\alpha_3-(6a+2)\alpha_4=\sqrt{2}(0,0,0,-1)$
\item $w_{17}^*(\beta_{17}^2)=-2a\alpha_1-2a\alpha_2-(2a+1)\alpha_3-2a\alpha_4=\sqrt{2}(a,b,0,-1/2)$
\item $w_{17}^*(\beta^3_{17})=-(2a+1)\alpha_1-(4a+1)\alpha_2-(6a+1)\alpha_3-(4a+2)\alpha_4=\sqrt{2}(-b,1/2,0,-a)$
\item $w^*_{17}(\beta_{17}^4)=-\alpha_1=\sqrt{2}(1/2,a,0,-b)$
\item $w_{17}^*(\beta_{17}^5)=-2a\alpha_1-(4a+1)\alpha_2-(6a+2)\alpha_3-(4a+2)\alpha_4=\sqrt{2}(0,-a,-b,-1/2)$
\item $w_{17}^*(\beta_{17}^6)=-\alpha_3=\sqrt{2}(a,-1/2,-b,0)$
\item $w_{17}^*(\beta_{17}^7)=-2a\alpha_2-4a\alpha_3-(2a+1)\alpha_4=\sqrt{2}(-b,-a,-1/2,0)$
\item $w_{17}^*(\beta_{17}^8)=2a\alpha_1+(2a+1)\alpha_2+(2a+1)\alpha_3+(2a+1)\alpha_4=\sqrt{2}(1/2,-1/2,-1/2,1/2)$
\item $w_{17}^*(\beta_{17}^9)=-\alpha_4=\sqrt{2}(-a,1/2,-b,0)$
\item $w^*_{17}(\beta_{17}^{10})=2a\alpha_1+(4a+1)\alpha_2+(6a+1)\alpha_3+(4a+1)\alpha_4=\sqrt{2}(0,a,-b,1/2)$
\item $w_{17}^*(\beta_{17}^{11})=\alpha_1+\alpha_2+\alpha_3=\sqrt{2}(-a,0,-1/2,b)$
\item $w_{17}^*(\beta_{17}^{12})=(2a+1)\alpha_1+(4a+2)\alpha_2+(6a+2)\alpha_3+(4a+2)\alpha_4=\sqrt{2}(0,b,-1/2,a)$
\end{enumerate}

\subsection{The set $Z(w_{17}^*(\chi_{17}))$} \label{ss stabilizer 17}
The set $Z(w_{17}^*(\chi_{17}))$ contains:
\begin{enumerate}
\item $w_{17}^*\alpha_1=\alpha_1+(2a+1)\alpha_2+(4a+1)\alpha_3+(2a+1)\alpha_4=\sqrt{2}(-1/2,a,0,b)$
\item $w_{17}^*\alpha_2=2a\alpha_1+2a\alpha_2+2a\alpha_3+2a\alpha_4=\sqrt{2}(0,-a,-b,1/2)$
\item $w_{17}^*\alpha_4=2a\alpha_1+(4a+1)\alpha_2+(6a+1)\alpha_3+(4a+2)\alpha_4=\sqrt{2}(a,b,0,1/2)$
\item $w_{17}^*(\alpha_1+\alpha_2)=(2a+1)\alpha_1+(4a+1)\alpha_2+(6a+1)\alpha_3+(4a+1)\alpha_4=\sqrt{2}(-1/2,0,-b,a)$
\end{enumerate}

\subsection{The stabilizer group $\mathrm{Stab}(w_{17}^*\chi_{17})$} \label{ss tabilizer 17}

$$\mathrm{Stab}(w_{17}^*\chi_{17})=\langle w_{17}^*s_1(w_{17}^*)^{-1},w_{17}^*s_2(w_{17}^*)^{-1},w_{17}^*s_4(w_{17}^*)^{-1} \rangle$$
The two non-trivial elements in  Lemma \ref{lem stabilizer group 17} are
$w_{17}^*s_2s_1(w_{17}^*)^{-1}$ and $w_{17}^*s_1s_2s_1s_4(w_{17}^*)^{-1}$ respectively.

\subsection{The set $P(w_{16}^*\chi_{16})$}

We use the notations in Section \ref{ss skew 16 ds}.  In SageMath code, it corresponds to $W[1040]$. The set $P(w_{16}^*(\chi_{16}))$ contains the following vectors:
\begin{enumerate}
\item $w_{16}^*(\beta_{16}^1)=-(2a+1)\alpha_1-(4a+1)\alpha_2-(6a+2)\alpha_3-(6a+1)\alpha_4=\sqrt{2}(0,b,-1/2,-a)$
\item $w_{16}^*(\beta_{16}^2)=-2a\alpha_3-2a\alpha_4=\sqrt{2}(0,0,-1,0)$
\item $w_{16}^*(\beta^3_{16})=-2a\alpha_1-4a\alpha_2-(4a+1)\alpha_3-(4a+1)\alpha_4=\sqrt{2}(-a,b,0,1/2)$
\item $w^*_{16}(\beta_{16}^4)=\alpha_1+\alpha_2+\alpha_3=\sqrt{2}(-a,0,-1/2,b)$
\item $w_{16}^*(\beta_{16}^5)=-2a\alpha_1-(2a+1)\alpha_2-(4a+1)\alpha_3-(2a+1)\alpha_4=\sqrt{2}(a,-b,0,-1/2)$
\item $w_{16}^*(\beta_{16}^6)=-\alpha_3=\sqrt{2}(a,-1/2,-b,0)$
\item $w_{16}^*(\beta_{16}^7)=-\alpha_2=\sqrt{2}(-1/2,-b,a,0)$
\item $w_{16}^*(\beta_{16}^8)=2a\alpha_1+2a\alpha_2+4a\alpha_3+(2a+1)\alpha_4=\sqrt{2}(-1/2,-1/2,1/2,1/2)$ 
\end{enumerate}

\subsection{The set $Z(w_{16}^*(\chi_{16}))$} \label{ss stabilizer chi16 others}
The set $Z(w_{16}^*(\chi_{16}))$ contains:
\begin{enumerate}
\item $w_{16}^*\alpha_1=\alpha_1+\alpha_2+(2a+1)\alpha_3+2a\alpha_4=\sqrt{2}(-a,0,1/2,b)$
\item $w_{16}^*\alpha_4=(2a+1)\alpha_1+(4a+1)\alpha_2+(4a+2)\alpha_3+(4a+1)\alpha_4=\sqrt{2}(0,-b,-1/2,a)$
\end{enumerate}




\subsection{Weight structure of $U''$}
The 4 weights of $\widetilde{M}\otimes\mathbb{C}_{-b\varpi_{1}}$ are:

$$w_{16}^*\chi_{16}=-\frac{1}{2}(4a\alpha_1+(8a+3)\alpha_2+(12a+5)\alpha_3+(12a+2)\alpha_4)$$
$$s_4w_{16}^*\chi_{16}=-\frac{1}{2}(4a\alpha_1+(8a+3)\alpha_2+(12a+5)\alpha_3+(10a+4)\alpha_4)$$
$$s_3s_4w_{16}^*\chi_{16}=-\frac{1}{2}(4a\alpha_1+(8a+3)\alpha_2+(14a+3)\alpha_3+(10a+4)\alpha_4)$$
$$s_2s_3s_4w_{16}^*\chi_{16}=-\frac{1}{2}(4a\alpha_1+10a\alpha_2+(14a+3)\alpha_3+(10a+4)\alpha_4)$$

\subsection{The stabilizer subgroup $\mathrm{Stab}(w_{16}^*\chi_{16})$} \label{ss stabilizer 16}

The stabilizer group $\mathrm{Stab}(w_{16}^*\chi_{16})$ is generated by two reflections $w_{16}^*s_1(w_{16}^*)^{-1}$ and $w_{16}^*s_4(w_{16}^*)^{-1}$, and is isomorphic to $\mathbb Z/2\mathbb Z\times \mathbb Z/2\mathbb Z$. The reduced decompositions of the non-trivial elements are:
\begin{itemize}
    \item $w_{(2)}:=w_{16}^*s_1(w_{16}^*)^{-1}=s_3s_4s_1s_2s_3s_4s_2s_3s_1$
    \item $w_{(3)}:=w_{16}^*s_4(w_{16}^*)^{-1}=s_1s_2s_3s_4s_3s_4s_1s_2s_3s_4s_1s_2s_3s_4s_2s_3s_4s_1s_2s_3s_4s_3s_1s_2s_1$
    \item $w_{(4)}:=w_{16}^*s_1s_4(w_{16}^*)^{-1}=s_1s_2s_3s_4s_3s_4s_2s_3s_4s_1s_2s_3s_4s_2s_3s_4s_2s_3s_4s_1s_2s_3s_4s_3s_2s_1$
    \end{itemize}
The last element is the non-trivial element in Proposition \ref{prop stabilizer chi16}.

 \end{document}